\documentclass[conf]{new-aiaa}
\usepackage[utf8]{inputenc}

\usepackage{graphicx}
\usepackage{amsmath}
\usepackage[version=4]{mhchem}
\usepackage{siunitx}
\usepackage{longtable,tabularx}

\usepackage{amsthm}
\usepackage{algorithm}
\usepackage{algpseudocode}

\usepackage{xcolor}
\usepackage[most]{tcolorbox}
\title{Constrained Lyapunov Stabilization based on Gauss Variational Equations: From Spacecraft Orbital Transfers to Rendezvous}

\author{Ilya Kolmanovsky\footnote{Professor, Department of Aerospace Engineering, AIAA Associate Fellow.}}
\affil{University of Michigan, 
	3038 Francois-Xavier Bagnoud Building, 1320 Beal Avenue, Ann Arbor, MI 48109-2140; {\tt ilya@umich.edu}} 

\author{Emanuele Garone\footnote{Professor, Service d’Automatique et Analyse des Systemes, AIAA Senior Member.}}
\affil{Université Libre de Bruxelles, Brussels, Belgium; {\tt emanuele.garone@ulb.be}}

\author{Grant Touchette\footnote{Graduate Student, Department of Aerospace Engineering.}
}
\affil{University of Michigan, 
	Francois-Xavier Bagnoud Building, 1320 Beal Avenue, Ann Arbor, MI 48109-2140; {\tt togrant@umich.edu}}

\theoremstyle{remark}

\begin{document}

\maketitle

\begin{abstract}
    Lyapunov feedback laws can be constructed for performing orbital transfer maneuvers based on Gauss Variations Equations (GVEs), while satisfying the specified state and control constraints.  These state and control constraints are enforced using barrier functions and saturation, respectively, while the reference governor, employed as a convergence governor, is utilized to avoid getting stuck at spurious equilibria that may be created by barrier functions.
    In this article, these Lyapunov feedback laws are extended to rendezvous maneuvers where not only five orbital elements are matched to the prescribed target values but also the true anomaly matches its (time-dependent) target value.  The modification involves altering the commanded semi-major axis with an outer-loop feedback law, also designed using Lyapunov techniques.  We illustrate the resulting safe closed-loop rendezvous maneuvers in simulations for the conventional thrust-based propulsion and for Lorenz force-based propulsion.  In the latter case, only the current through the tether is controlled subject to the current limits to accomplish prescribed orbital transfer and rendezvous maneuvers.  
\end{abstract}

\begin{center}
\begin{tcolorbox}[
    colback=gray!15,
    colframe=gray!55!black,
    width=0.94\textwidth,
    arc=2mm,
    boxrule=0.6pt,
    left=10pt,
    right=10pt,
    top=8pt,
    bottom=8pt
]
\textbf{This document is an author's original manuscript (preprint) of the following conference paper:}

\medskip

\noindent
I. Kolmanovsky, E. Garone, and G. Touchette,
``Constrained Lyapunov Stabilization Based on Gauss Variational Equations:
From Spacecraft Orbital Transfers to Rendezvous,''
\emph{2026 AIAA SciTech Forum}, Orlando, Florida,
Paper AIAA 2026-0636, Session: Low-Thrust Trajectories,
published online Jan. 8, 2026.
DOI: \href{https://doi.org/10.2514/6.2026-0636}{10.2514/6.2026-0636}.
\end{tcolorbox}
\end{center}

\section{Nomenclature}

{\renewcommand\arraystretch{1.0}
\noindent\begin{longtable*}{@{}l @{\quad=\quad} l@{}}
$a$  & semi-major axis \\
$e$ &    eccentricity \\
$i$& inclination \\
$\omega$ & argument of periapsis \\
$\Omega$ & RAAN \\
$\theta$   & true anomaly \\
$M$ & mean anomaly \\
$\psi$ & eccentric anomaly \\
$\vec{F}$ & force vector \\
$m$ & spacecraft mass \\
$X$ & state vector \\
$\tilde{X}$ & extended state vector\\ 
$(\cdot)_{\tt des}$ & target value \\
$a_{\tt alt}$ & altered semi-major axis command by outer loop controller \\
$U$ & control input vector \\
$\mathcal{U}$ & set defining control constraints\\
$F(X,\theta), G(X,\theta)$   & model functions \\
$c_i$  & constraints \\
$V_0$, $V$, $\bar{V}$, $\tilde{V}$ & Lyapunov functions \\
$B_i$ & Barrier functions \\
$q_i$ & barrier function weights \\
$\epsilon_i$ & barrier function parameters \\
$\varepsilon_0$ & terminal condition parameter used by refeerence governor \\
$P$, $\bar{P}$ & Lyapunov function weight matrices \\
$\mathcal{B}_\epsilon$ & ball of radius $\epsilon$ \\
$S$, $T$, $W$ & perturbation force (thrust of Lorentz force)
components \\
$B_x$, $B_y$, $B_z$ & components of magnetic field \\
$\vec{l}$ & vector along tether direction with magnitude equal to tether length \\
$I$ & current \\
$t$  & time  \\
$\hat{e}_r$, $\hat{e}_\theta$, $\hat{e}_h$ & unit vectors of STW (LVLH) frame \\
$X_{\tt des,alt}$, $\kappa$, $\varepsilon_0$ & modified by reference governor target state, and reference governor parameters
\end{longtable*}}

\section{Introduction}

\lettrine{L}yapunov feedback  laws for orbital transfers can be defined \cite{gurfil2016celestial} based on quadratic Lyapunov functions and Gauss Variational Equations (GVEs).  These feedback laws can steer the spacecraft to an orbit prescribed by the target semi-major axis $a$ [km], eccentricity $e$, inclination $i$ [rad], right ascension of asending node (RAAN) $\Omega$ [rad], and argument of periapsis $\omega$ [rad].  Their potential advantage over open-loop optimized trajectories is robustness to uncertainties and disturbances such as thrust errors or lunisolar perturbations.
  
More recently \cite{garone2024constrained}, such feedback laws have also been proven to be stabilizing under arbitrary tight and  even  one-sided constraints on the radial and out-of-the orbital plane components of the 
thrust induced accelerations.
Interestingly, by invoking LaSalle’s invariance principle, one can establish closed-loop asymptotic stability without imposing additional persistence-of-excitation conditions \cite{garone2024constrained}. However, analytically demonstrating exponential convergence rates still appears to require such conditions \cite{semeraro2024reference}.

Furthermore, it has been shown that state constraints on the minimum radius of periapsis and minimum eccentricity can be handled by the augmentation of the Lyapunov function by barrier functions~\cite{semeraro2023constrained} or by using reference governors~\cite{semeraro2024reference,Kolmanovsky2024CommandGovernors}.  Since convergence properties in \cite{garone2024constrained} hold for any positive definite quadratic Lyapunov function and associated control law, onboard optimization of the weights in the Lyapunov function can be employed to speed-up the closed-loop maneuvers; the speed-up is achieved by aligning the sublevel sets of the Lyapunov function with the constraints \cite{Kolmanovsky2024CommandGovernors}.

In  this  paper we demonstrate that the above feedback laws  can be interconnected with an outer loop feedback law that adjusts the commanded semi-major axis to additionally match the true anomaly $\theta$ [rad] to the (time-varying) target true anomaly, $\theta_{\tt des}$. The combined control law is thus able to perform the rendezvous maneuver on the target orbit based on GVEs.
The proposed approach is applicable to large maneuvers and handles state and control constraints through the use of saturation and barrier functions.  Notably, the use of barrier functions to enforce state constraints is important in the setting where a commanded semi-major axis to the inner-loop controller is adjusted by an outer-loop controller. Otherwise, the orbit manipulated by the outer-loop controller could lead to unacceptable altitude and collisions with the primary body.  As in \cite{semeraro2023constrained}, we employ the reference governor as a convergence governor to ensure that the trajectory avoids spurious equilibria that could potentially be created by the barrier functions.

In addition to simulations for the case of the conventional thrust we show that our approach can also be applied to rendezvous maneuvers which exploit Lorenz force actuation \cite{cosmo1997tethers,levin2007dynamic}.  In this case, the control input is the current passed through the tether.  The tether is assumed to remain stretched along the radial direction throughout the maneuver.  In this case, constrained rendezvous maneuvers can be accomplished by modulating the single control input, the current through the tether.

\section{Gauss Variational Equations}

The evolution of the classical osculating orbital elements  is described by the Gauss-Euler Variational Equations (GVEs) \cite{gurfil2016celestial}: 
\begin{align}\label{equ:gve}
	\frac{da}{dt} &= \frac{2a^2}{\sqrt{\mu p}}e\sin{\theta}S+\frac{2a^2}{\sqrt{\mu p}}\frac{p}{r}T, \nonumber \\
	\frac{de}{dt} &= \frac{p\sin\theta}{\sqrt{\mu p}}S + \frac{p(\cos\psi+\cos\theta)}{\sqrt{\mu p}}T, \nonumber\\
	\frac{di}{dt} &= \frac{r}{\sqrt{\mu p}}\cos{(\theta+\omega)}W, \\
	\frac{d\Omega}{dt} &= \frac{r}{\sqrt{\mu p}}\frac{\sin{(\theta+\omega)}}{\sin i}W,
	\nonumber \\
	\frac{d\omega}{dt} &= -\frac{p\cos\theta}{e\sqrt{\mu p}}S+\frac{(r+p)\sin\theta}{e\sqrt{\mu p}}T-\frac{r\sin(\theta+\omega)\cot i}{\sqrt{\mu p}}W,\nonumber \\
	\frac{d\theta}{dt} &= \frac{\sqrt{\mu p}}{r^2}+\frac{p\cos\theta}{e}\frac{S}{\sqrt{\mu p}}-\frac{p+r}{e}\cos\theta\frac{T}{\sqrt{\mu p}}. \nonumber
\end{align}
In (\ref{equ:gve}), $\mu$ is the gravitational parameter of the primary, $r=\frac{p}{1+e\cos\theta}$ the orbit radius (distance from the gravity center to the spacecraft center of mass), $p=a(1-e^2)=h^2/\mu$ is the orbit parameter, $h$ is the magnitude of the angular momentum, and $\psi = \arccos\left(\frac{1}{e}-\frac{r}{ae}\right)$ is the eccentric anomaly.  The variables  
$S$, $T$ and $W$ represent thrust/perturbation induced relative acceleration [km/s$^2$] components defined so that the perturbation force $\vec{F}$ can be represented as
\begin{align*}
	\frac{\vec{F}}{m} = S\hat{e}_r + T\hat{e}_\theta + W\hat{e}_h,
\end{align*}
where $m$ denotes the mass of the spacecraft while $\hat{e}_r$, $\hat{e}_\theta$ and $\hat{e}_h$ designate unit vectors satisfying
$$\hat{e}_r=\frac{\vec{r}}{r},~\hat{e}_h=\frac{\vec{h}}{h},~ \vec{h} =\vec{r} \times \vec{v},~
\mbox{and}~ \hat{e}_\theta=\hat{e}_h \times \hat{e}_r.$$
Here
$\vec{r}$ denotes the spacecraft position vector in an inertial frame centered at the primary body, $r=\|\vec{r}\|$,  $\vec{v}$ is the spacecraft velocity vector, and  $\times$ denotes the vector product.

With
\begin{align*}
	X = [a \ \ e \ \ i \ \ \Omega \ \ \omega]^{\sf T}, \ \ U = [S \ \ T \ \ W]^{\sf T},
\end{align*}
the first five of GVEs (\ref{equ:gve}) take the form of a drift-free nonlinear system which is affine in control,
\begin{equation} \label{equ:main1}
	\dot X(t) = G\big(X(t),\theta(t)\big)U(t),
\end{equation}
where $G(X(t),\theta(t)) \in \mathbb{R}^{5 \times 3}$.
Note that the true anomaly $\theta(t)$ is not included into the state vector $X(t)$ and is treated as a time-varying parameter the evolution of which is determined by the sixth equation
in (\ref{equ:gve}). 

The control constraints are expressed as 
\[U(t) \in \mathcal{U} \mbox{ for all $t \geq 0$},\]
where $\mathcal{U}$ is a prescribed set.

\section{Lyapunov Control Design for Orbital Transfers}

Let $P=P^{\sf T} \succ 0$ be {\em any} $5 \times 5$ positive-definite weight matrix and select a Lyapunov function candidate,
\begin{equation}\label{equ:main2}
	V_0(X, X_{\tt des})=\frac{1}{2}\left(X-X_{\tt des}\right)^{\sf T}P\left(X-X_{\tt des}\right).
\end{equation}
Consider a minimization of  $V_0(X(t+\Delta t), X_{\tt des})$ for small $\Delta t>0$ and of the control effort involved, i.e., of
\begin{gather}
V_0(X(t+\Delta t), X_{\tt des}) + \frac{1}{2}\int_t^{t+\Delta t} U^{\sf T}(\tau) U(\tau) d \tau \nonumber \\
\approx V_0(X(t), X_{\tt des})+\frac{\partial V_0}{\partial X} [X(t), X_{\tt des}]G\big(X(t),\theta(t)\big)U(t) \Delta t+
\frac{1}{2} U^{\sf T}(t) U(t) \Delta t,
\label{equ:appro}
\end{gather}
with respect to $U(t) \in \mathcal{U}$. 
A single iteration of a primal projected gradient algorithm applied to minimization of the function on the right hand side of (\ref{equ:appro}) yields a Lyapunov feedback law,  
	\begin{equation}\label{equ:u_nocon}
			U(t) = \mathrm{Proj}_\mathcal{U} \left[U_{\tt nom}(t) \right],\quad  U_{\tt nom}=-G^\top(X(t), \theta(t)) \left(\frac{\partial V_0}{\partial X} [X(t), X_{\tt des}]\right)^{\sf T}
            =-G^\top(X(t), \theta(t)) P (X - X_{\text{des}}).
			\end{equation}
 The perspective of obtaining (\ref{equ:u_nocon}) by dynamic optimization of (\ref{equ:appro}) is useful in tuning \cite{kolmanovsky2002speed}, e.g., faster convergence of a particular state can be induced by increasing the corresponding weight (i.e., the corresponding diagonal element of $P \succ 0$, assuming $P$ is a diagonal matrix).

The closed-loop stability and convergence properties, i.e., \[
X(t) \to X_{\text{des}} \quad \text{as } t \to \infty \text{ for constant } X_{\text{des}},
\] hold			
under the following assumptions \cite{garone2024constrained}: \\
			\begin{enumerate}
				\item $0 \in \mathrm{int} \, \mathcal{U}, \; \mathcal{U}$ is closed and convex;\\[0.2em]
				\item $\mathrm{Proj}$ is  either the Euclidean or the radial projection; \\[0.2em]
				\item  The closed-loop trajectory remains within the model validity range, which can be ensured, e.g., by the condition
                $\{X: V_0(X,X_{\tt des}) \leq V_0(X(0),X_{\tt des})\} \subset \mathcal{R}$,
				$\mathcal{R} = \{X : a > 0, \; 0 < e < 1, \; i \neq 0 \}$. \\
			\end{enumerate}

Note that the condition $0 \in \mathrm{int} \, \mathcal{U}$  can be replaced \cite{garone2024constrained} by a weaker condition: 
			\[
			\mathcal{U} \supseteq \{[S, T, W]^\top \mid S \geq 0, W \geq 0\} \cap \mathcal{B}_\varepsilon,~\mbox{for some $\varepsilon>0$},
			\]
			or by a similar condition with $S \leq 0$, $W \leq 0$ instead of $S\geq 0$, $W \geq 0$.
Thus the Lyapunov feedback law is stabilizing under arbitrary tight control constraints that on $S$ and $W$ can be even one-sided.  
Assuming the spacecraft needs to be oriented to generate net thrust force in the prescribed direction, restricting thrust directions with one-sided constraints could be advantageous in terms of simplifying the attitude control.  Other advantages of adhering to one-sided thrust constraints and following more complex trajectories could include reducing the possibility of inferring the maneuver-intent (in an adverserial setting).

The above Lyapunov feedback law is related to the classical $Q$-laws \cite{Petropoulos2004QLaw,Petropoulos2005Refined}, and, as shown in \cite{garone2024constrained},  it also possesses rigorous (mathematically provable) closed-loop asymptotic stability guarantees, even under stringent control constraints satisfying the above assumptions.

\section{Handling State Constraints with Barrier Functions}

\subsection{Barrier functions}
Following \cite{semeraro2023constrained,kolmanovsky2002speed,
dongare2022integrated},
barrier functions are used to enforce state constraints.  To accommodate $n_c \geq 1$ constraints of the form,
\begin{equation}\label{equ:cons}
c_i(X) \geq 0,~i=1,\cdots,n_c,
\end{equation}
we augment the Lyapunov function with the barrier terms:
$$V(X,X_{\tt des})=V_0(X,X_{\tt des})+\sum_{i=1}^{n_c}B_i(X),$$
where
\begin{equation}\label{equ:main2-1}
	B_i(X) = \left\{
	\begin{array}{cc}
		\frac{1}{2} q_i(c_i(X) - \epsilon_i)^2 
		\hfill & 
		\text{if } c_i((X)<\epsilon_i, \\
		0 \hfill & \text{otherwise},
	\end{array} \right. 
\end{equation}
and where $q_i>0$ are weight parameters penalizing the respective constraint violations, while $\epsilon_i \geq 0$ are safety margins introduced so that the barrier functions become non-zero prior to constraints becoming violated.

The nominal feedback law (\ref{equ:u_nocon}) is modified to
\begin{align} \label{equ:controllaw1}
	U_{\tt nom}  = -G^\top(X(t), \theta(t)) \left(\frac{\partial V}{\partial X} [X(t), X_{\tt des}]\right)^{\sf T}= -\left[(X-X_{\tt des})^{\sf T} PG(X,\theta) + C^{\sf T}(X)G(X,\theta)\right]^{\sf T},  
\end{align}
which is saturated to enforce the control constraints as
\begin{equation}\label{equ:control_sat}
U=\mathrm{Proj}_\mathcal{U} \left[ U_{\tt nom} \right], \end{equation}
where
\begin{equation}
	C(X) = \sum_{i=1}^{n_c} \left(\dfrac{\partial B_i}{\partial X} (X) \right)^{\sf T}.
\end{equation}
With this feedback law,
$$\frac{d}{dt} V(X(t),X_{\tt des}) \leq 0,$$ 
along the closed-loop system trajectories implying that the sublevel sets of $V$ are positively invariant. 


Consider now a trajectory $X(t)$ emanating from a state, $X(t_0)$, at the time instant, $t_0$, that does not violate the constraints plus safety margins at the time instant $t_0$, so that $B_i(X(t_0)) = 0$, $i=1,\cdots,n_c$.
Note that $V(X(t_0),X_{\tt des})=V_0(X(t_0),X_{\tt des})$.
Since $V(X(t)),X_{\tt des})$ is a non-increasing function of time $t$ and $B_i(X(t)) \geq 0$, $i=1,\cdots,n_c$,
it follows that
\begin{align*}
	\max_{i=1,\cdots,n_c} &\{B_i(X(t))\} \leq \sum_{i=1}^{n_c}
    B_i(X(t)) \leq V(X(t),X_{\tt des}) \leq V(X(t_0),X_{\tt des}) \leq  V_0(X(t_0),X_{\tt des}).
\end{align*}
If $q_i$ satisfy
\begin{equation}\label{equ:logic}
	q_i \geq \frac{2}{ \epsilon_i^2}V_0(X(t_0),X_{\tt des}), \quad i=1,\cdots,n_c,
\end{equation}
then $B_i(X(t)) \leq \frac{1}{2} q_i \epsilon_i^2$,
for all $t \geq t_0$ and $i=1,\cdots,n_c$; hence,
 the constraints (\ref{equ:cons})  are enforced along the closed-loop trajectory.
  
The conditions (\ref{equ:logic})  can be used to reset $q_i$ 
along the trajectory.  Specifically if $B_i(X(t_k))=0$ for all $i=1,\cdots,n_c$ at any given time instant $t_k$,
then resetting
\begin{equation}\label{equ:reset}
	q_i(t_k) = \frac{2}{ \epsilon_i^2}V_0(X(t_k),X_{\tt des}),~i=1,\cdots,n_c,
\end{equation}
does not change the value of the Lyapunov function, $V(X(t_k),X_{\tt des})$, and still ensures that the constraints are enforced.  As $V_0(X(t),X_{\tt des})$ is typically decreasing with time $t$, the periodic reset (\ref{equ:reset}) lowers the values of $q_i$'s  and facilitates recovering the unconstrained controller performance.  Note also that the values of $q_i$, $i=1,\cdots,n_c$, will need to be reset whenever $X_{\tt des}$ changes.

\subsection{Constraints}

In the sequel, we consider spacecraft orbital maneuvers under the following five constraints.  

The first constraint has the form,
\begin{equation}\label{equ:maincon1}
	c_1(X) = r_{\tt p} - r_{\tt min} \geq 0,\quad r_{\tt p}=a(1-e),
\end{equation}
and ensures that the radius of the periapsis of the spacecraft orbit, $r_{\tt p}$,
is larger than the specified value, $r_{\tt min}$.  The constraint (\ref{equ:maincon1})
protects not only against the distance to the primary $r$ falling below $r_{\tt min}$ at any given time instant (and thus avoiding being too close to/colliding with the primary) but also that $r$ will stay above $r_{\tt min}$ even if there is a thruster failure and thrust becomes zero.

The second constraint, 
\begin{equation}\label{equ:maincon3}
	c_2(X) = e-e_{\tt min} \geq  0,
\end{equation}
ensures that the eccentricity of the spacecraft orbit is maintained above a specified minimum value, $e_{\tt min} \geq 0$.
This constraint ensures that the eccentricity does not approach zero too closely where GVEs have a singularity therefore preserving the validity of the model.  Similarly, we impose the constraint on the the maximum eccentricity,
\begin{equation}\label{equ:maincon3_1}
	c_3(X) = e_{\tt max}-e \geq  0,
\end{equation}
where $e_{\tt max}<1$ so that the osculating orbit remains elliptic.
The fourth constraint is imposed on the minimum inclination to preserve the validity of GVE-based model:
\begin{equation}\label{equ:maincon3_2}
	c_4(X) = i-i_{\tt min} \geq  0.
\end{equation}

\subsection{Gradients of barrier functions}

The gradients of the barrier functions corresponding to these constraints are given by explicit expressions:
\[
 \left(\dfrac{\partial B_1}{\partial X}\right)^{\sf T} =
\begin{cases}
\begin{bmatrix}
q_1(1-e)\big(a(1-e)-r_{\min}-\epsilon_1\big)\\[4pt]
-q_1a\big(a(1-e)-r_{\min}-\epsilon_1\big)\\[4pt]
0\\[2pt]0\\[2pt]0
\end{bmatrix}, & \text{if } a(1-e)<r_{\min}+\epsilon_1,\\[14pt]
\mathbf{0}, & \text{otherwise},
\end{cases}
\]

\[
\left(\dfrac{\partial B_2}{\partial X}\right)^{\sf T} = 
\begin{cases}
\begin{bmatrix}
0\\[2pt]
q_2\big(e-e_{\min}-\epsilon_2\big)\\[2pt]
0\\[2pt]0\\[2pt]0
\end{bmatrix}, & \text{if } e<e_{\min}+\epsilon_2,\\[14pt]
\mathbf{0}, & \text{otherwise},
\end{cases}
\qquad 
\left(\dfrac{\partial B_3}{\partial X}\right)^{\sf T} = 
\begin{cases}
\begin{bmatrix}
0\\[2pt]
q_3\big(e-e_{\max}+\epsilon_3\big)\\[2pt]
0\\[2pt]0\\[2pt]0
\end{bmatrix}, & \text{if } e>e_{\max}-\epsilon_3,\\[14pt]
\mathbf{0}, & \text{otherwise},
\end{cases}
\]

\[
\left(\dfrac{\partial B_4}{\partial X}\right)^{\sf T}  = 
\begin{cases}
\begin{bmatrix}
0\\[2pt]
0\\[2pt]
q_4\big(i-i_{\min}-\epsilon_4\big)\\[2pt]
0\\[2pt]0
\end{bmatrix}, & \text{if } i<i_{\min}+\epsilon_4,\\[14pt]
\mathbf{0}, & \text{otherwise},
\end{cases}
\]
where ${\bf 0}$ is a $5 \times 1$ zero vector.

\section{Reference Governor as Convergence Governor}

As  noted in \cite{semeraro2023constrained}, formal convergence guaranties when the Lyapunov feedback law is modified with the barrier terms are lacking, and hypothetically the closed-loop trajectory could ``get stuck'' at spurious equilibria that could be created by the use of the barrier functions. 

To rigorously ensure convergence, in \cite{semeraro2023constrained} we proposed the use of the reference governor~\cite{garone2017reference} as a convergence governor.  This convergence governor replaces the reference command, $X_{\tt des}$, with a modified reference command, $X_{\tt des, alt}$, which is as close as possible to $X_{\tt des}$ subject to the condition that the predicted closed-loop trajectory enters the constraint-admissible invariant terminal set centered around $X_{\tt des, alt}$ in which the barrier function terms remain inactive.

Towards this end, we denote the control law in (\ref{equ:control_sat})
as $U(X,X_{\tt des},q,\theta)$ and we let the predicted state at the time instant $t_k+\tau$ 
with this control law applied and for the given $X(t_k)$ and $\theta(t_k)$ at the time instant $t_k$ be  denoted by $\hat{X}(t_k+\tau,X(t_k),X_{\tt des},q,\theta(t_k))$.
In this paper, we define the convergence governor so that it only requires a single trajectory prediction of $\hat{X}(t_k+\tau,X(t_k),X_{\tt des},q,\theta(t_k))$ at each $t_k$, which can be performed through the onboard simulations.

Consider the given initial state, $X(t_0)$, at the initial time instant $t_0$, and the prescribed target state, $X_{\tt des}$. Assume that there exists $\varepsilon_0>0$ such that 
$$V_0(X,X(t_0)\alpha+(1-\alpha)X_{\tt des}) \leq \varepsilon_0~\Rightarrow~B_i(X) \leq 0,$$
for all $i=1,\cdots,n_c$ and all $0 \leq \alpha \leq 1$. Note that due to this assumption the line segment connecting the initial state, $X(t_0)$, with the target state, $X_{\tt des}$, strictly satisfies the constraints.

At each time instant, $t_k$, at which $B(X(t_k))=0$, the convergence governor computes a candidate target state, $\tilde{X}_{\tt des, alt}$, as
\begin{equation}\label{equ:rg1}
    \tilde{X}_{\tt des, alt}(t_k)=\kappa_k X_{\tt des}+(1-\kappa_k)X_{\tt des,alt}(t_{k-1}),
\end{equation}
and barrier function weight vector, $\tilde{q}(t_k)$, with components,
$$\tilde{q}_i(t_k)=\dfrac{2}{\epsilon_i^2}{V_0(X(t_k),\tilde{X}_{\tt des, alt}(t_k)}),~i=1,\cdots,n_c.$$
If the {\em target acceptance condition,} $$V_0(\hat{X}(t_k+\tau,X(t_k),\tilde{X}_{\tt des,alt},\tilde{q},\theta(t_k)),\tilde{X}_{\tt des, alt}(t_k))\leq \varepsilon_0,$$
is satisfied for some time instant $\tau \in \mathcal{T}$, where  $\mathcal{T}$ is the time interval over which the trajectory is predicted, we set ${X}_{\tt des, alt}(t_k)=
\tilde{X}_{\tt des, alt}(t_k)$,
$q(t_k)=\tilde{q}(t_k)$.    
 In (\ref{equ:rg1}),  the parameter $\kappa_k$, $0 \leq \kappa_k \leq 1$, is chosen so that
$$\kappa_k=\dfrac{1}{k-k_{\tt last}},$$
where $t_{k_{\tt last}}$ is the last time instant at which the target acceptance condition was satisfied.  If at a time instant $t_k$ the target acceptance condition is not satisfied then ${X}_{\tt des, alt}(t_k)=
{X}_{\tt des, alt}(t_{k-1})$. The control law
$U(X,X_{\tt des,alt}(t_k),q(t_k),\theta)$ is applied between the time instants $t_k$ and $t_{k+1}$ to the spacecraft.  



\section{Rendezvous Using Conventional Thrust Actuation}

In order to perform a rendezvous maneuver, in addition to matching five orbital elements to their respective target values, the phasing has to be also matched.  

\subsection{Incorporating $M$ into the state vector approach}

An approach pursued initially was to add the mean anomaly, $M$, to the state vector, i.e., define,
\[ \tilde{X}=[a,e,i,\omega,\Omega, M]^{\sf T },\]
note that the GVE for $M$ has the form,
\begin{equation}\label{equ:++}
\frac{dM}{dt} = \sqrt{\frac{\mu}{a^3}}+
 \sqrt{1-e^2} \frac{(p \cos(\theta)-2 r e)}{ \sqrt{p\mu} e} S-
\frac{ \sqrt{1-e^2} (p+r) \sin(\theta)}{ \sqrt{p\mu}  e} T,
\end{equation}
so that
$$\dot{\tilde{X}}=\tilde{F}(\tilde{X})+\tilde{G}(\tilde{X},\theta) U,$$
and define
$$\tilde{V}_0(\tilde{X})=\frac{1}{2} (\tilde{X}-\tilde{X}_{\tt des}(t))^{\sf T} \tilde{P} (\tilde{X}-\tilde{X}_{\tt des}(t)).$$
The mean anomaly, $M$, is used in $\tilde{X}$ because the corresponding desired mean anomaly is a linear function of time $M_{\tt des}(t) = \sqrt{\frac{\mu}{a_{\tt des}^3}}(t-t_0)+M_{\tt des}(t_0)$ and behaves more regularly compared to  the true anomaly; in addition,  the drift term $\tilde{F}$ is simpler (does not depend on $\theta$).

Considering the Lyapunov function derivative along the trajectories,
$$\frac{d\tilde{V}_0}{d t}=(\tilde{X}-\tilde{X}_{\tt des})^{\sf T} \tilde{P} (\tilde{F}(\tilde{X})-\dot{\tilde{X}}_{\tt des})+
(\tilde{X}-\tilde{X}_{\tt des})^{\sf T} \tilde{P} \tilde{G} U,
$$
leads to the following feedback law replacing (\ref{equ:u_nocon}):
\begin{equation}\label{equ:Lyapunovmeananomalycontroller}U=-\mathrm{Proj}_\mathcal{U} \left[\tilde{G}^{\sf T}(\tilde{X},\theta) \tilde{P} (\tilde{X}-\tilde{X}_{\tt des})  \right],\end{equation}
to which the barrier function terms are similarly added as in (\ref{equ:controllaw1}) to handle state constraints.

\begin{figure}[h!]
	\centering
	\includegraphics[width=0.4\linewidth]{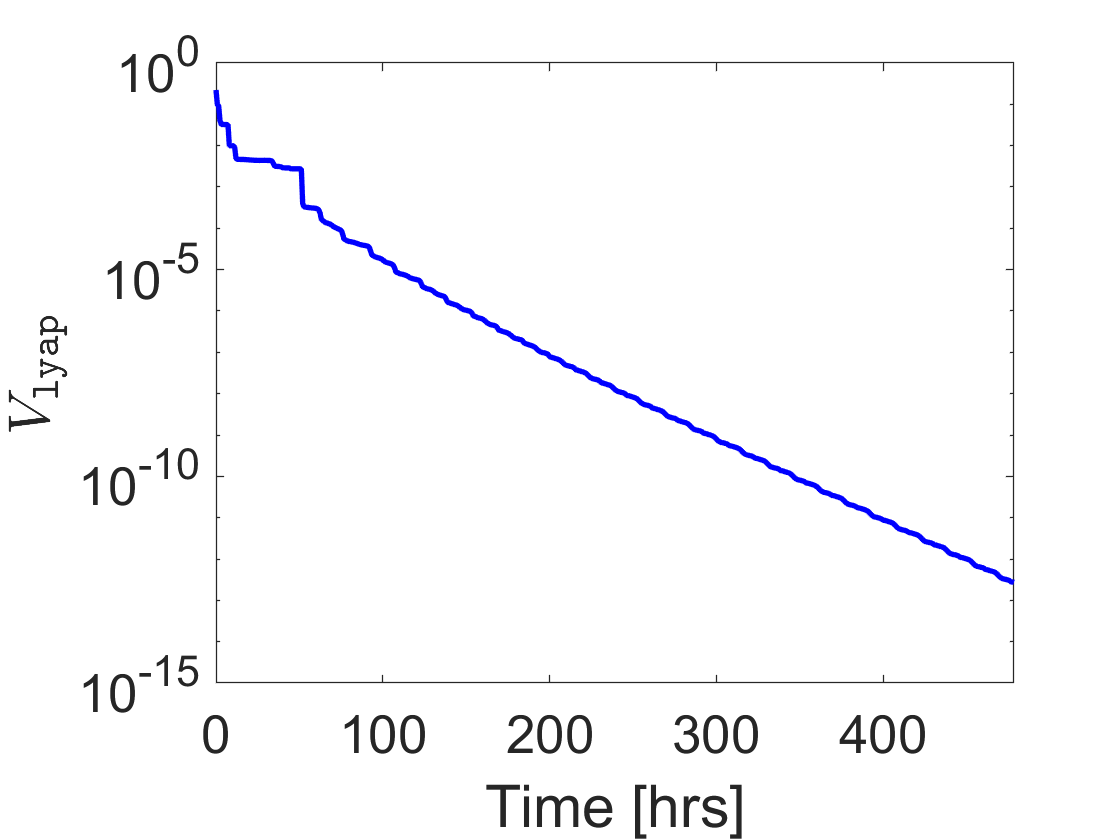}
	\includegraphics[width=0.4\linewidth]{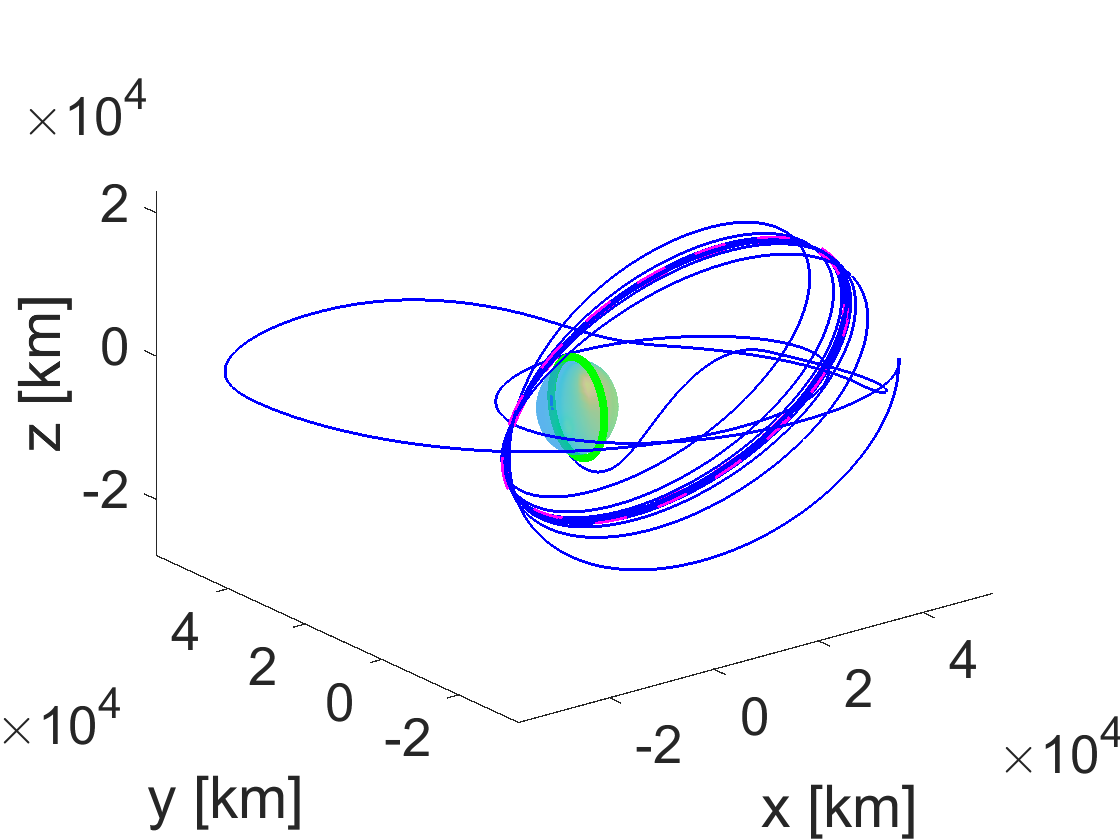}
	\caption{Rendezvous maneuver to higher semi-major axis orbit: Lyapunov function $\tilde{V}_0(\tilde{X}(t))$ (left) and three dimensional trajectory (right) with $\mathcal{U}$ defined by two-sided constraints on $S$ and $W$. Initial orbit is shown in green and target orbit in dashed magenta.}
	\label{fig:twosided}
\end{figure}

\begin{figure}[h!]
	\centering
	\includegraphics[width=0.4\linewidth]{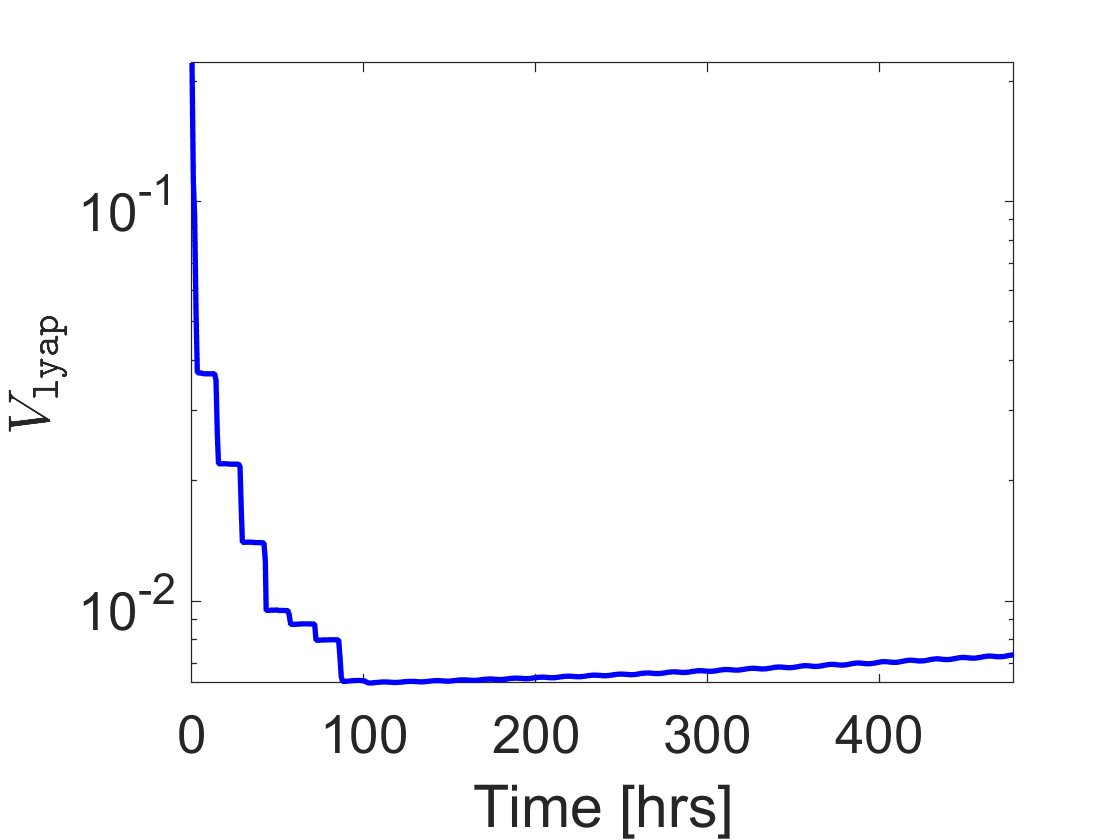}
	\includegraphics[width=0.4\linewidth]{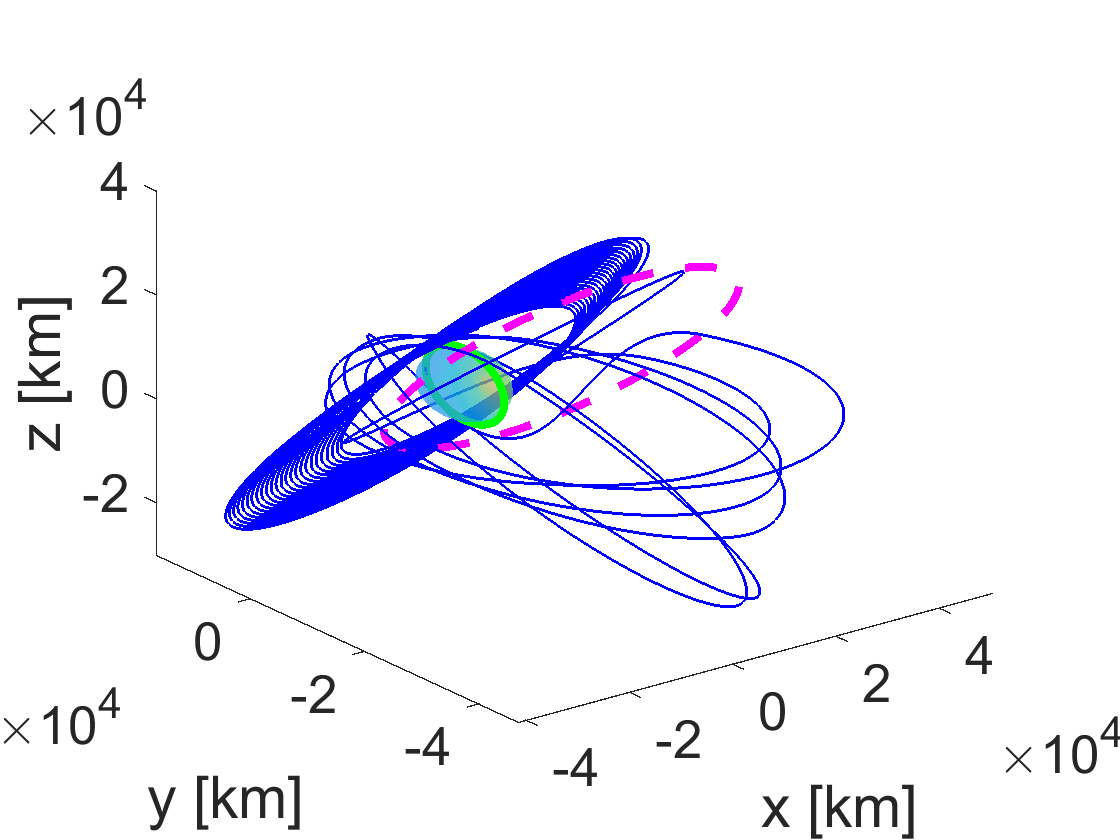}
	\caption{Rendezvous maneuver to higher semi-major axis orbit: Lyapunov function $\tilde{V}_0(\tilde{X}(t))$  (left) and three dimensional trajectory (right) with $\mathcal{U}$ defined by one-sided constraints on $S$ and $W$. Initial orbit is shown in green and target orbit in dashed magenta.}
	\label{fig:onesided}
\end{figure}

Figure~\ref{fig:twosided} shows that when $\mathcal{U}=\{U:~\|U\|\leq 0.001\}$ is defined by two-sided constraints on $S$ and $W$, and the rendezvous maneuver to an orbit with a larger semi-major axis is considered, the trajectory (including mean anamoly) converges to the desired one as $\tilde{V}_0(\tilde{X}(t)) \to 0$ as $t \to \infty$; but in the case of $\mathcal{U}$ defined by one-sided constraints on $S$ and $W$, 
$\mathcal{U}= \{U:~\| U \|\leq 0.001\} 
\bigcap \{U:~S \geq 0,~W \geq 0\}$, 
the convergence does not take place, see Figure~\ref{fig:onesided}. 
The parameters in these simulations were:
$\tilde{P}=diag( 4.0 \times 10^{-10},0.01, 0.005, 0.0075, 0.0005, 0.0005)$, $\mu    = 398600.4405$ (km$^3$/s$^2$), 
$R_e = 6378$ (km),
$r_{\tt min} = R_e+250$, $e_{\tt min}=10^{-3}$,
$\epsilon_1=50$, $\epsilon_2=5 \times 10^{-4}$,
$\tilde{X}_{\tt des}(t)=[R_e+25000,0.6, 
-\frac{\pi}{3},
0, \pi, \sqrt{\frac{\mu}{(R_e+25000)^3}}t]^{\sf T}$,
$\tilde{X}(0)=[R_e+500,       0.02,      \frac{\pi}{2},  \frac{3 \pi}{2},  \pi, 0]^{\sf T}$. In these simulations, the constraints (\ref{equ:maincon3_1}), (\ref{equ:maincon3_2}) were not imposed.  Similar conclusions have been reached for other choices of $\tilde{P} \succ 0$ that were tried. 

Note that even in the case of two-sided control constraints,  closed-loop stability guarantees with the controller (\ref{equ:Lyapunovmeananomalycontroller}) are lacking.  

\subsection{Cascade Inner Loop-Outer Loop Control Architecture}

\begin{figure}[h!]
    \centering
    \includegraphics[width=0.7\linewidth]{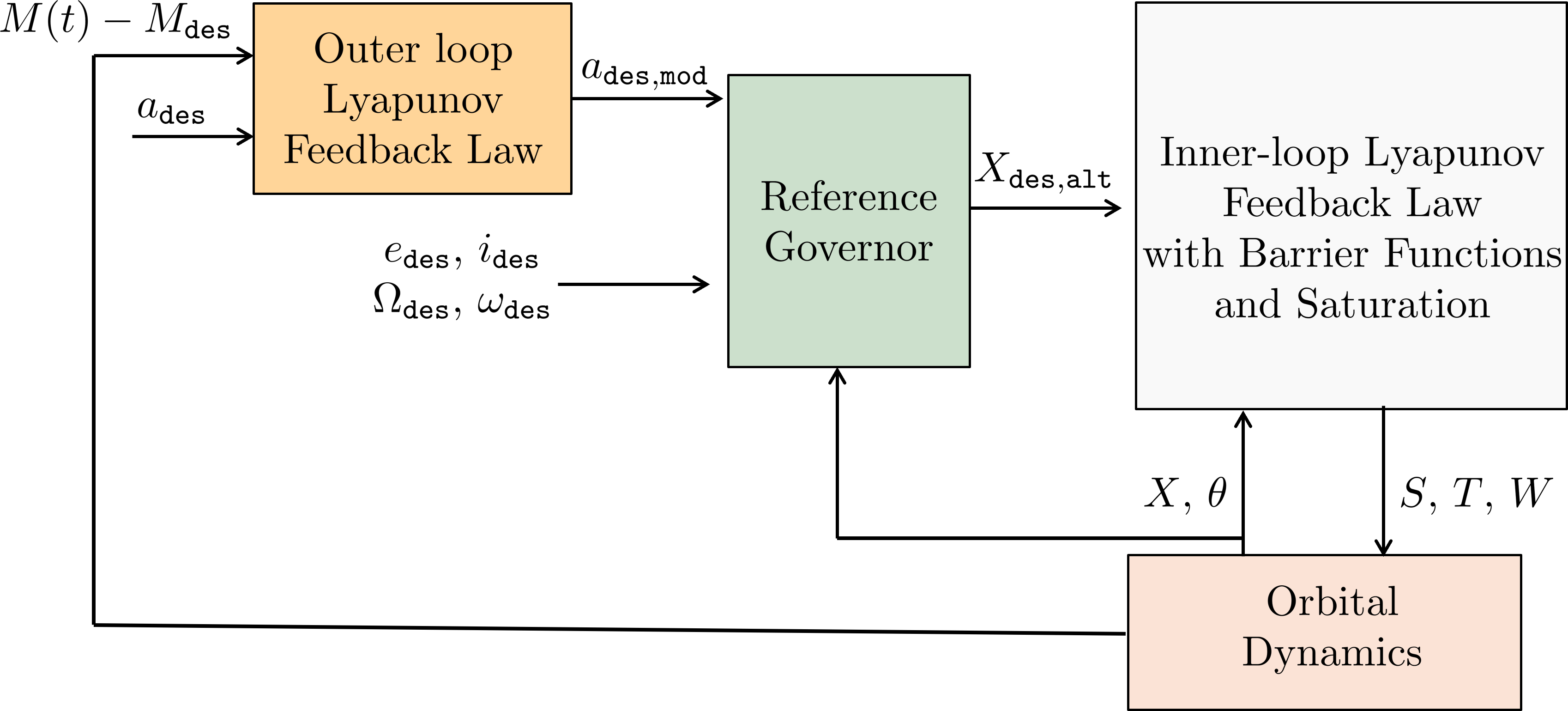}
    \caption{Overall control scheme schematics.}
    \label{fig:controlschemeschematics}
\end{figure}

Hence an alternative approach is pursued that relies on proven convergence properties of the Lyapunov feedback law under control constraints for the five orbital elements and an outer loop controller which adjusts the semi-major axis command, $a_{\tt des}$, based on the difference between the desired mean anomaly and actual mean anomaly.   The approach is partly inspired by \cite{kolmanovsky1996switched} (and references therein), where the control of fiber variables ($M$ is our case) is through the control of base variables ($X$ in our case) to induce the desired geometric phase.
Utilizing the error between the desired true anomaly and actual true anomaly for feedback has also been considered, but it was found that the mean anomaly error leads to better responses in simulations, in particular, because of simpler (linear) dependence of the desired mean anomaly on time.  Notably, a similar approach is adopted in \cite{Lantukh2017EnhancedQLaw,naasz2002classical} where the target for the semi-major axis is adjusted based on the true longitude error for the case of the $Q$-law; however, in our case the use of an outer loop feedback has become necessary because of the one-sided constraints on $S$ and $W$  (which were not considered in \cite{Lantukh2017EnhancedQLaw,naasz2002classical}).  The schematics of the overal scheme are shown in Figure~\ref{fig:controlschemeschematics}.

As $M-M_{\tt des}$ is an angle, its dynamics can be viewed as evolving on a circle, and the state can be represented as \cite{marley2020kinematic}
$$z=\begin{bmatrix}
    \cos(M-M_{\tt des}) \\
    \sin(M-M_{\tt des})
\end{bmatrix}.$$
Then,
\begin{equation}\label{equ:zdot}
\dot{z}=\begin{bmatrix}
    - \sin (M-M_{\tt des}) \\
    \cos(M-M_{\tt des}) 
\end{bmatrix} (\dot{M}-\dot{M}_{\tt des})=
 S z (\dot{M}-\dot{M}_{\tt des}),\qquad S=\begin{bmatrix}
    0 & -1 \\ 1 & 0 
\end{bmatrix} .
\end{equation}
Defining  a Lyapunov function candidate,
\begin{equation}\label{equ:AA}
\bar{V}=\frac{1}{2}(z-z_{\tt des})^{\sf T} \bar{P} (z-z_{\tt des}),\quad 
z_{\tt des}=\begin{bmatrix}
    1 \\ 0
\end{bmatrix}, \end{equation}
where $$\bar{P}=\begin{bmatrix} 
\bar{p}_1 & 0 \\ 0 & \bar{p}_2 \end{bmatrix}
 \succ 0,$$
it follows that 
\begin{equation}\label{equ:BB}
\dot{\bar{V}}=(z-z_{\tt des})^{\sf T} \bar{P} S z (\dot{M}-\dot{M}_{\tt des}).
\end{equation}
If the control law for $\dot{M}$ is prescribed as 
\begin{equation}\label{equ:zcontrol}
\dot{M} = \dot{M}_{\tt des}-\lambda z^{\sf T} S^{\sf T} \bar{P}  (z-z_{\tt des}),~\lambda>0,
\end{equation}
then $\dot{\bar{V}} \leq 0$.

The properties of the closed-loop system (\ref{equ:zdot})-(\ref{equ:zcontrol}) are summarized as follows.

{\bf Proposition:} 
Assume that  $\bar{p}_2 \leq 2 \bar{p}_1$.  Then the equilibrium at $z_{\tt des}$ of the closed-loop system (\ref{equ:zdot})-(\ref{equ:zcontrol}) is locally exponentially stable with the region of attraction being the set $\{z:~z=[\cos(a),\sin(a)]^{\sf T},~-\pi<a<\pi\}$. Furthermore, the closed-loop system is locally input-to-state stable (LISS)
at $z_{\tt des}$.

{\bf Proof:}
Consider 
\[
z=\begin{bmatrix}\cos a\\ \sin a\end{bmatrix}, \]
and let
$\rho=\sin^2\left(\dfrac{a}{2}\right)$. 
Suppose that $\rho \in [0, 1-\epsilon]$ for $\epsilon>0$.
Then, from (\ref{equ:AA}), $$\bar{V}=2 \rho \left(\bar{p}_2+(\bar{p}_1-\bar{p}_2)\rho \right),$$
and from (\ref{equ:BB}),
$$\dot{\bar{V}}=4 \lambda \rho (\rho-1) 
\left(\bar{p}_2+2(\bar{p}_1-\bar{p}_2) \rho \right)^2,$$
leading to
\begin{eqnarray*}
-\dfrac{\dot{\bar{V}}}{\bar{V}} &= &
\dfrac{2 \lambda (1-\rho) 
\left(\bar{p}_2+2(\bar{p}_1-\bar{p}_2) \rho \right)^2}{
\bar{p}_2+(\bar{p}_1-\bar{p}_2) \rho}.
\end{eqnarray*}

The conditions, $\bar{p}_2 \leq 2 \bar{p}_1$, $0 \leq \rho\leq 1-\epsilon$, ensure that
$\bar{p}_2+2(\bar{p}_1-\bar{p}_2)\rho>0$
and
\begin{equation}\label{equ:liss}
-\dfrac{\dot{\bar{V}}}{\bar{V}} \leq
\bar{c}(\epsilon)\stackrel{\Delta}{=}\dfrac{2 \lambda  
\max\{\bar{p}_2^2, \bar{p}_2+2(\bar{p}_1-\bar{p}_2)(1-\epsilon)^2\}}
{\min\{\bar{p}_2,\bar{p}_2+(\bar{p}_1-\bar{p}_2)(1-\epsilon)\}}>0.\end{equation}
We also note that the condition $\bar{p}_2 \leq 2 \bar{p}_1$ ensures that $\bar{V}=2 \rho (\bar{p}_2+(\bar{p}_1-\bar{p}_2)\rho)$ 
is monotonically non-decreasing with respect to $\rho \in [0,1-\epsilon]$ as
$$ 
\frac{\partial \bar{V} } {\partial \rho} 
= 2 \bar{p}_2 + 4(\bar{p}_1
-\bar{p_2})\rho \geq \min\{ 2 \bar{p}_2,~ 2(2\bar{p_1}-\bar{p}_2) \}\geq 0.$$
Thus $\dot{\bar{V}}(t) \leq 0$ along the closed-loop trajectories implies that ${\rho}(t)$ is nonincreasing along the closed-loop trajectories and $\rho(0) \in [0,1-\epsilon]$ implies $\rho(t) \in [0,1-\epsilon]$ for all $t \geq 0$. 
  Thus (\ref{equ:liss}) holds along trajectories implying local exponential stability.
If $$
z(0)=\begin{bmatrix}\cos a\\ \sin a\end{bmatrix},~-\pi<a<\pi,$$
then $\rho(0)=\sin^2\left(\dfrac{a}{2}\right)<1-\epsilon$ for some $\epsilon>0$ implying that $z(0)$ is in the region of attraction.
The LISS at the origin for the closed-loop system with a disturbance input, $d$,
$$\dot{z}=\bar{S} z \left((\dot{M}-\dot{M}_{\tt des})+d\right),$$
 follows as 
 $\bar{V}$ can be shown to be a (local) ISS-Lyapunov function.  \qed

If $|S|$ and $|T|$ are small, it follows from  (\ref{equ:++}) that
$$\dot{M} \approx n = \sqrt{\dfrac{\mu}{a^3}}.$$
At the same time, the desired mean motion on the target orbit satisfies,
$$\dot{M}_{\tt des}=n_{\tt des} = \sqrt{\dfrac{\mu}{a_{\tt des}^3}}.$$
Therefore,
based on (\ref{equ:zcontrol}), we define
the mean motion, $n_{\tt mod}$, aimed at adjusting the phasing and achieving rendezvous, as follows:
\begin{equation}\label{equ:nmod}
n_{\tt mod}=\min\!\Big(
4\,n_{\mathrm{des}},\;
\max\!\big(
0.25\,n_{\tt {des}},\;
n_{\tt {des}}
-
\lambda\, z^\top \bar{S}^\top \bar{P} (z - z_{\tt des})
\big)
\Big).
\end{equation}
The modified semi-major axis command, $a_{\tt des,mod}$, can now be computed to yield the prescribed $n_{\tt mod}$: 
\begin{equation}\label{equ:amod}
a_{\tt {des},mod}
=
\max\!\Big(
0.8 a_{\tt des},\;
\min\!\big(
1.2 a_{\tt des},\;
\left(\frac{\mu}{{n}_{\tt mod}^{2}}\right)^{\!\tfrac{1}{3}}
\big)
\Big).
\end{equation}
Note that $n_{\tt mod}$, $a_{\tt {des},mod}$ are saturated between suitably defined minimum and maximum values; these limits can be refined through simulations  for a particular maneuver. 

\subsection{Simulation Case Studies}

In the subsequent simulations, the following values are used. The gravitational parameter is
$\mu    = 398600.4405$ (km$^3$/s$^2$), Earth radius is
$R_e = 6378$ (km), the minimum radius of periapsis is $r_{\tt min} = R_e+250$ (km), the minimum eccentricity is
$e_{\tt min} =10^{-3}$, the maximum eccentricity is $e_{\tt max}=0.85$, and the minimum inclination is $i_{\tt min}=0.1$.

After manual tuning, we adopted
$P=diag( 4.0 \times 10^{-10},0.01, 0.01, 0.0075, 5 \times 10^{-4})$.  The small value of the first element of $P$ is primarily dictated by the  units of $a$ being km and values of $a$ three orders of magnitude higher than the rest of the orbital elements. 

The
barrier functions parameters are set as 
$\epsilon_1 = 50,~\epsilon_2= 
5 \times 10^{-4},~ \epsilon_3 = 5 \times 10^{-4}, \epsilon_3 =5 \times 10^{-4}$.  The thrust magnitude saturation limits are assumed to be such that
$U_{\tt max}=10^{-3}$ km/sec$^2$.

In the sequel, we firstly consider the maneuver from the lower semi-major axis orbit with 
${X}(0)=[6828,0.01,\frac{\pi}{2},0,\pi]^{\sf T} $
to the higher semi-major axis orbit with
${X}_{\tt des}=[42164,0.8,\frac{\pi}{20},\pi,0]^{\sf T}$.
We then consider the transition from the higher semi-major axis orbit with 
${X}(0)=[42164,0.8,\frac{\pi}{20},\pi,0]^{\sf T}$ to the
lower semi-major axis orbit with
${X}_{\tt des}=[6828,0.01,\frac{\pi}{2},0,\pi]^{\sf T}$.
In both cases, the desired true anomaly trajectory is computed based on
$$\dot{\theta}_{\tt des}(t)=
\frac{\sqrt{\mu} (1+e_{\tt des} \cos(\theta_{des}(t)))^2}{(a_{\tt des} (1-e_{\tt des}^2))^{\frac{3}{2}}},
\quad 
{\theta}_{\tt des}(0)=0.$$


For the outer loop feedback law (\ref{equ:zcontrol}), we set
$$\bar{P}=\begin{bmatrix}
    0.05 & 0 \\ 0 & 0.05
\end{bmatrix},$$
and we used $\lambda=2 \times 10^{-4}$ for the transition from the higher semi-major axis orbit to the lower semi-major axis orbit and a more aggressive value of
$\lambda=8 \times 10^{-4}$ for the transition from the lower semi-major axis orbit to the higher semi-major axis orbit.
For the reference governor implementation, 
we used $\varepsilon_0=0.005$ for the transition to the lower semi-major axis orbit and $\varepsilon_0=0.0005$ for the transition to the higher semi-major axis orbit.  The prediction time interval was $\mathcal{T}=[0,22500]$ sec in both cases.

\subsection{Orbital Maneuvers to Higher Orbit with Conventional Thrust Actuation}

\begin{figure}[h!]
	\centering
	\includegraphics[width=0.4\linewidth]{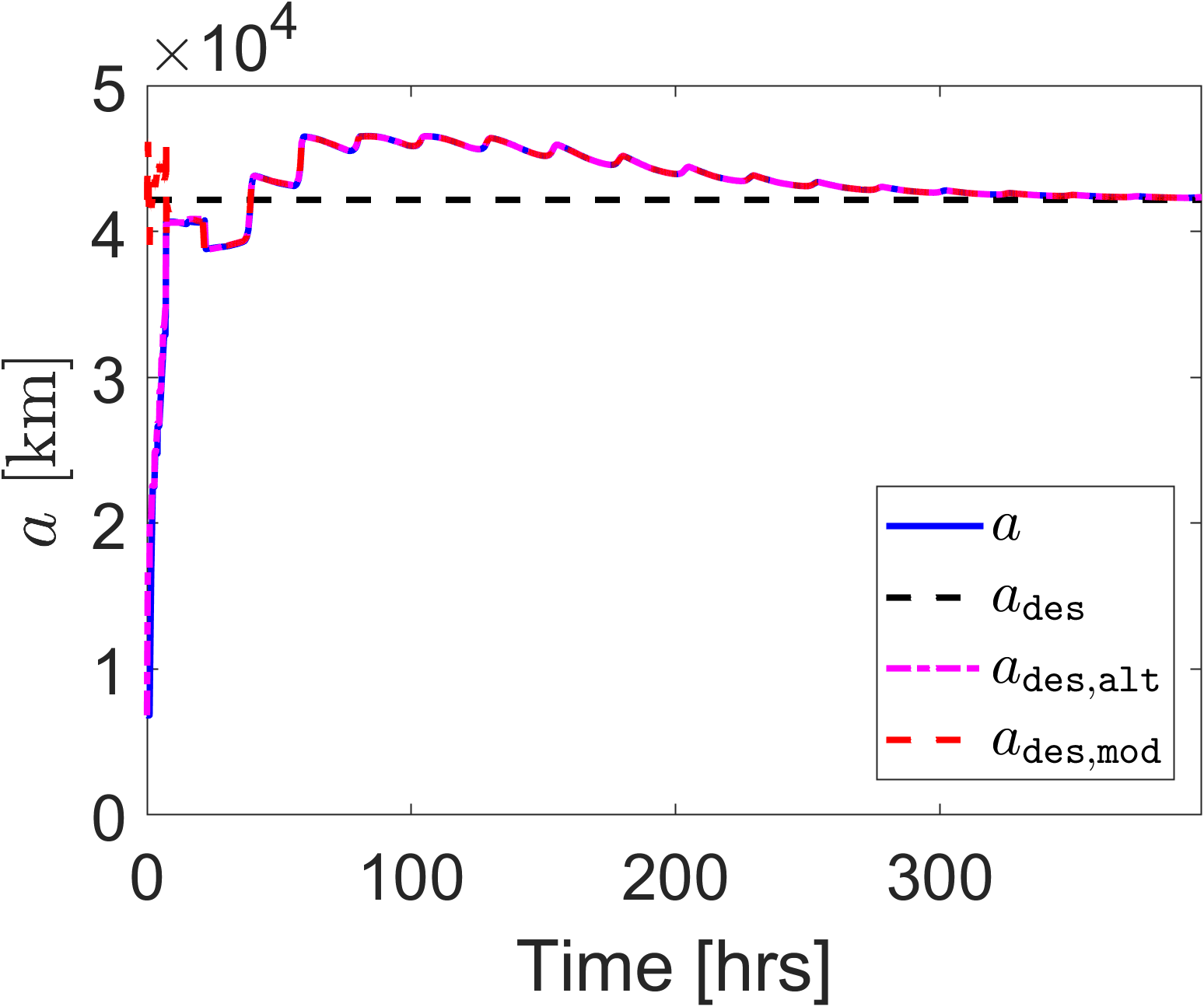}
	\includegraphics[width=0.4\linewidth]{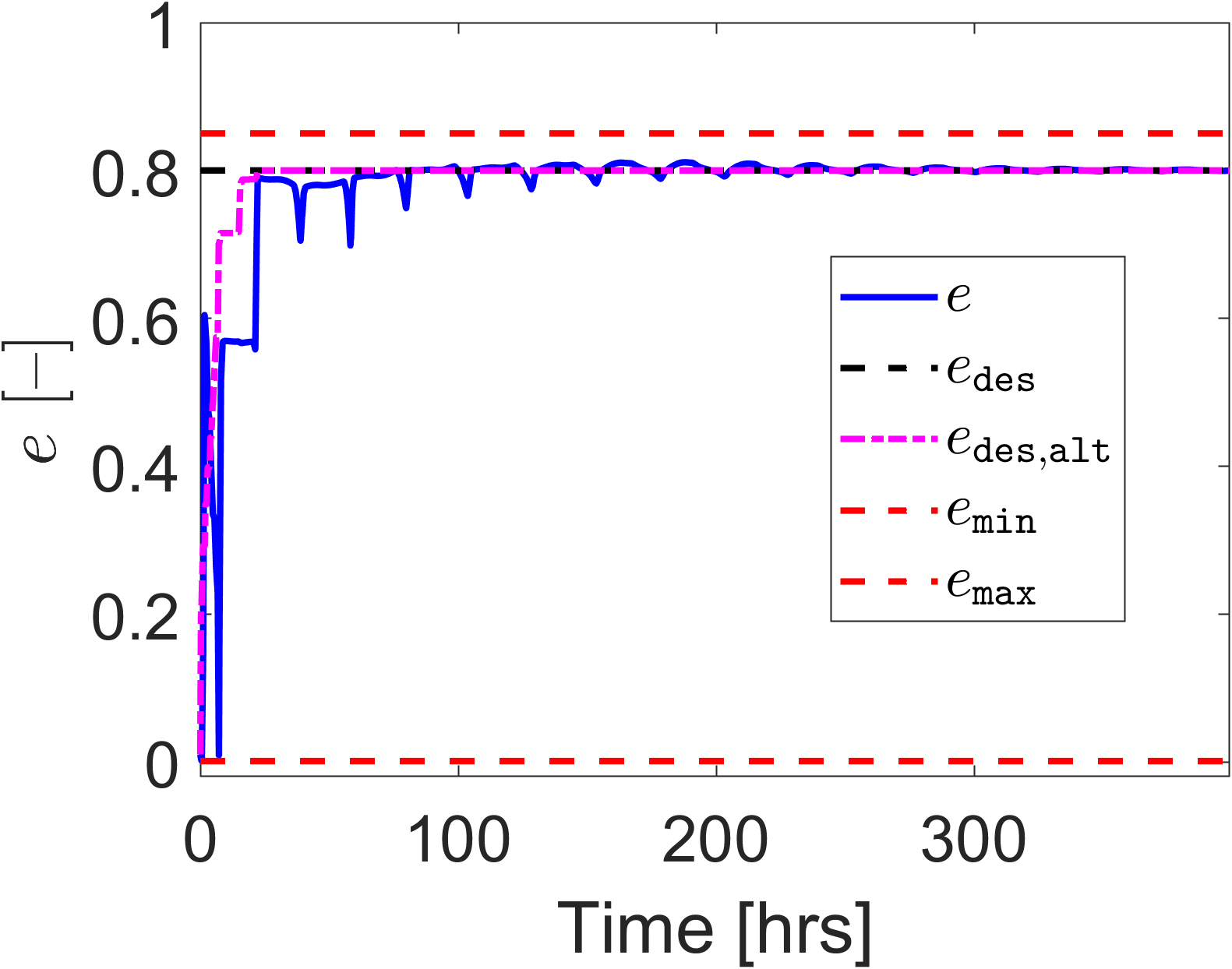}
	\caption{Semi-major axis (left) and eccentricity (right) during transition from a lower semi-major axis orbit to a higher semi-major axis orbit with conventional thrust actuation.}
	\label{fig:aup}
\end{figure}

\begin{figure}[h!]
	\centering
	\includegraphics[width=0.4\linewidth]{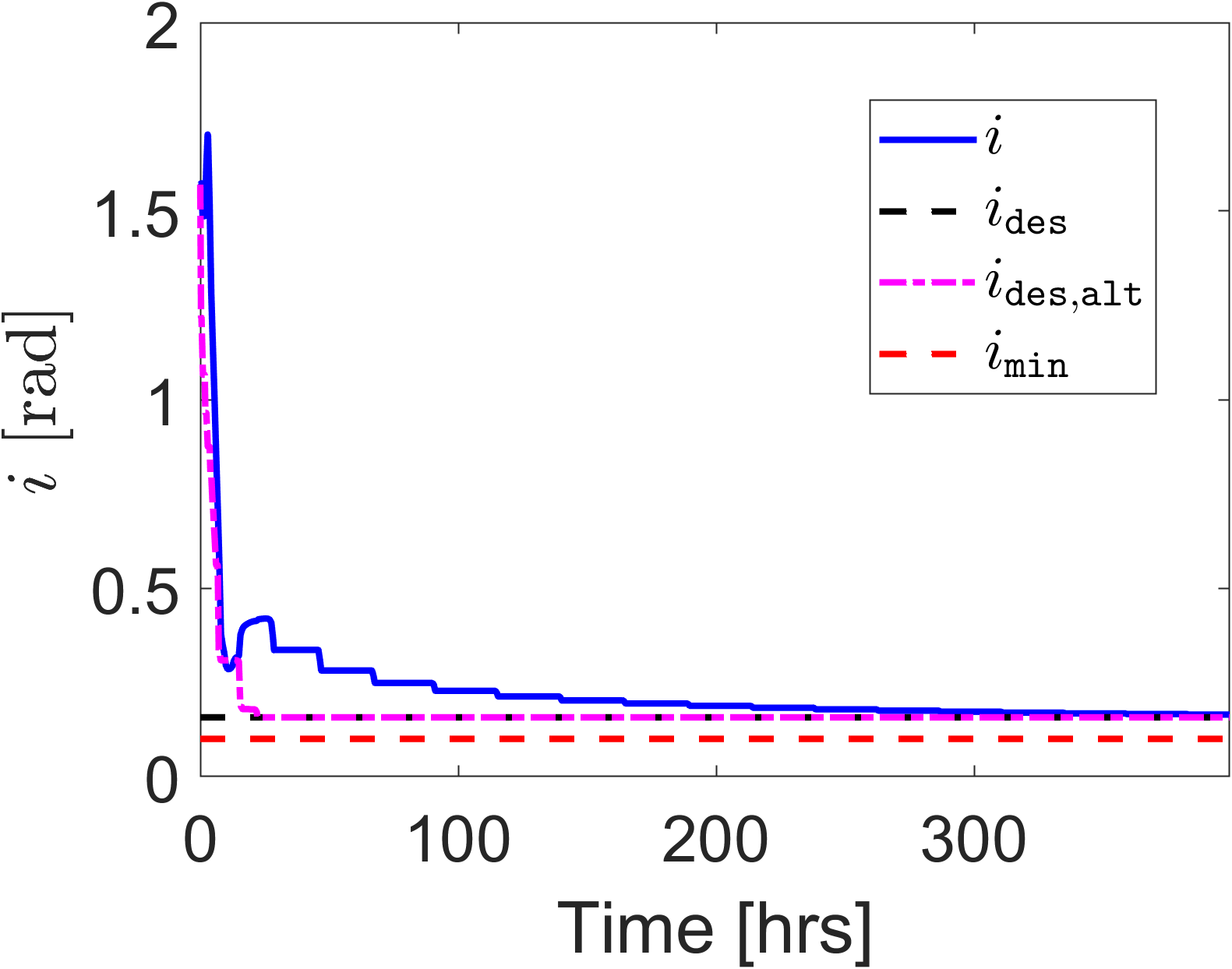}
	\includegraphics[width=0.4\linewidth]{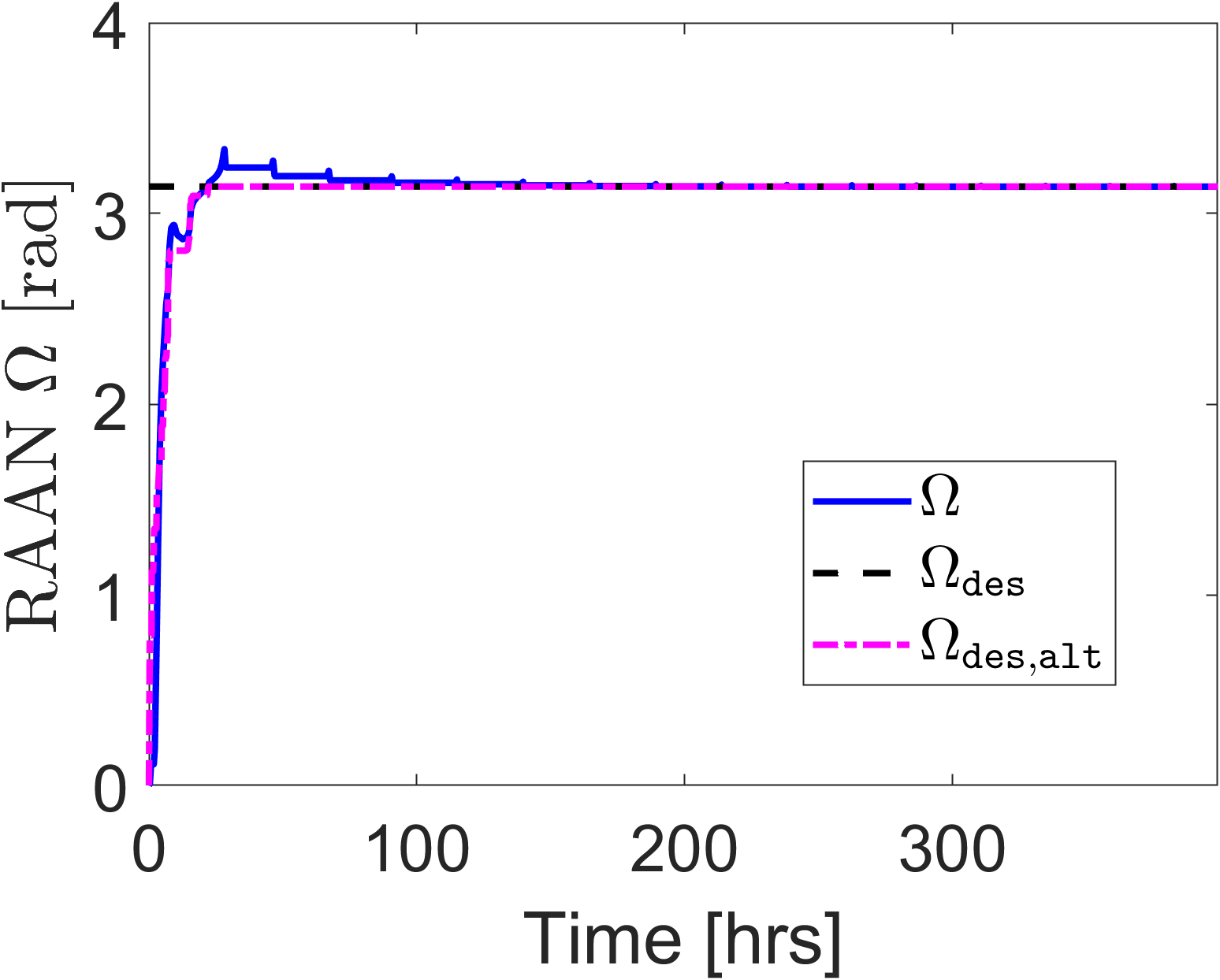}
	\caption{Inclination (left) and RAAN (right) during transition from a lower semi-major axis orbit to a higher semi-major axis orbit with conventional thrust actuation.}
	\label{fig:iup}
\end{figure}

\begin{figure}[h!]
	\centering
	\includegraphics[width=0.4\linewidth]{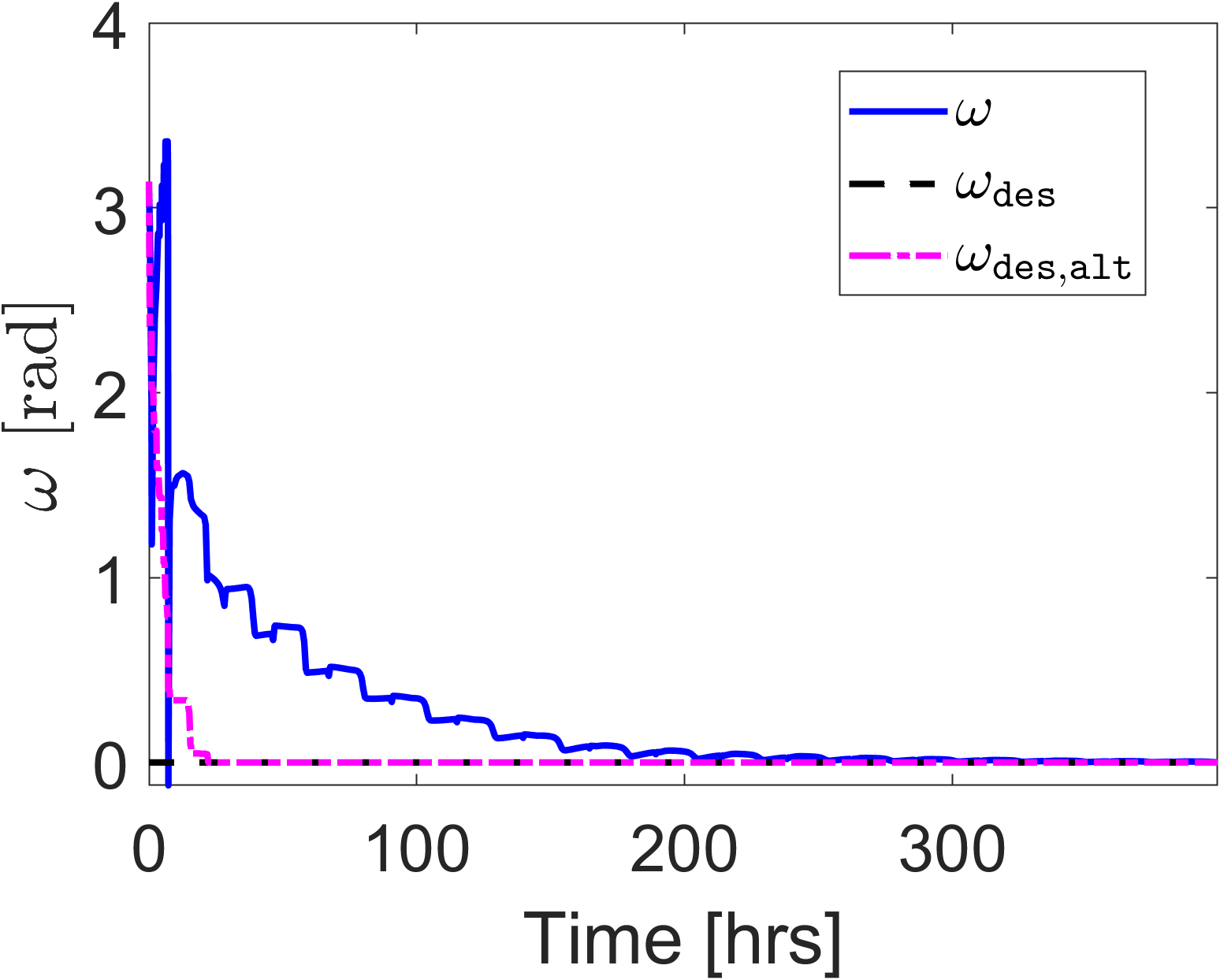}
	\includegraphics[width=0.4\linewidth]{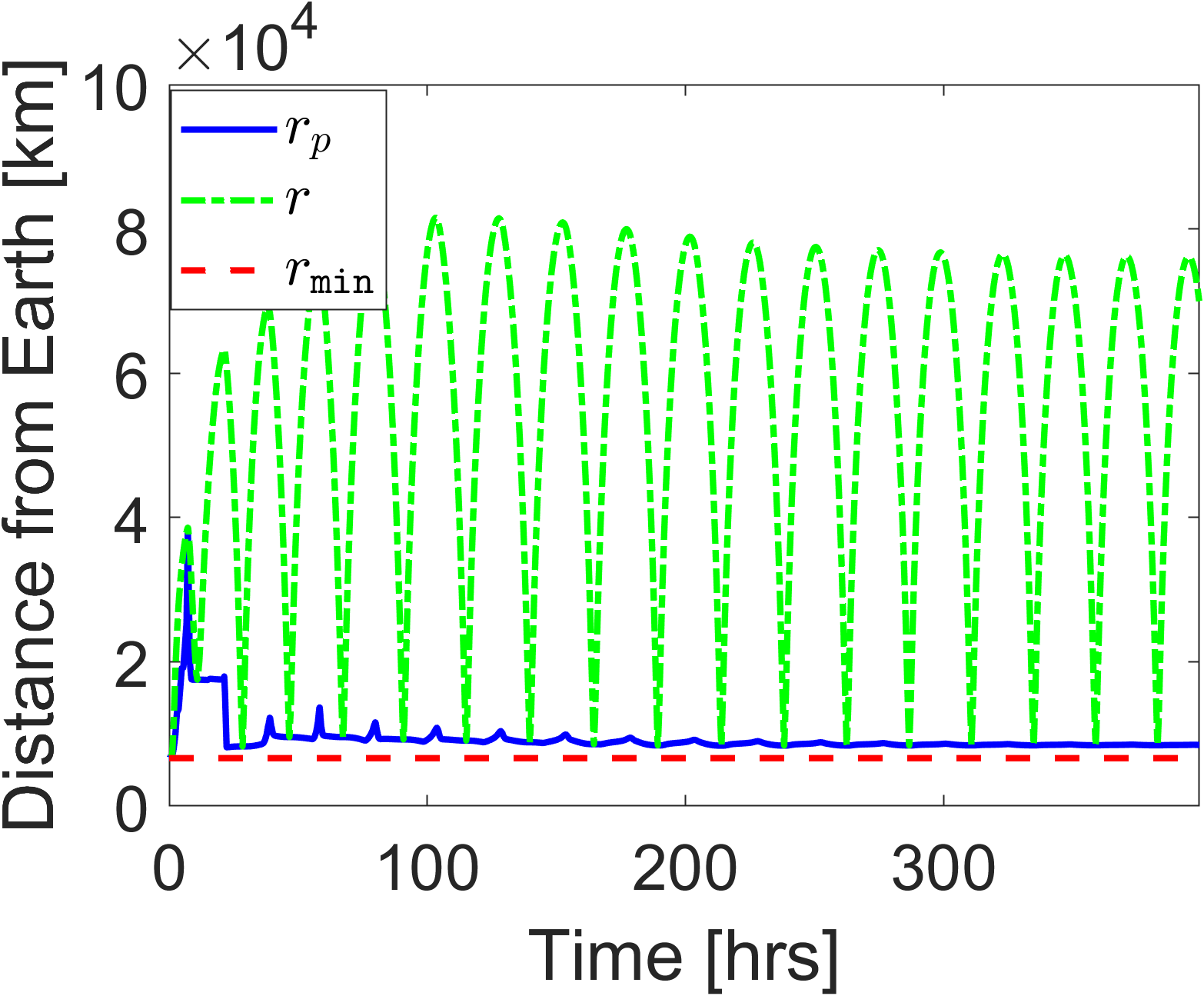}
	\caption{Argument of periapsis (left) and radius of periapsis (right) during transition from a lower semi-major axis orbit to a higher semi-major axis orbit with conventional thrust force actuation.}
	\label{fig:argumentofperiapsisup}
\end{figure}

\begin{figure}[h!]
	\centering
	\includegraphics[width=0.4\linewidth]{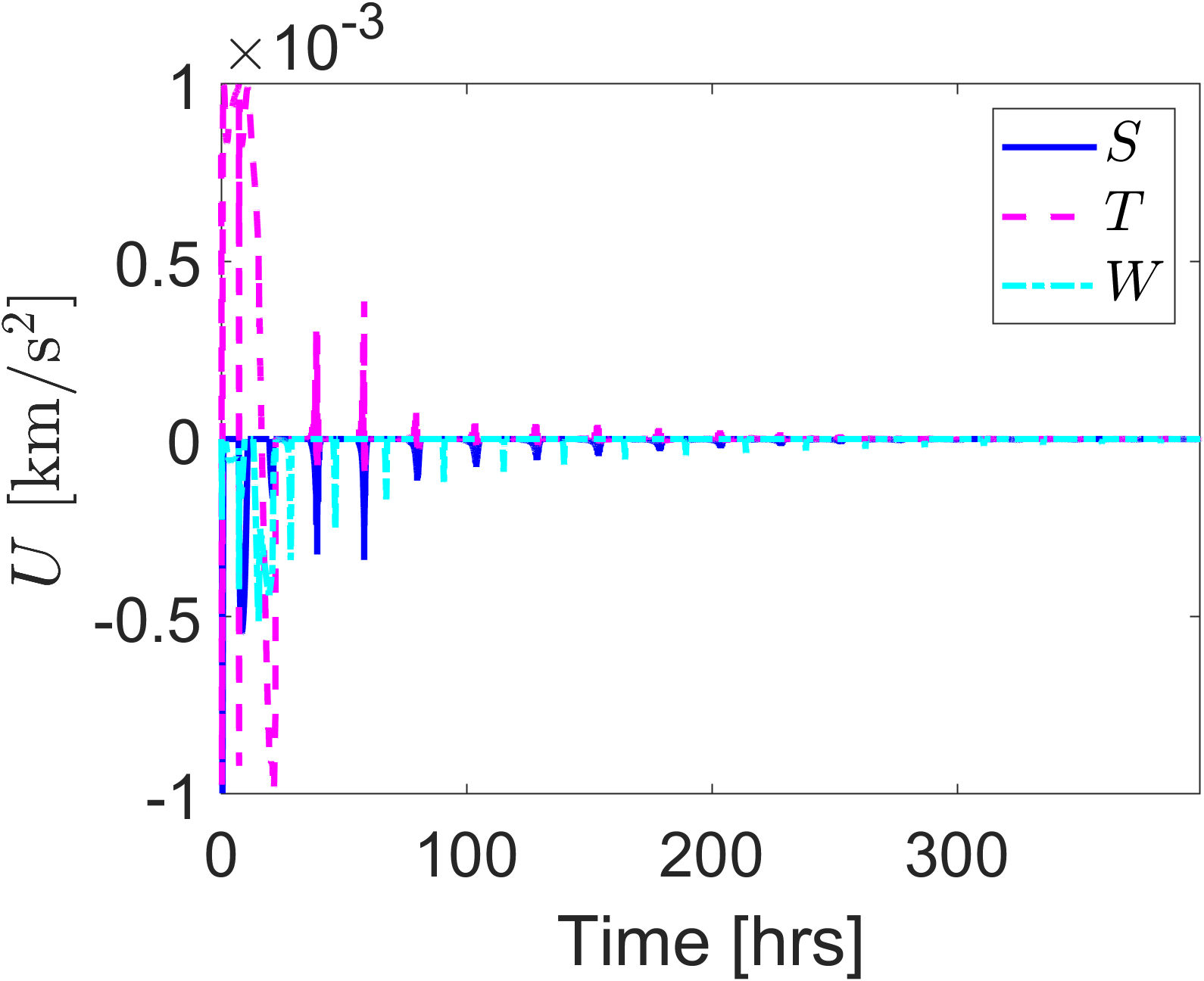}
	\includegraphics[width=0.4\linewidth]{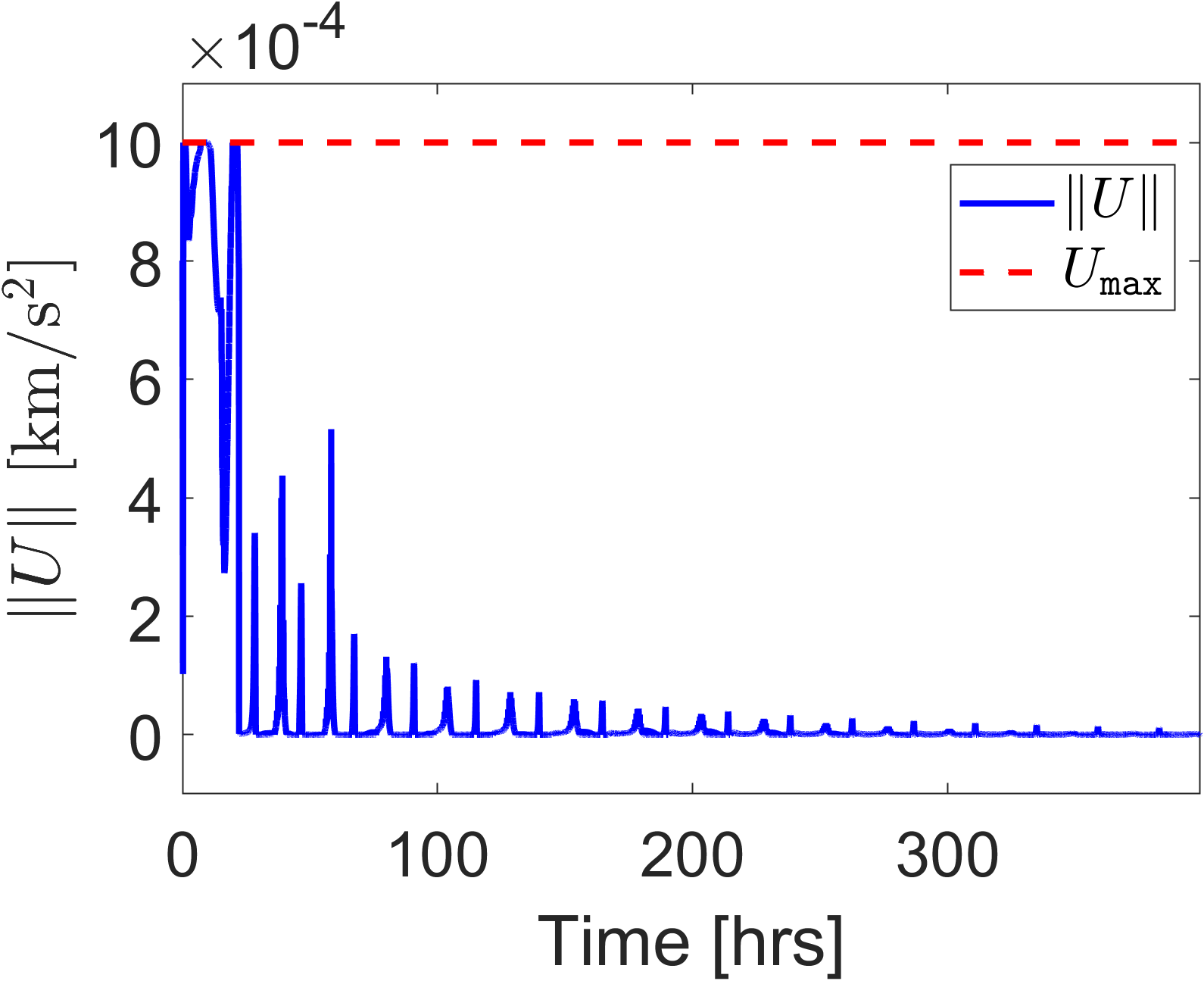}
	\caption{$S$, $T$ and $W$ components (left) and $\|U\|$ (right) during transition from a lower semi-major axis orbit to a higher semi-major axis orbit with conventional thrust actuation.}
	\label{fig:STWup}
\end{figure}

\begin{figure}[h!]
	\centering
	\includegraphics[width=0.4\linewidth]{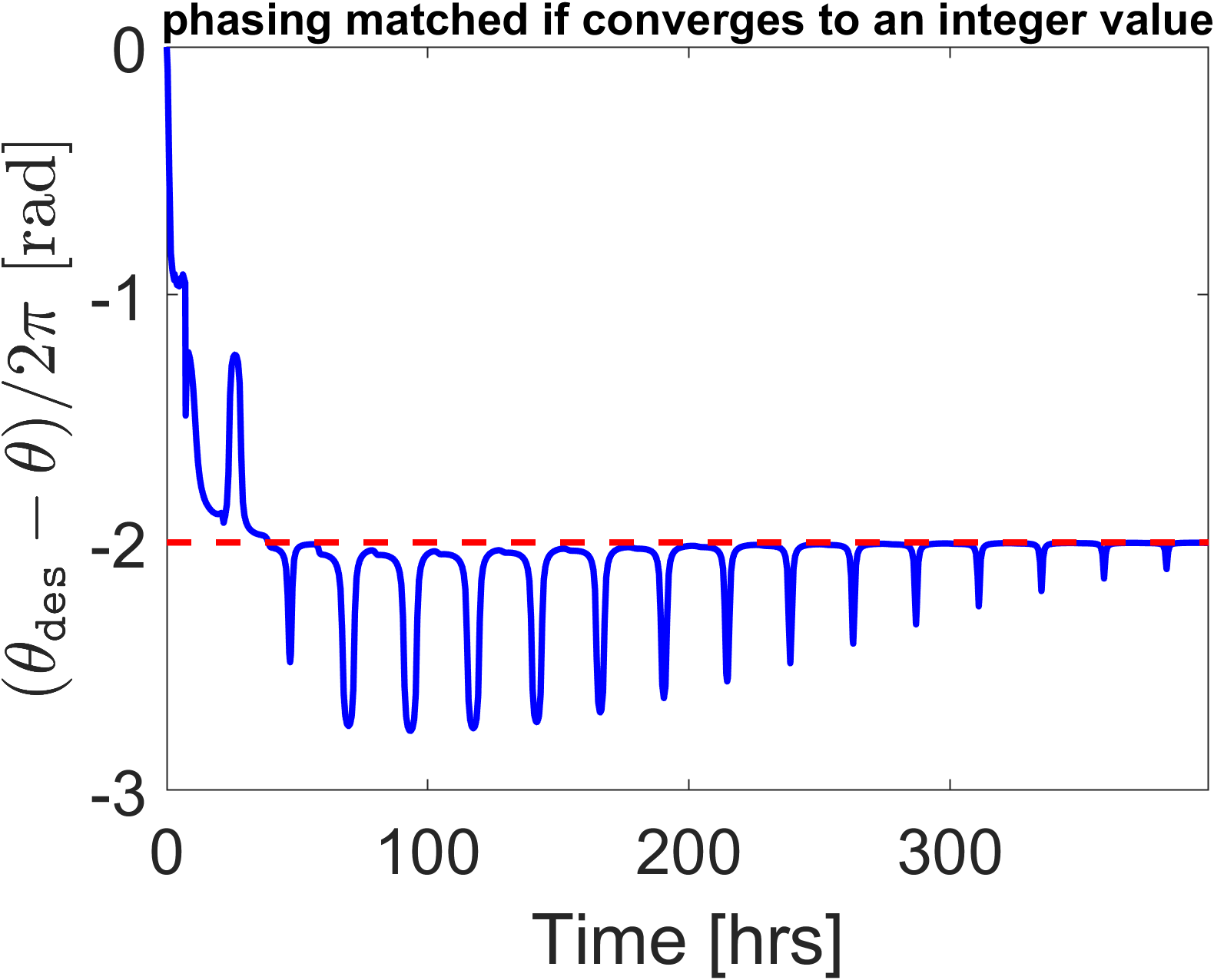}
    \includegraphics[width=0.4\linewidth]{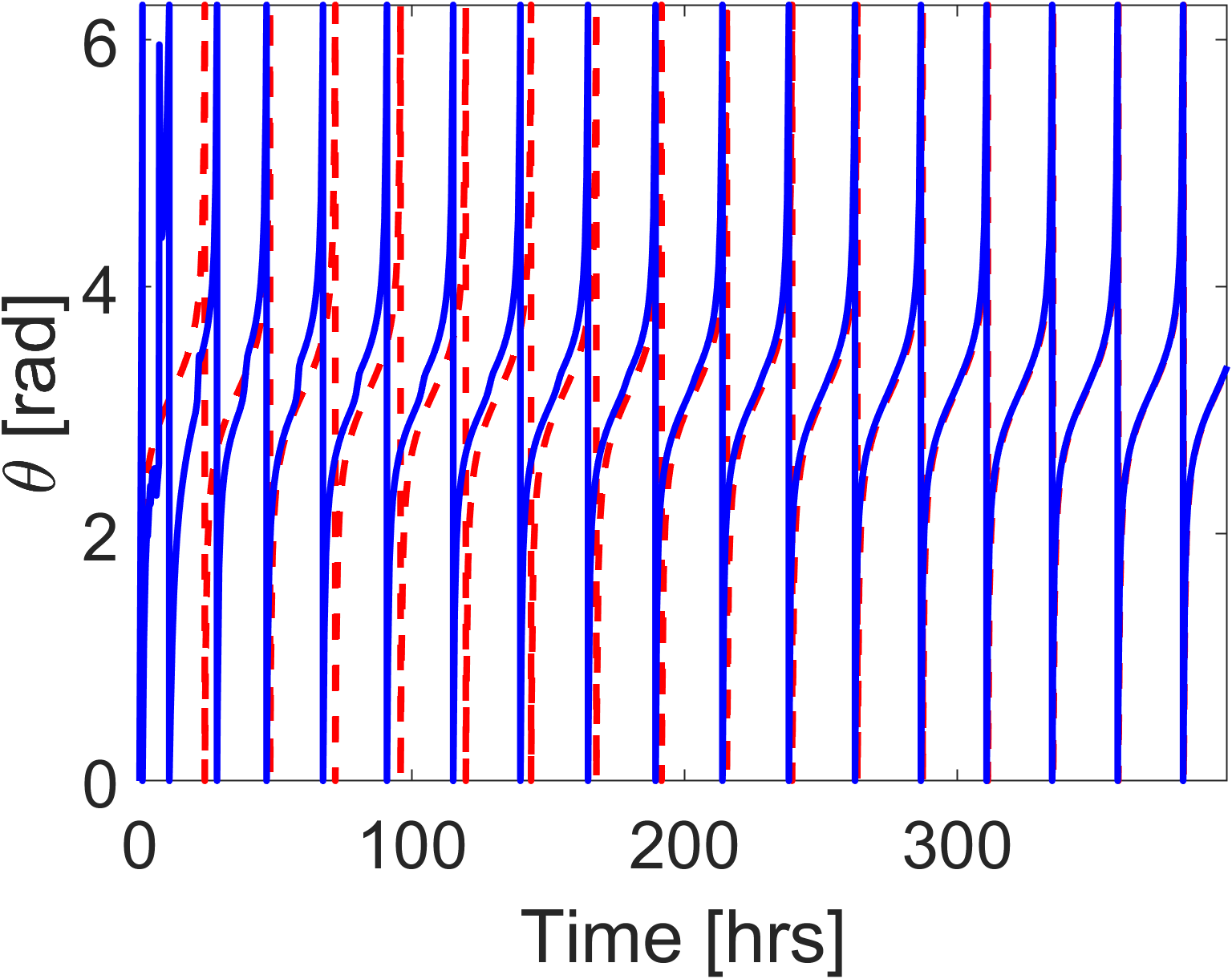}
	\caption{True anomaly difference (left) and true anomalies versus time (right) during transition from a lower semi-major axis orbit to a higher semi-major axis orbit with conventional thrust actuation.
		}
	\label{fig:thetaup}
\end{figure}

\begin{figure}[h!]
	\centering		\includegraphics[width=0.5\linewidth]{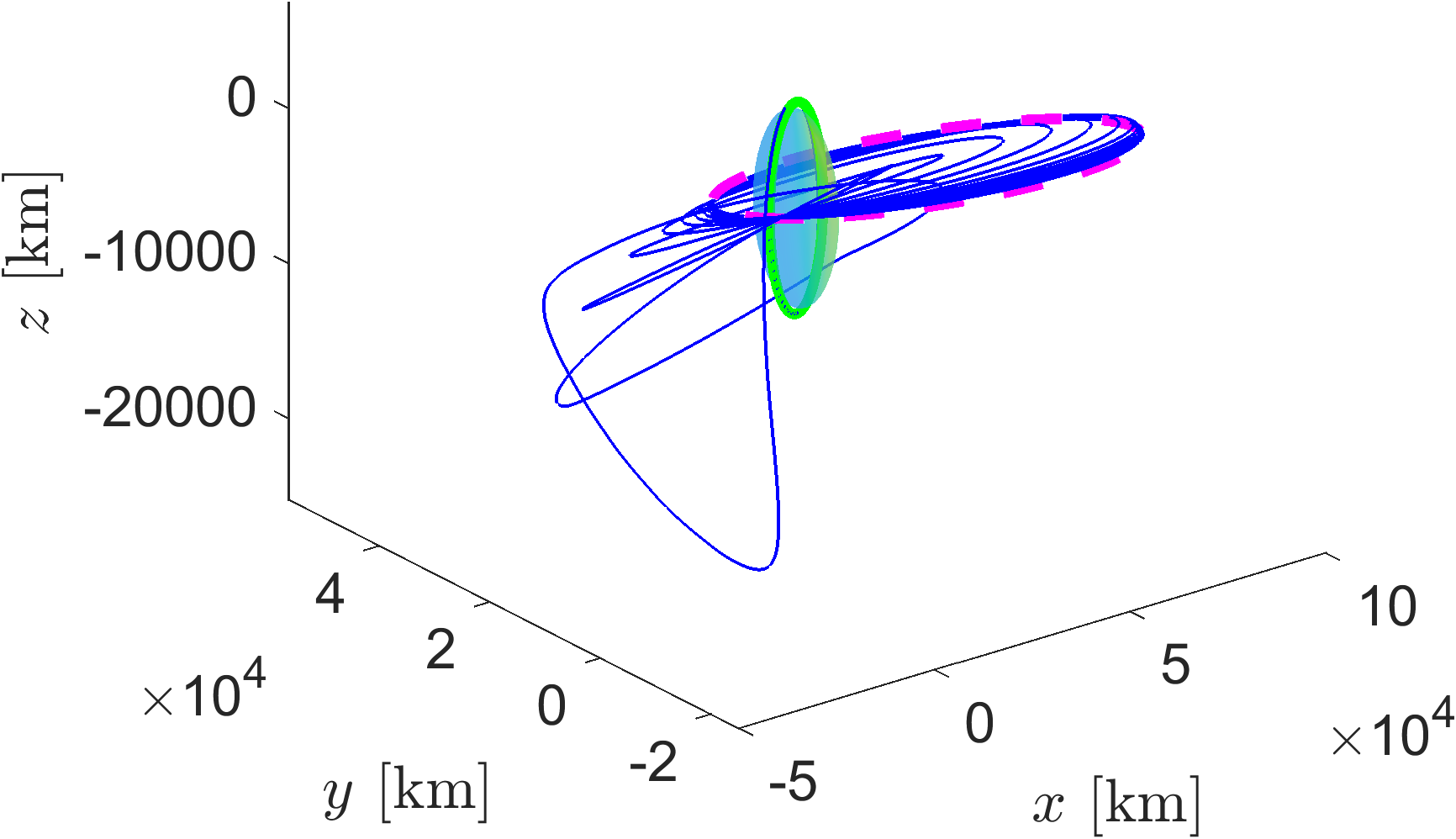}
	\caption{Three dimensional initial orbit (green), transfer orbit (blue) and final (desired) orbit (dashed magenta) during transition from a lower semi-major axis orbit to a higher semi-major axis orbit with conventional thrust actuation. }
	\label{fig:thetaup_1}
\end{figure}

Figures~\ref{fig:aup}-\ref{fig:thetaup_1} illustrate the simulated maneuvers for the spacecraft transferring from a lower semi-major axis orbit and performing rendezvous with a spacecraft on the higher semi-major axis orbit.  


Notably, all the constraints are enforced during the maneuver, including $S \geq 0$, $W \geq 0$ and $\|U \| \leq U_{\tt max}$. See Figure~\ref{fig:STWup}.  

Figure~\ref{fig:thetaup}-right shows that the time history of the true anomaly converges to that of the target thereby achieving rendezvous.
Figure~\ref{fig:thetaup}-left shows the time history of $\frac{\theta_{\tt des}(t)-\theta(t)}{2 \pi}$, which converges to an integer value of $-2$; this also confirms the true anomaly matching by the end of the maneuver, achieved by altering the semi-major axis reference command indicated in Figure~\ref{fig:aup}-left by the red dashed line.  The activation of the outer loop controller and the tracking of the altered semi-major axis command by the Lyapunov feedback controller can also be  observed from Figure~\ref{fig:aup}.  

\begin{figure}[h!]
	\centering
	\includegraphics[width=0.4\linewidth]{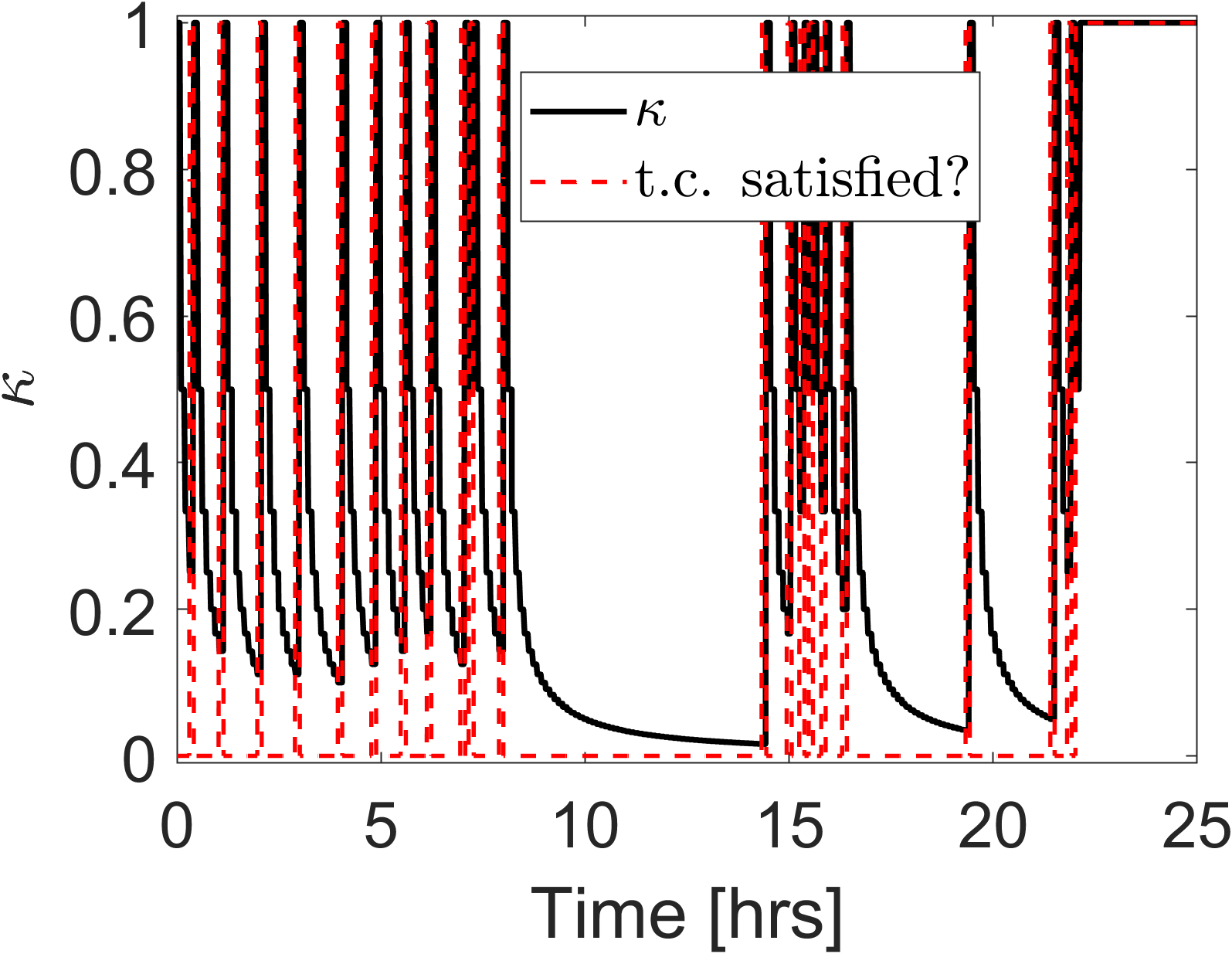}
\includegraphics[width=0.4\linewidth]{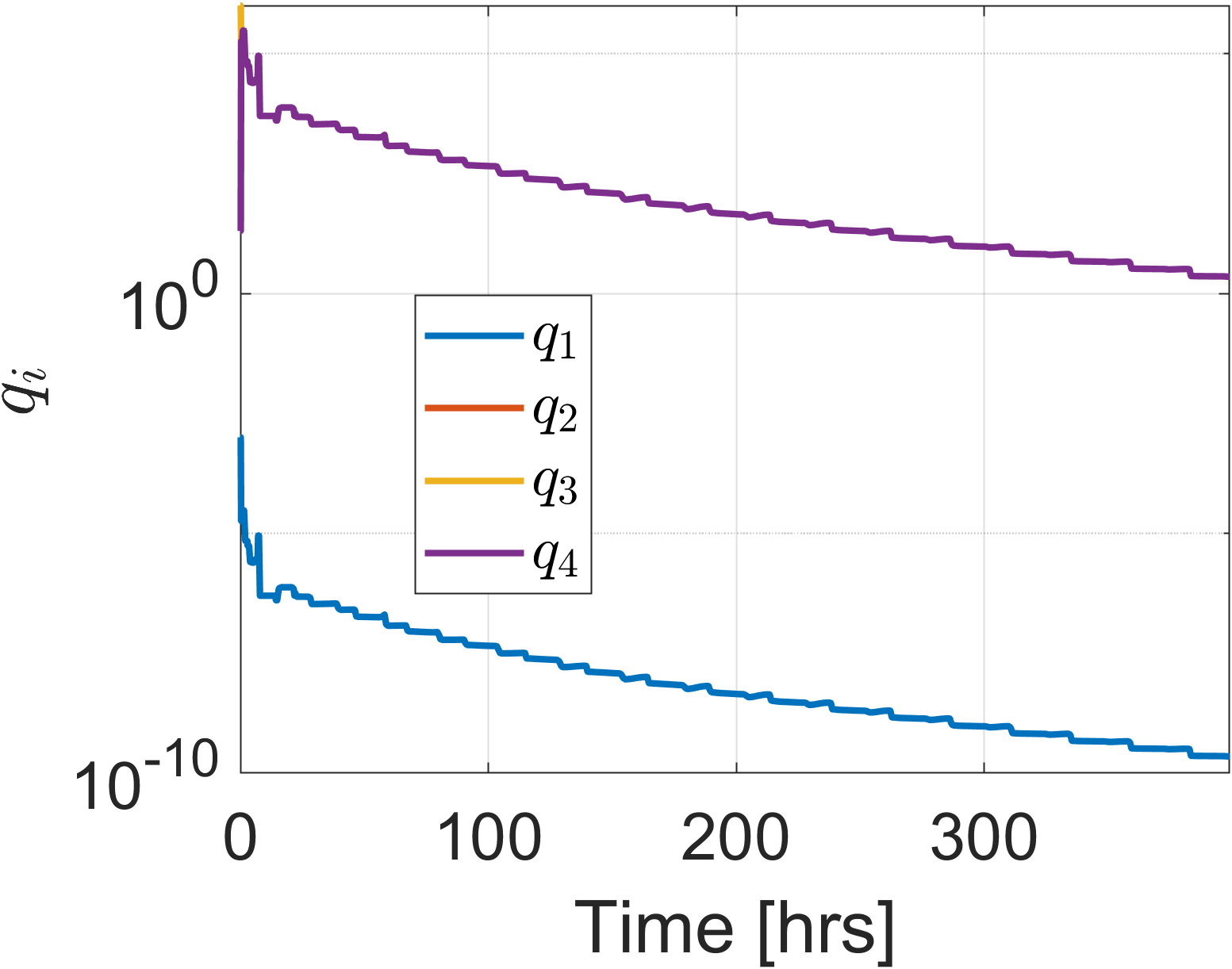}
	\caption{Time histories of reference governor parameter $\kappa_k$ for the initial 25 hours (left) and barrier function parameters $q_i$ (right) during transition from a lower semi-major axis orbit to a higher semi-major axis orbit with conventional thrust actuation.
		}
	\label{fig:kappa_up}
\end{figure}

Figure~\ref{fig:kappa_up} shows repeated instances of the gradual reduction in the reference governor parameter $\kappa_k$ till the terminal condition on the predicted trajectory is satisfied.  


\clearpage
\subsection{Orbital Maneuvers to Lower Orbit with Conventional Thrust Actuation}

\begin{figure}[h!]
	\centering
	\includegraphics[width=0.4\linewidth]{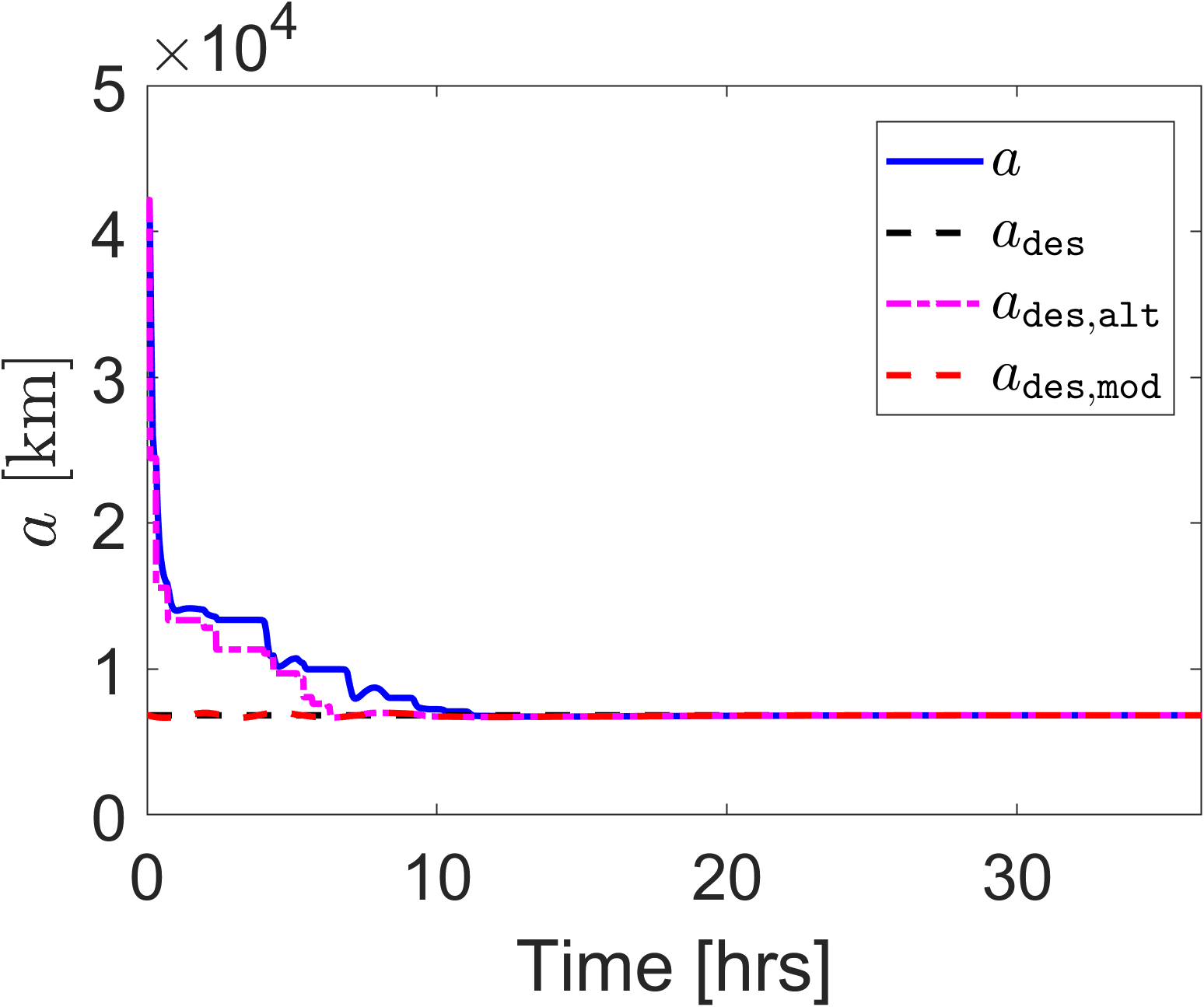}
	\includegraphics[width=0.4\linewidth]{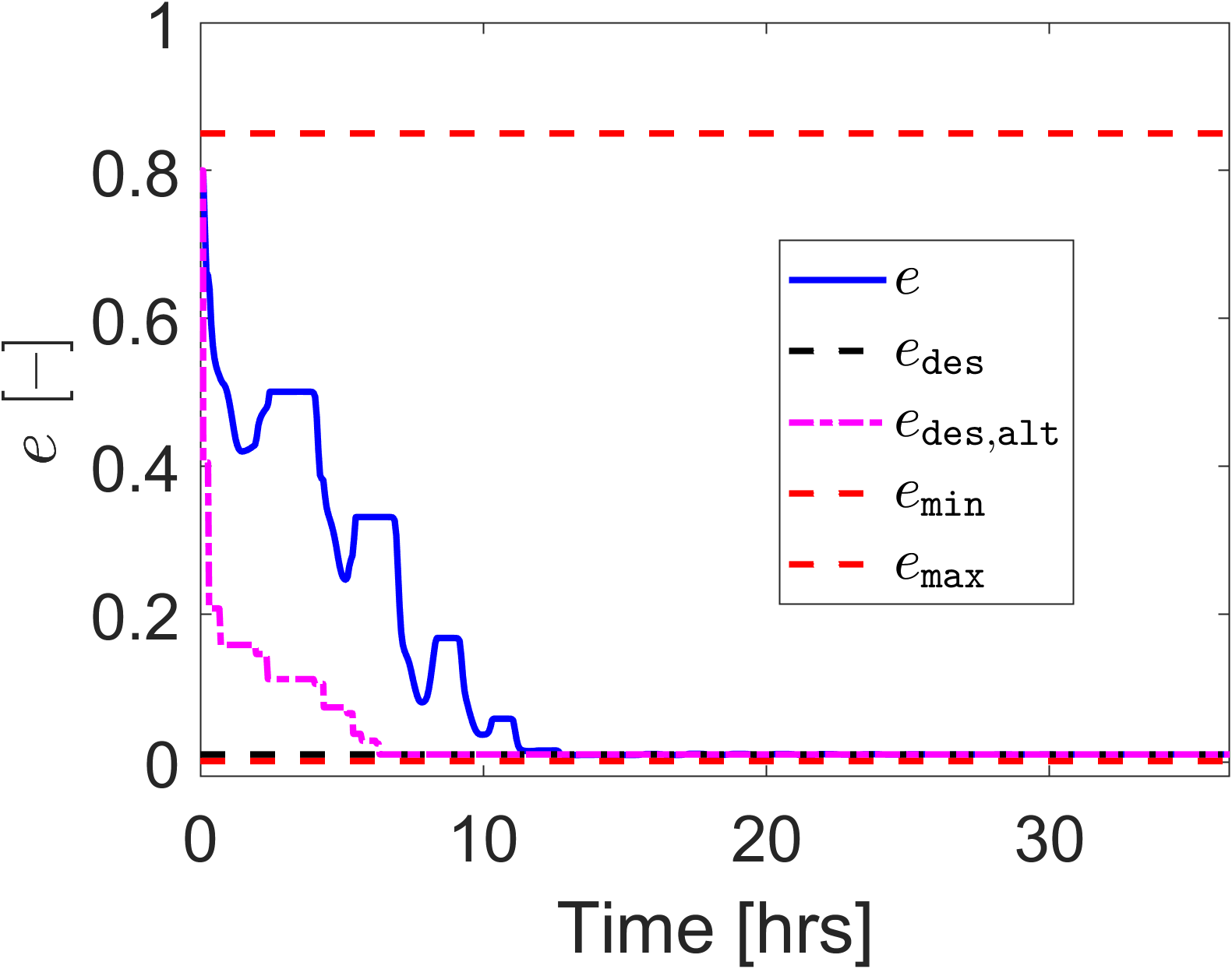}
	\caption{Semi-major axis (left) and eccentricity (right) during transition from a higher semi-major axis orbit to a lower semi-major axis orbit with conventional thrust actuation.}
	\label{fig:adown}
\end{figure}

\begin{figure}[h!]
	\centering
	\includegraphics[width=0.4\linewidth]{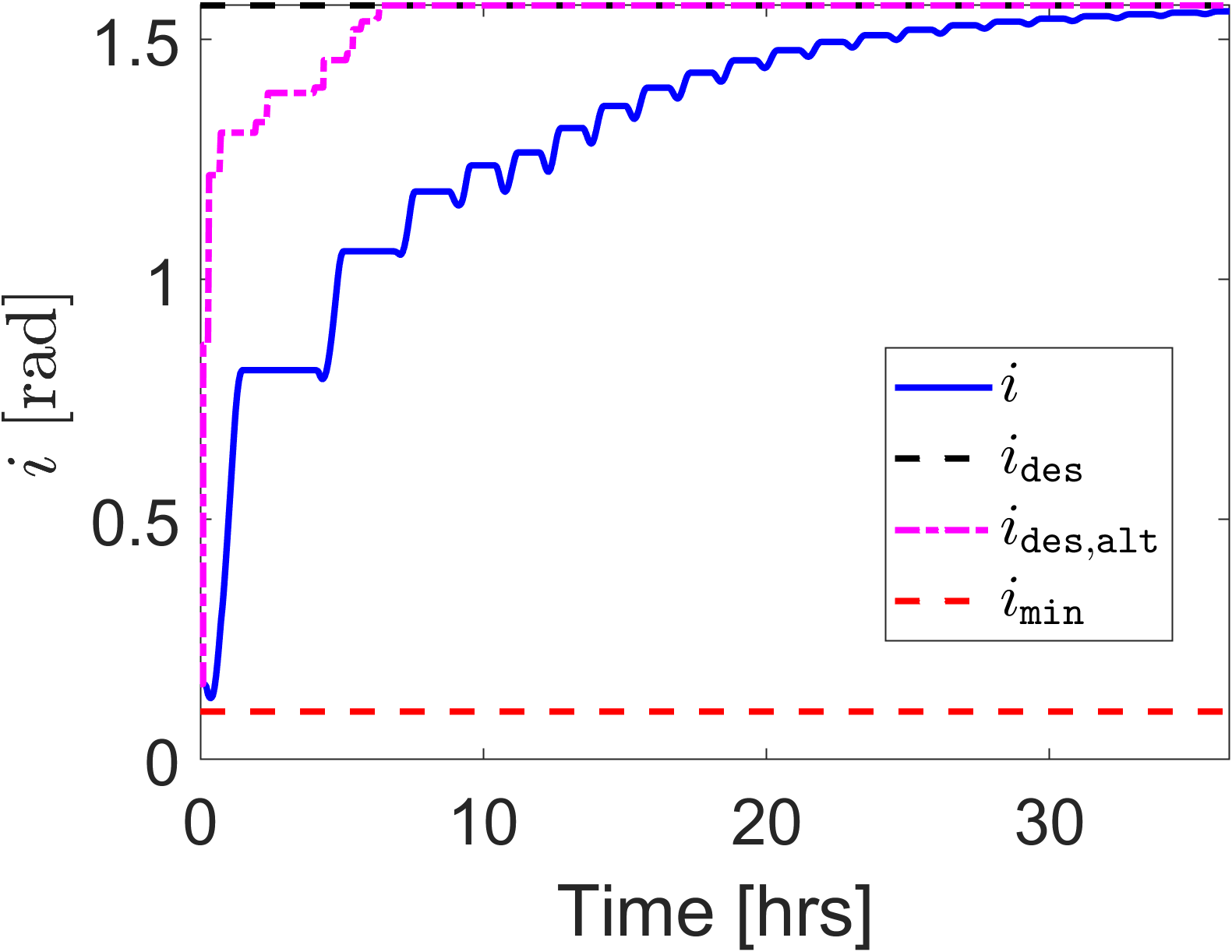}
	\includegraphics[width=0.4\linewidth]{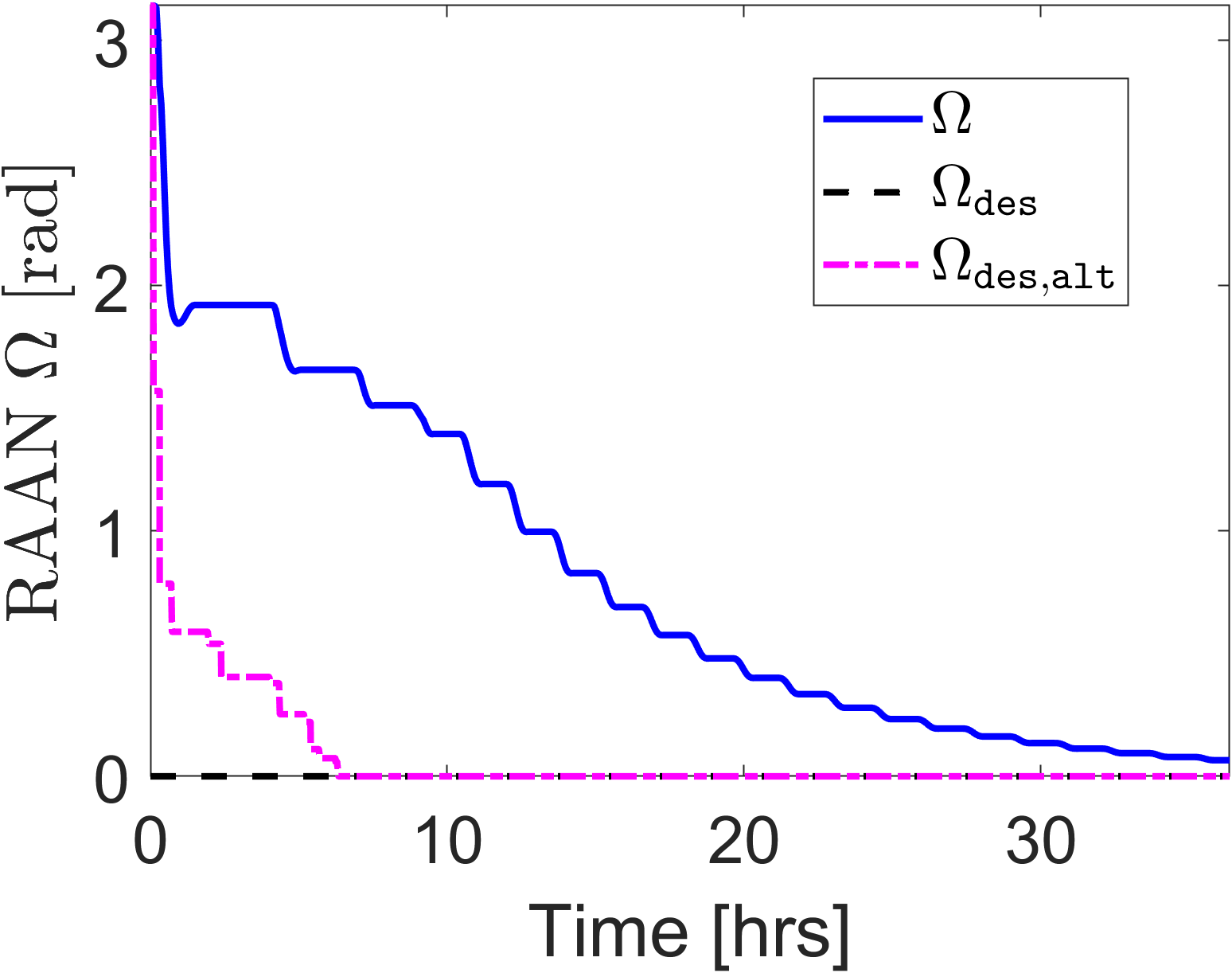}
	\caption{Inclination (left) and RAAN (right) during transition from a higher semi-major axis orbit to a lower semi-major axis orbit with conventional thrust actuation.}
	\label{fig:idown}
\end{figure}

\begin{figure}[h!]
	\centering
	\includegraphics[width=0.4\linewidth]{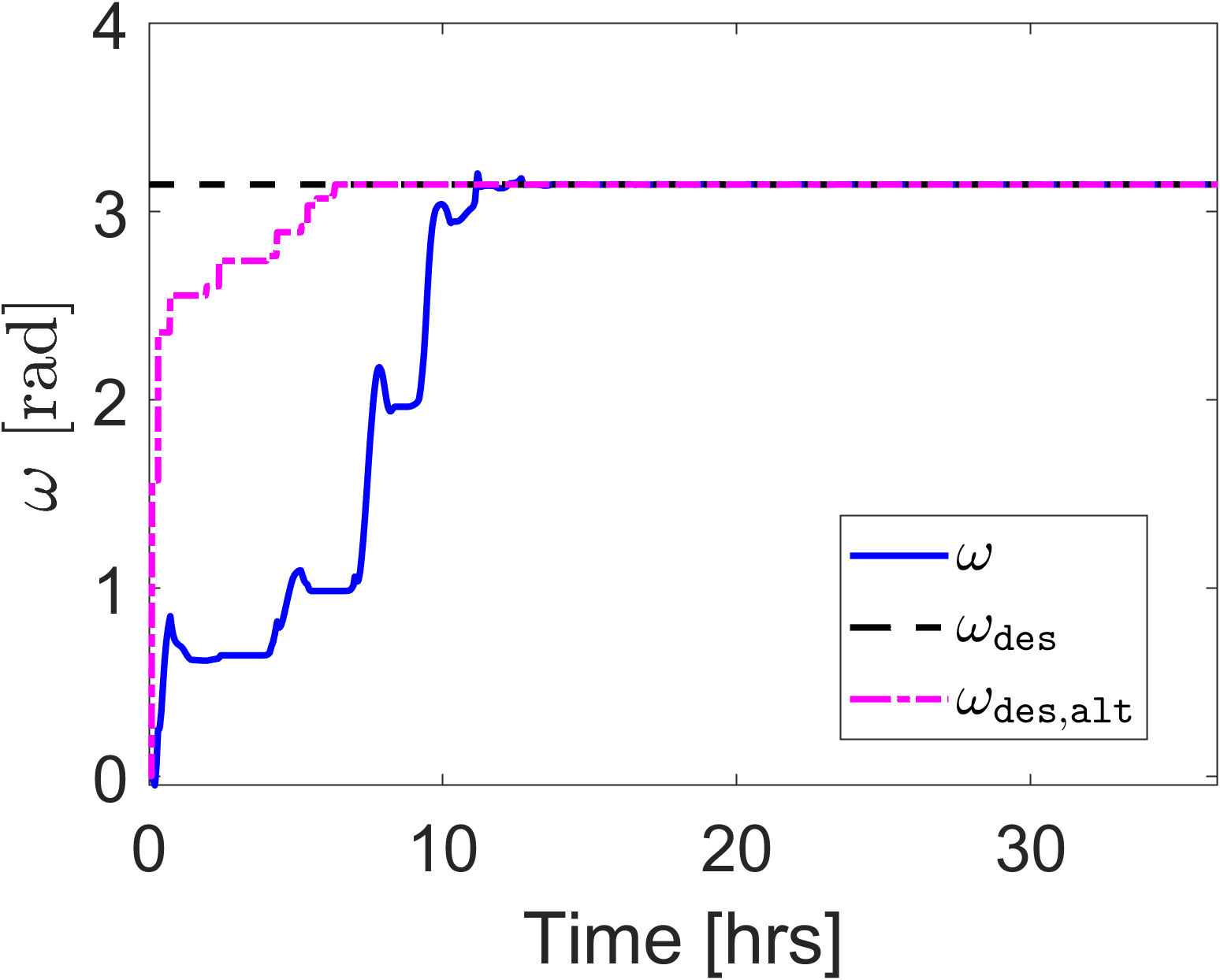}
	\includegraphics[width=0.4\linewidth]{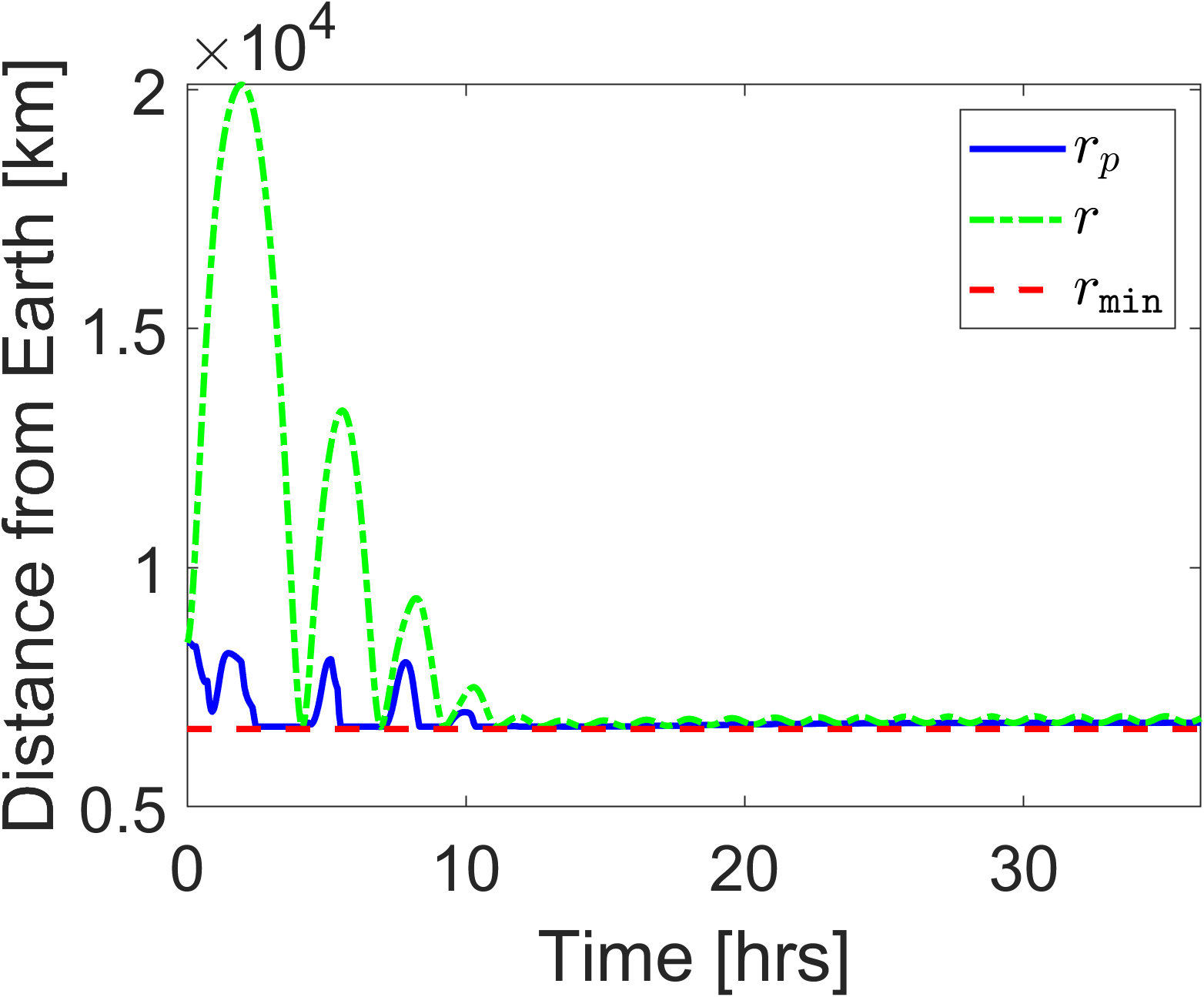}
	\caption{Argument of periapsis (left) and radius of periapsis (right) during transition from a higher semi-major axis orbit to a lower semi-major axis orbit with conventional thrust actuation.}
	\label{fig:argumentofperiapsisdown}
\end{figure}

\begin{figure}[h!]
	\centering
	\includegraphics[width=0.4\linewidth]{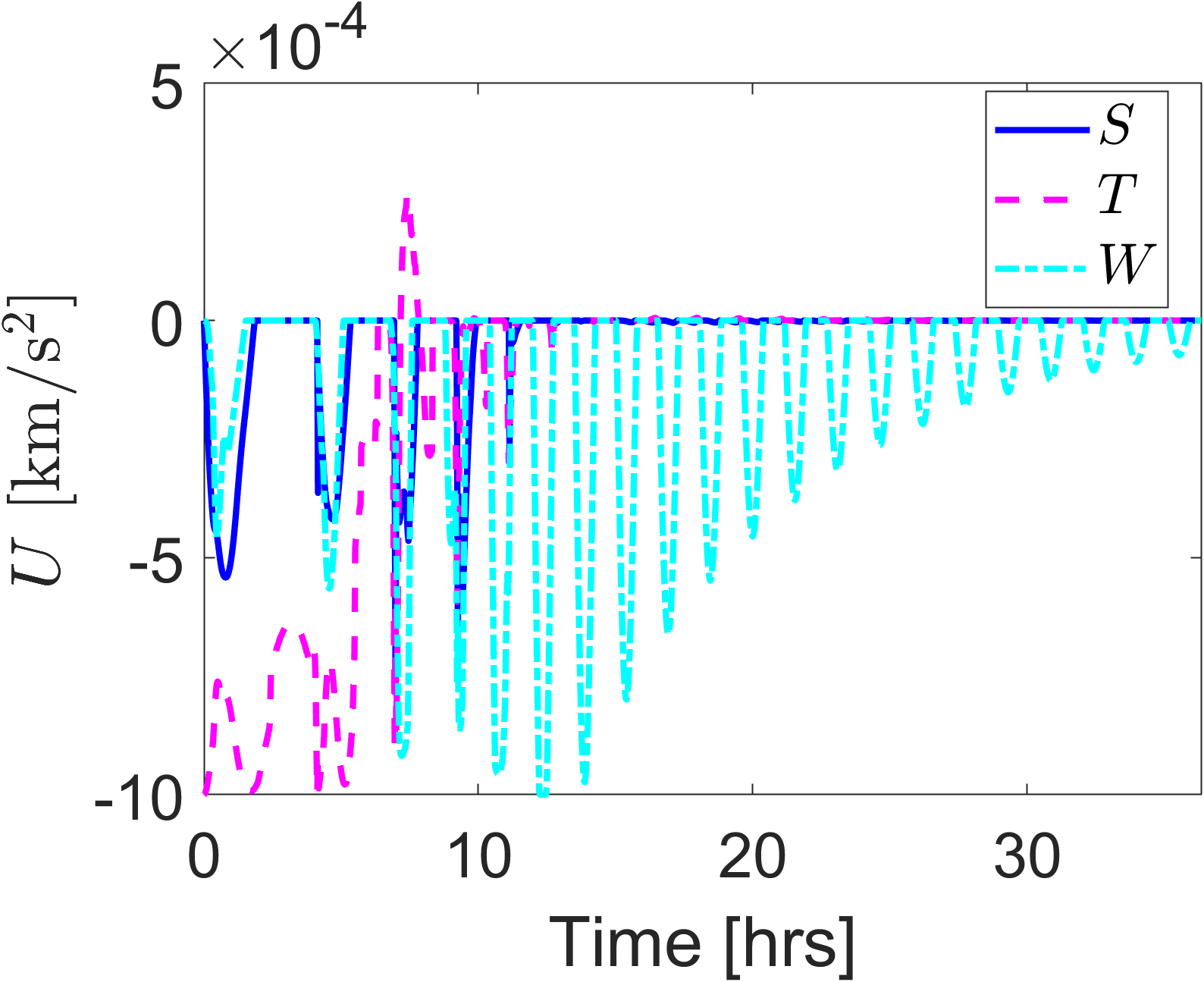}
	\includegraphics[width=0.4\linewidth]{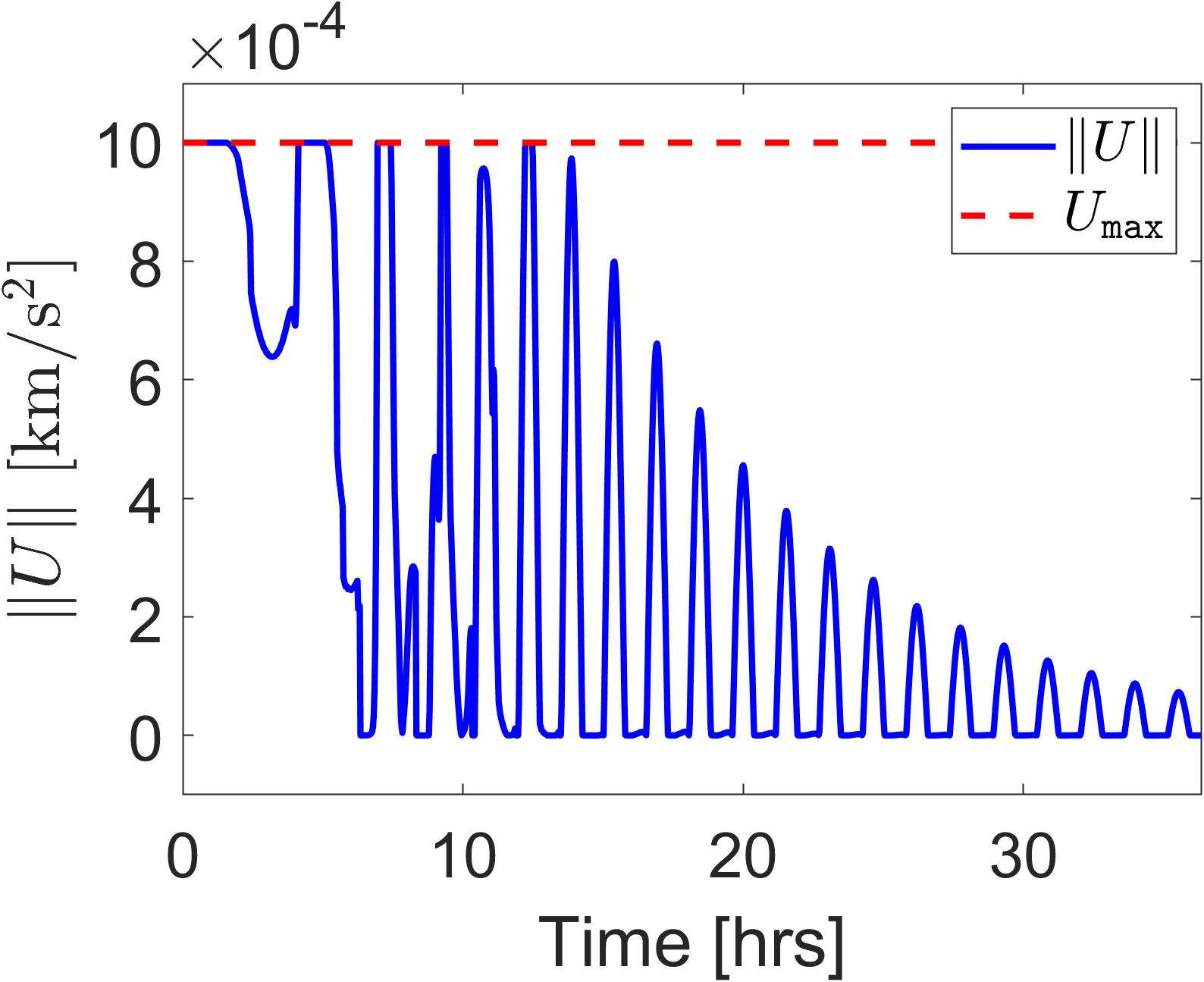}
	\caption{$S$, $T$ and $W$ components (left) and $\|U\|$ (right) during transition from a higher semi-major axis orbit to a lower semi-major axis orbit with conventional thrust actuation.}
	\label{fig:STWdown}
\end{figure}

\begin{figure}[h!]
	\centering
	\includegraphics[width=0.4\linewidth]{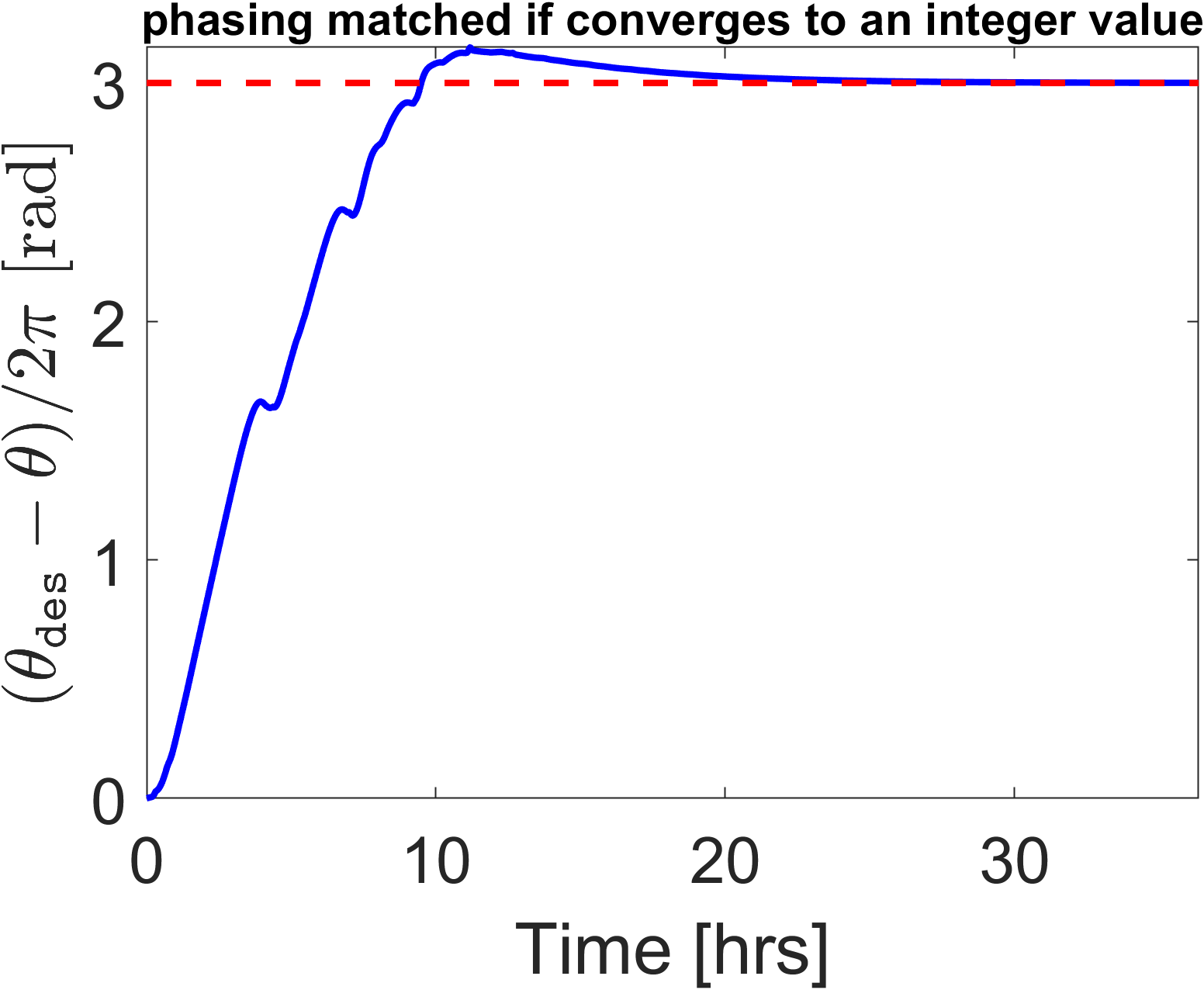}
    \includegraphics[width=0.4\linewidth]{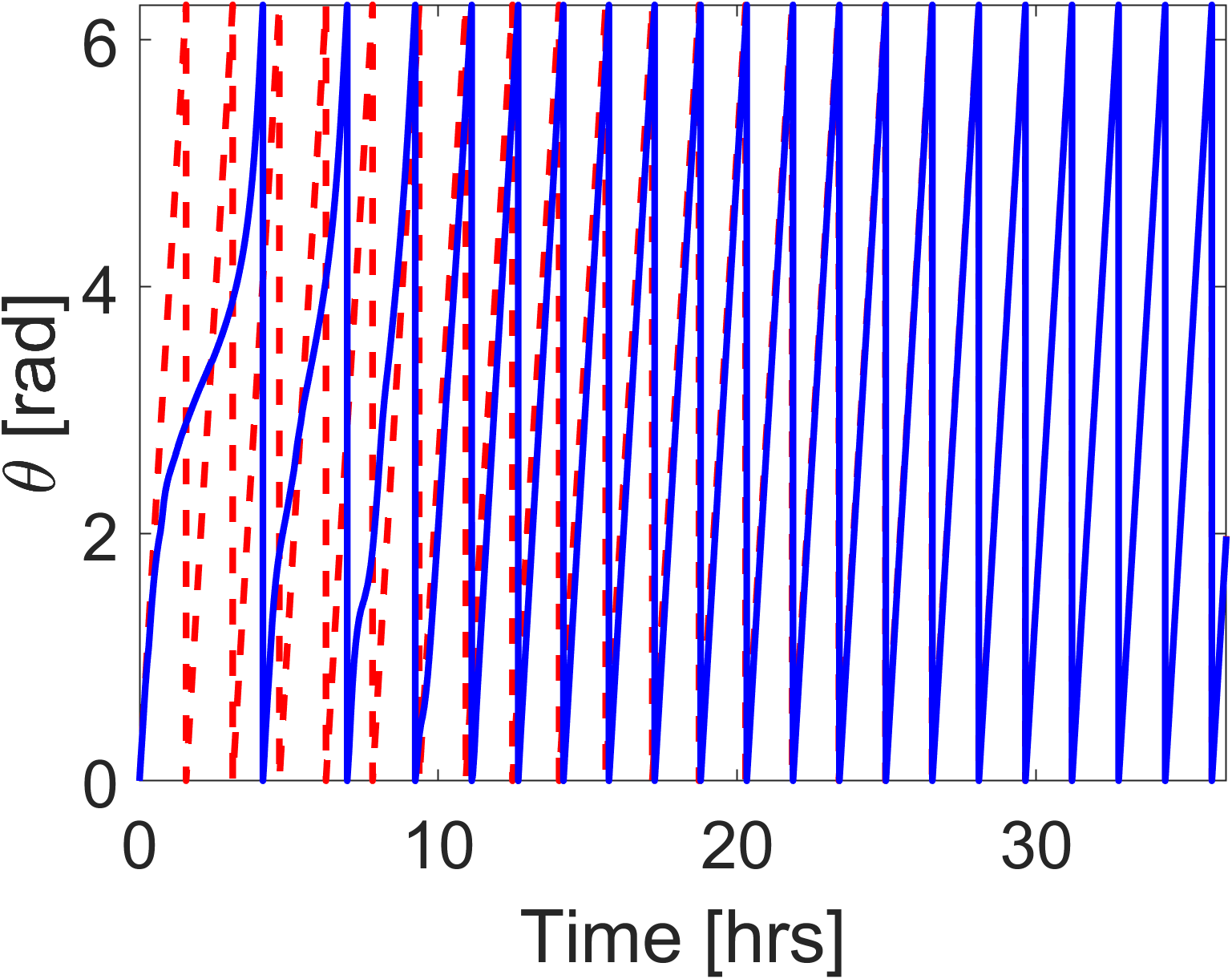}
	\caption{True anomaly difference (left) and true anomalies versus time (right) during transition from a higher semi-major axis orbit to a lower semi-major axis orbit with conventional thrust actuation.
		}
	\label{fig:thetadown}
\end{figure}

Figures~\ref{fig:adown}-\ref{fig:thetadown} show the transition from higher semi-major axis orbit to lower semi-major axis orbit with the conventional thrust actuation.  Note that the maneuver time is considerably shorter as compared to the opposite transfer while all the control constraints, $S \geq 0$ and $W \geq 0$ and $\|U \| \leq U_{\tt max}$,  and state constraints are enforced.  


Figure~\ref{fig:thetadown}-left shows the time history of $\frac{\theta_{\tt des}(t)-\theta(t)}{2 \pi}$, which converges to an integer value of $3$; this indicates the true anomaly matching by the end of the maneuver, achieved by altering the semi-major axis reference command indicated in Figure~\ref{fig:adown}-left by the red dashed line.

\begin{figure}[h!]
	\centering
	\includegraphics[width=0.4\linewidth]{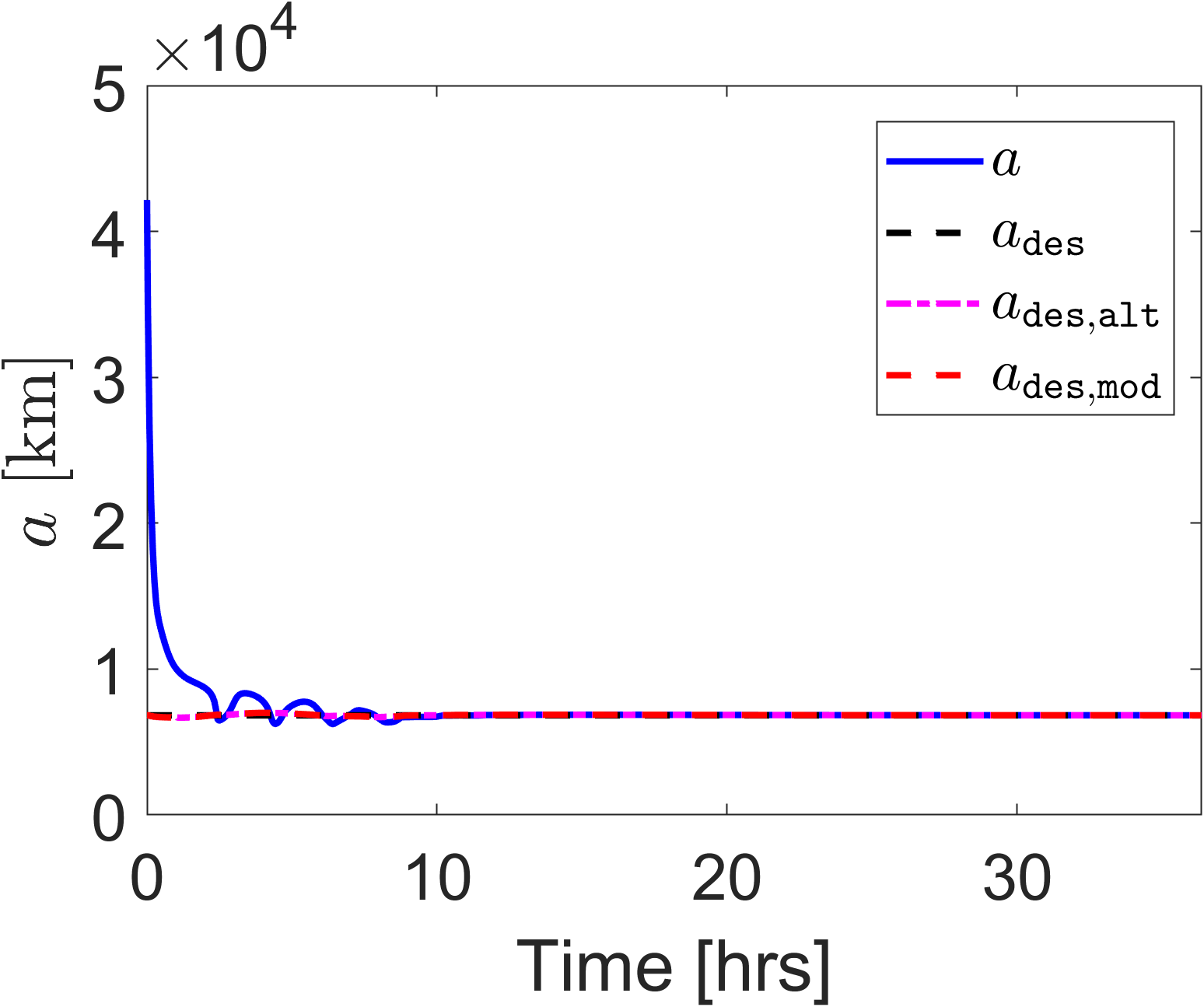}
	\includegraphics[width=0.4\linewidth]{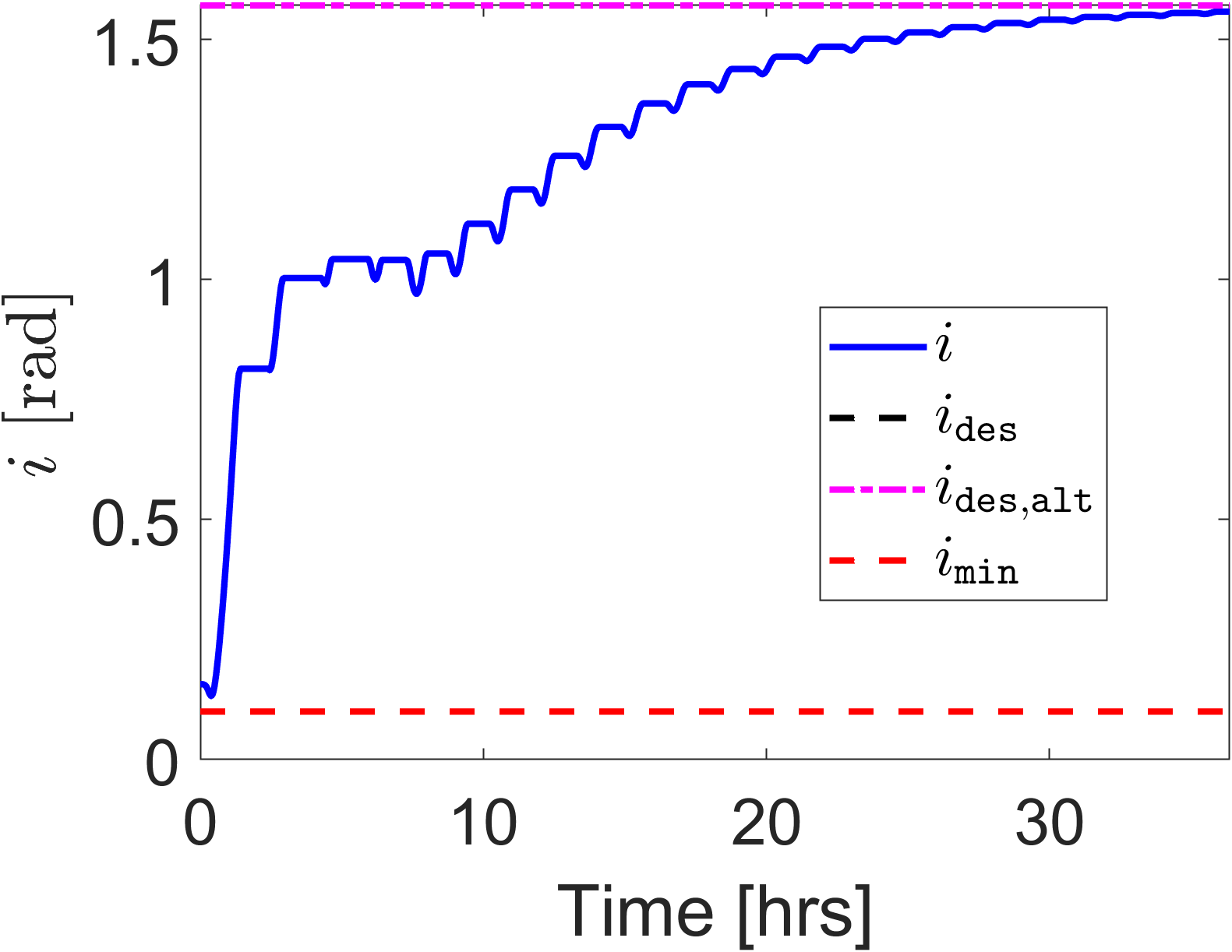}
	\caption{Semi-major axis (left) and inclination (right) during transition from higher semi-major axis to lower semi-major axis orbit with conventional thrust actuation without reference governor and without barrier functions.}
	\label{fig:ablation1}
\end{figure}

\begin{figure}[h!]
	\centering
	\includegraphics[width=0.4\linewidth]{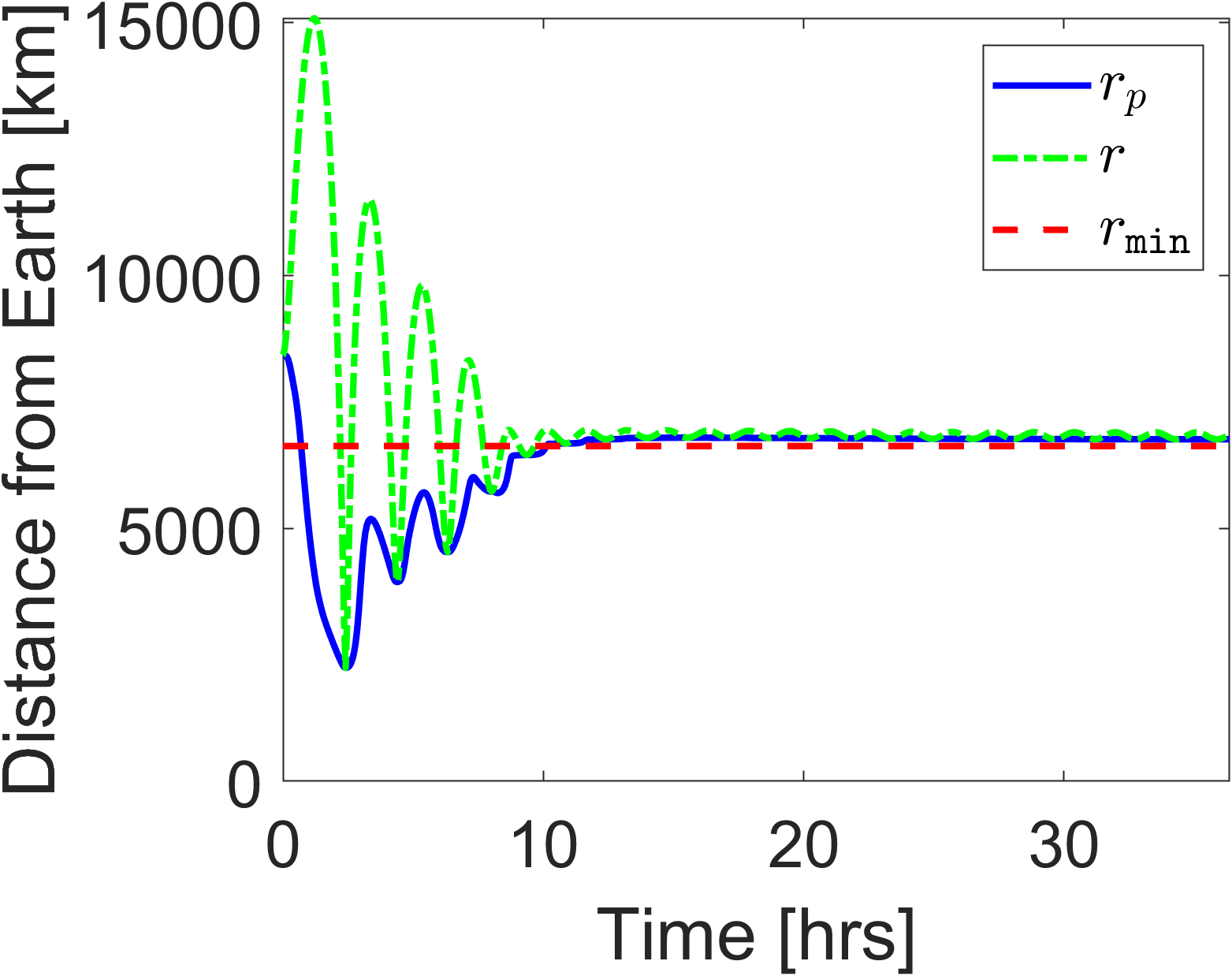}
	\includegraphics[width=0.4\linewidth]{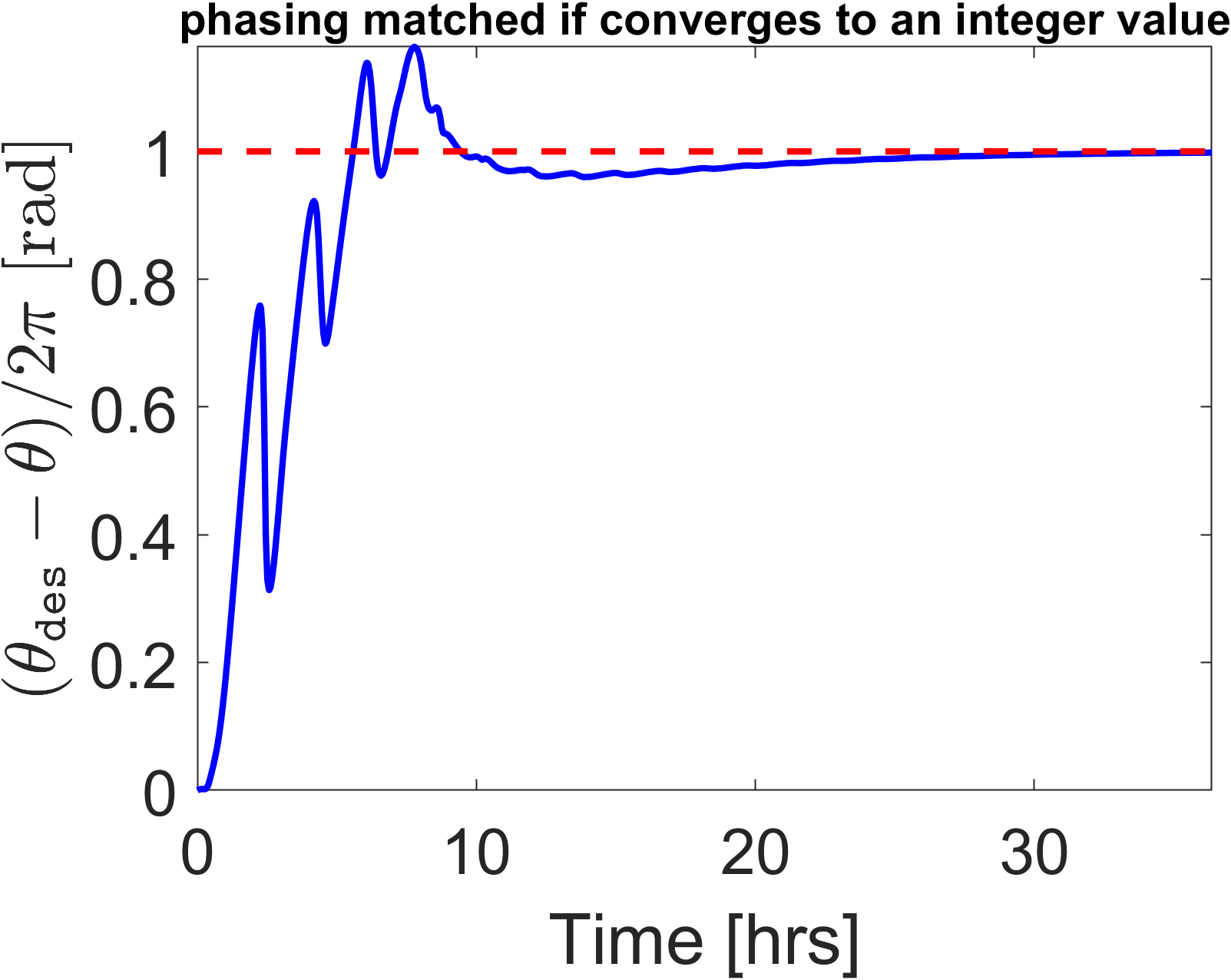}
	\caption{Radius of periapsis (left) and difference between actual and desired true anomalies (right) during transition from a higher semi-major axis to a lower semi-major axis orbit with conventional thrust actuation without reference governor and without barrier functions.}
	\label{fig:ablation2}
\end{figure}

Finally, Figures~\ref{fig:ablation1}-\ref{fig:ablation2} report simulation results without the reference governor and barrier functions.
Notably, the constraint on the radius of periapsis is violated by a large amount.

\begin{figure}[h!]
	\centering		\includegraphics[width=0.5\linewidth]{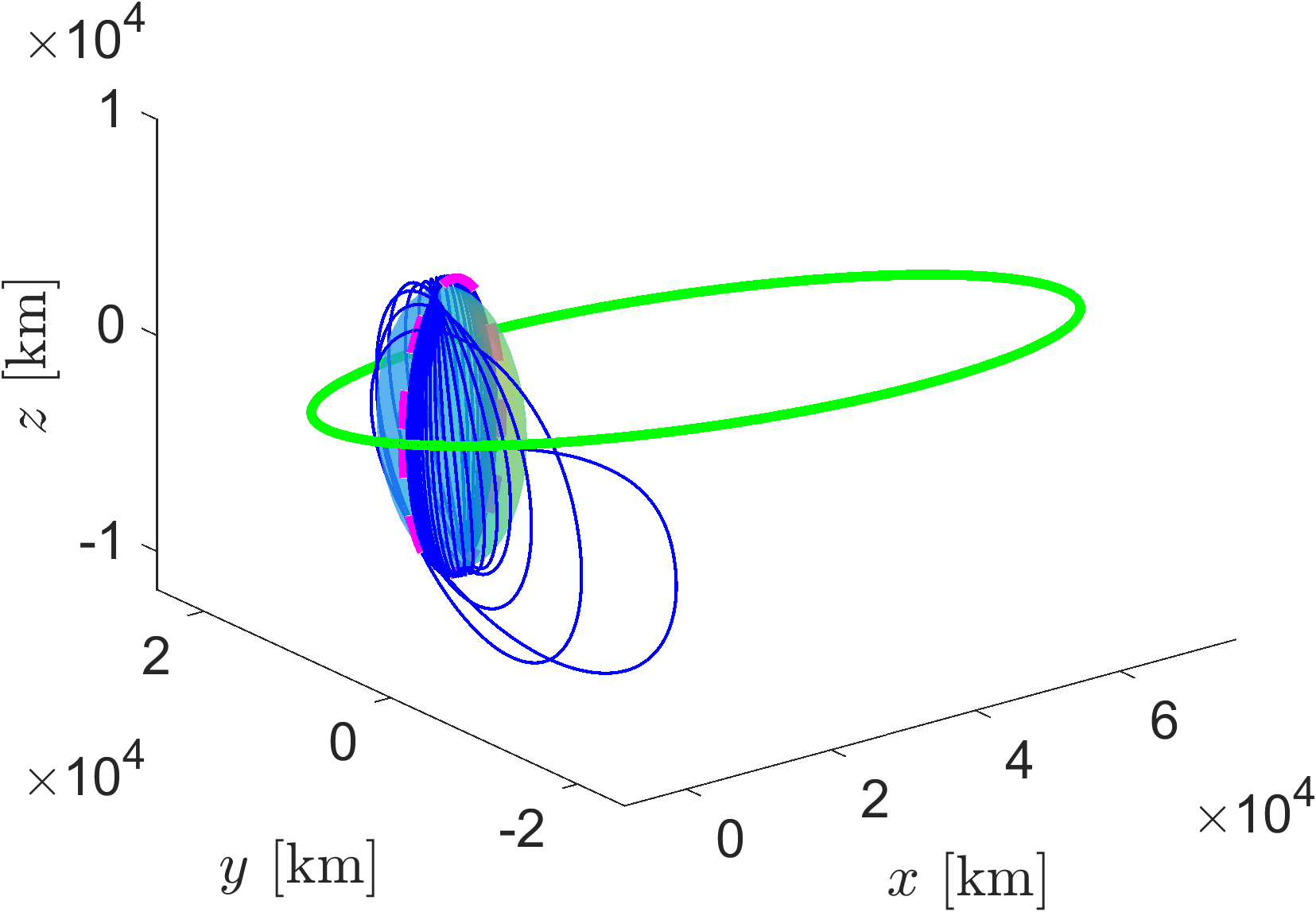}
	\caption{Three dimensional initial orbit (green), transfer orbit (blue) and final (desired) orbit (dashed magenta) during transition from a higher semi-major axis orbit to  a lower semi-major axis orbit with conventional thrust actuation. }
	\label{fig:thetadown_1}
\end{figure}

\begin{figure}[h!]
	\centering
	\includegraphics[width=0.4\linewidth]{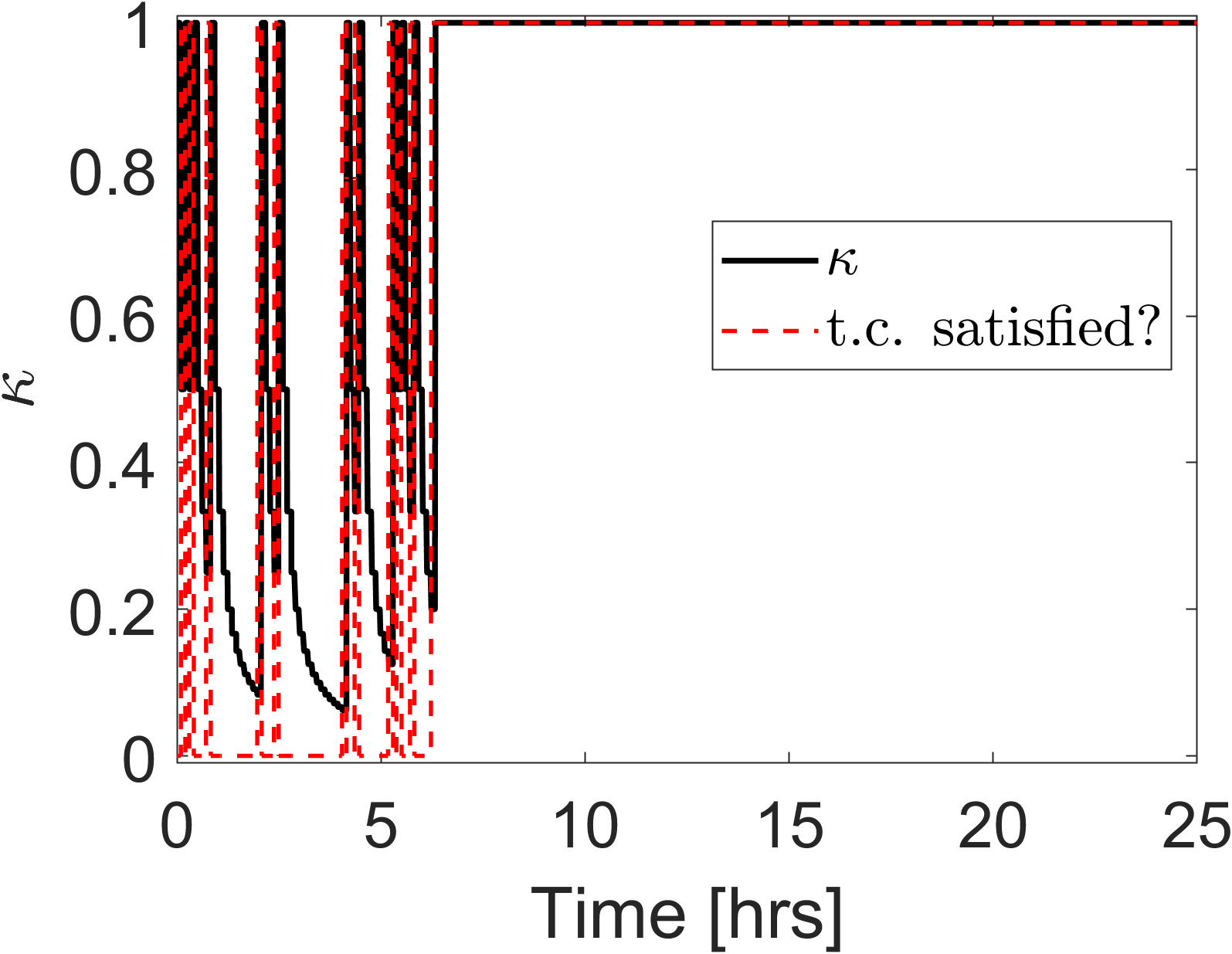}
\includegraphics[width=0.4\linewidth]{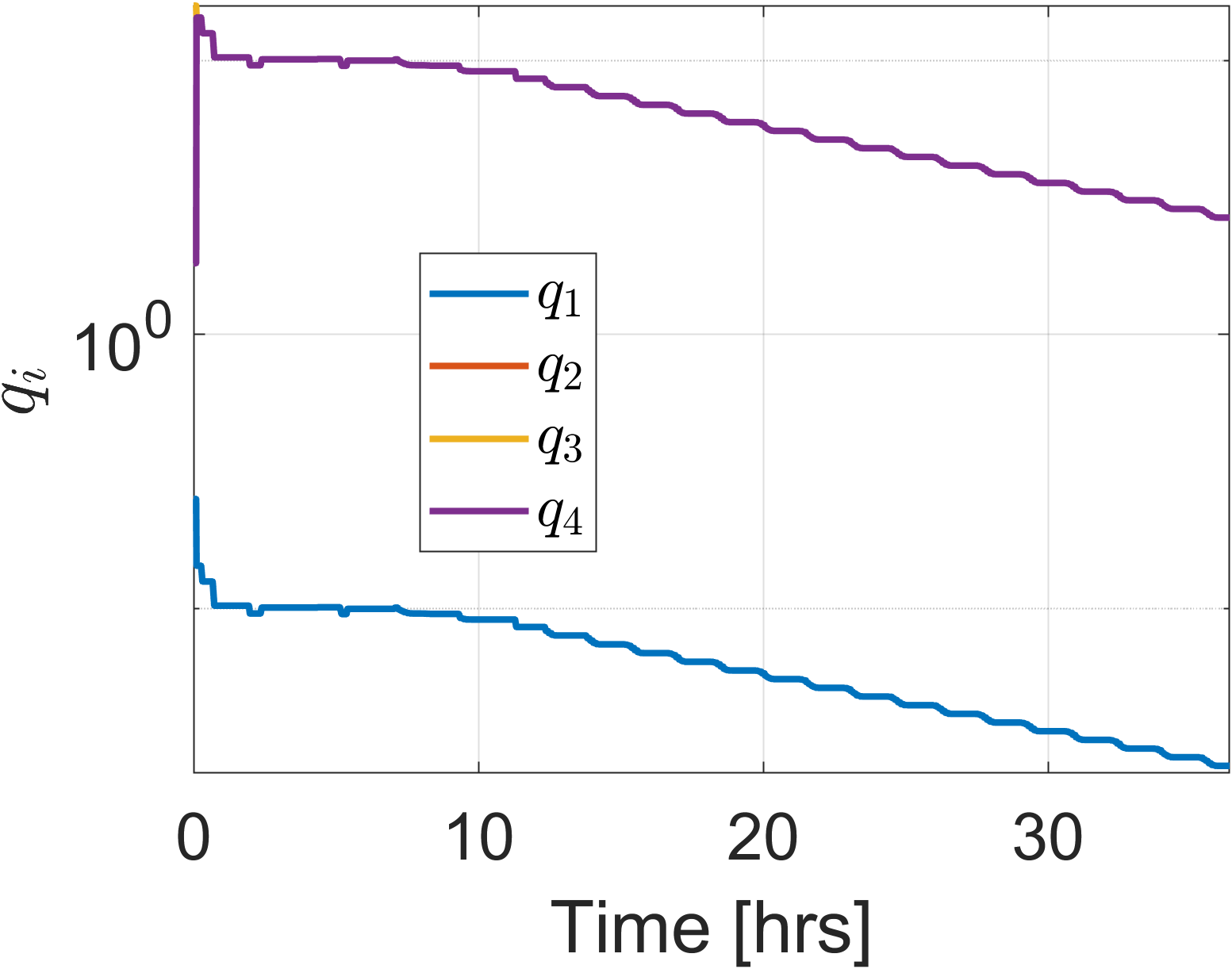}
	\caption{Time histories of reference governor parameter $\kappa_k$ for the initial 25 hours (left) and barrier function parameters $q_i$ (right) during transition from a lower semi-major axis orbit to a higher semi-major axis orbit with conventional thrust actuation. 
		}
	\label{fig:kappa_down}
\end{figure}

\section{Rendezvous Using Lorentz Force}

We approximate the  Lorentz force acting on the spacecraft with a tether as
$$\frac{\vec{F}}{m}=\frac{I}{m} \vec{l}\times \vec{B},~\vec{l}=l \hat{e}_r,$$
 where $I$ is the current through the tether and $l=0.5$ km. We assume that the tether remains stretched along the local vertical.  This assumption is consistent, e.g., with \cite{bombardelli2013deorbiting}. 

Electrodynamic tethers are considered
in \cite{levin2007dynamic} and in
\cite{cosmo1997tethers}.
 The table  on p.~184 of \cite{cosmo1997tethers}  highlights current patterns that
result in various orbital element changes one at a time. Here we demonstrate that we can change orbital elements 
simultaneously with a feedback law and under state and control constraints.

The magnetic field vector of the Earth, $\vec{B}$, is computed using a dipole-like model in which components of the Earth's magnetic field, $B_x$, $B_y$, $B_z$ (Tesla) in STW frame are expressed as functions of the distance to Earth center, $r$ (km), Earth's magnetic dipole value, $\mu_m=8.0 \times 10^6$ (Tesla$\cdot$km$^3$), inclination $i$ (rad), RAAN $\Omega$ (rad), Greenwich Mean Hour at $t=0$  GMH$_0$ (angle 
  in rad of Earth rotation around $z$-axis at $t=0$):
$$
B_x = -2 \mu_m  \frac{1}{r^3} ( \sin(\theta) \sin(i) \cos(\beta_m) + 
(\sin(\eta) cos(i) sin(\nu) + \cos(\eta) cos(\theta) ) \sin(\beta_m) ),
$$
$$
B_y = \mu_m \frac{1}{r^3}(\cos(\theta) \sin(i) \cos(\beta_m) +
( \sin(\eta) \cos(i) cos(\theta) - cos(\eta) sin(\theta) ) 
\sin(\beta_m) ),
$$
$$
B_z = \mu_m \frac{1}{r^3}
(\cos(i) \cos(\beta_m) - \sin(\eta) \sin(i) \sin(\beta_m) ),
$$
where
$\alpha_m = 4.4680429$ (rad) is  the right ascension of the magnetic axis, 
$\beta_m  = 0.2042035$ (rad) is  the dipole tilt angle,
$\theta_s = \mbox{GMH$_0$} +  2\pi t/86164$ (rad) is the Greenwich sidereal time and
$\eta =  \alpha_m - \Omega + \theta_s.
$

The GVEs with the Lorentz force actuation can be written in the form,
$$\dot{X}=\bar{G}(X,\theta)\bar{U},$$
where 
$$\bar{G}(X,\theta)=G(X,\theta) \left[\begin{array}{c} 0 \\ -\frac{l B_z}{m} \\ \frac{l B_y}{m} \end{array} \right],$$ and $\bar{U}=I\in \bar{\mathcal{U}}=[-I_{\tt max},I_{\tt max}]$ is the current through the tether.  Notably, only the current through the tether can be controlled. As the orientation of the tether is maintained along the local vertical, the radially induced acceleration is zero, $S(t)=0$ for all $t$.  Thus the spacecraft dynamics are underactuated: Only a single current input is available.

The following feedback law is adopted:
\[\bar{U}=-\mathrm{Proj}_{\bar{\mathcal{U}}} \left[ \frac{5 \times 10^5}{\frac{l}{m} \sqrt{B_y^2+B_z^2}} \bar{G}^{\sf T}({X},\theta) {P} ({X}-{X}_{\tt des})  \right],
\]
where a scaling factor is applied as without this scaling the values of the current that are generated are too low (alternatively, this scaling could be absorbed into the elements of $P$). 
This feedback law was augmented by barrier function terms and by the outer-loop controller for adjusting the reference command for the semi-major axis.

After a manual tuning we chose, $P=diag(6.4 \times 10^{-8}, 2, 2, 1.5, 0.01),$ and $\bar{P}=diag(0.05,0.05)$.  
Generally, the maneuvers with only Lorentz force actuation take much longer and are more limited in terms of their reach as opposed to the use of regular propulsion.  The reference governor was not used in the subsequent simulations due to the limited range of the maneuvers.

In the sequel, we simulate orbital maneuvers from a lower semi-major axis orbit to a higher semi-major axis orbit and from a higher semi-major axis orbit to a lower semi-major axis orbit.  In the former case,
$X(0)=[R_e+450, 0.02, 57^\circ, 18^\circ, 20^\circ]^{\sf T}$, 
$X_{\tt des}=[R_e+850, 0.06, 64^\circ, 22^\circ, 30^\circ]^{\sf T},$ and in the latter case,
$X(0)=[R_e+850, 0.06, 64^\circ, 22^\circ, 30^\circ]^{\sf T},$
$X_{\tt des}=[R_e+450, 0.02, 57^\circ, 18^\circ, 20^\circ]^{\sf T}$.
The desired true anomaly trajectory is computed based on
$$\dot{\theta}_{\tt des}(t)=\frac{\sqrt{\mu} (1+e_{\tt des} \cos(\theta_{des}(t)))^2}{
(a_{\tt des} (1-e_{\tt des}^2))^{\frac{3}{2}}}, \quad {\theta}_{\tt des}(0)=0.$$
Other parameters used in the simulations were $\mu=398600.4405$  (km$^3$/s$^2$), $R_e=6378$ (km), $r_{\tt min}=R_e+250$ (km),
$e_{\tt min}=1 \times 10^{-3}$, 
$e_{\tt max}=0.85$, $i_{\tt min}=0.1$,
$\epsilon_1=50$, $\epsilon_2=5 \times 10^{-4}$, $\epsilon_3=5 \times 10^{-4}$,
$\epsilon_4=5 \times 10^{-4}$, $I_{\tt max} = 10$ (A).
Tighter saturation limits are used versus (\ref{equ:amod}) when computing 
$a_{\tt {des},mod}$: 
\begin{equation}\label{equ:amod_lorentz}
a_{\tt {des},mod}
=
\max\!\Big(
0.9 a_{\tt des},\;
\min\!\big(
1.1 a_{\tt des},\;
\left(\frac{\mu}{{n}_{\tt mod}^{2}}\right)^{\!\tfrac{1}{3}}
\big)
\Big).
\end{equation}

\subsection{Orbital Maneuvers to Higher Orbit with Lorentz Force Based Propulsion}

\begin{figure}[h!]
	\centering
	\includegraphics[width=0.4\linewidth]{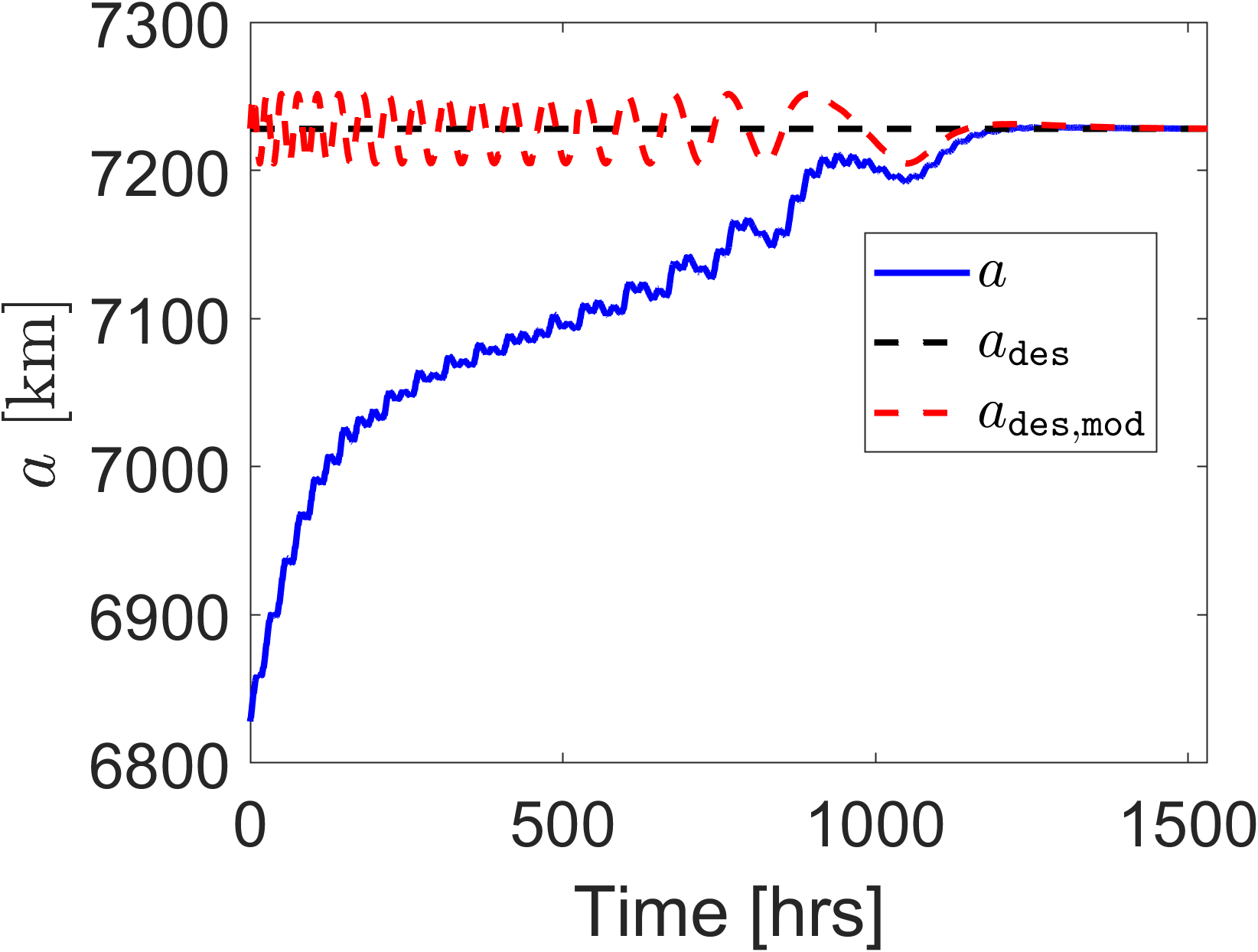}
	\includegraphics[width=0.4\linewidth]{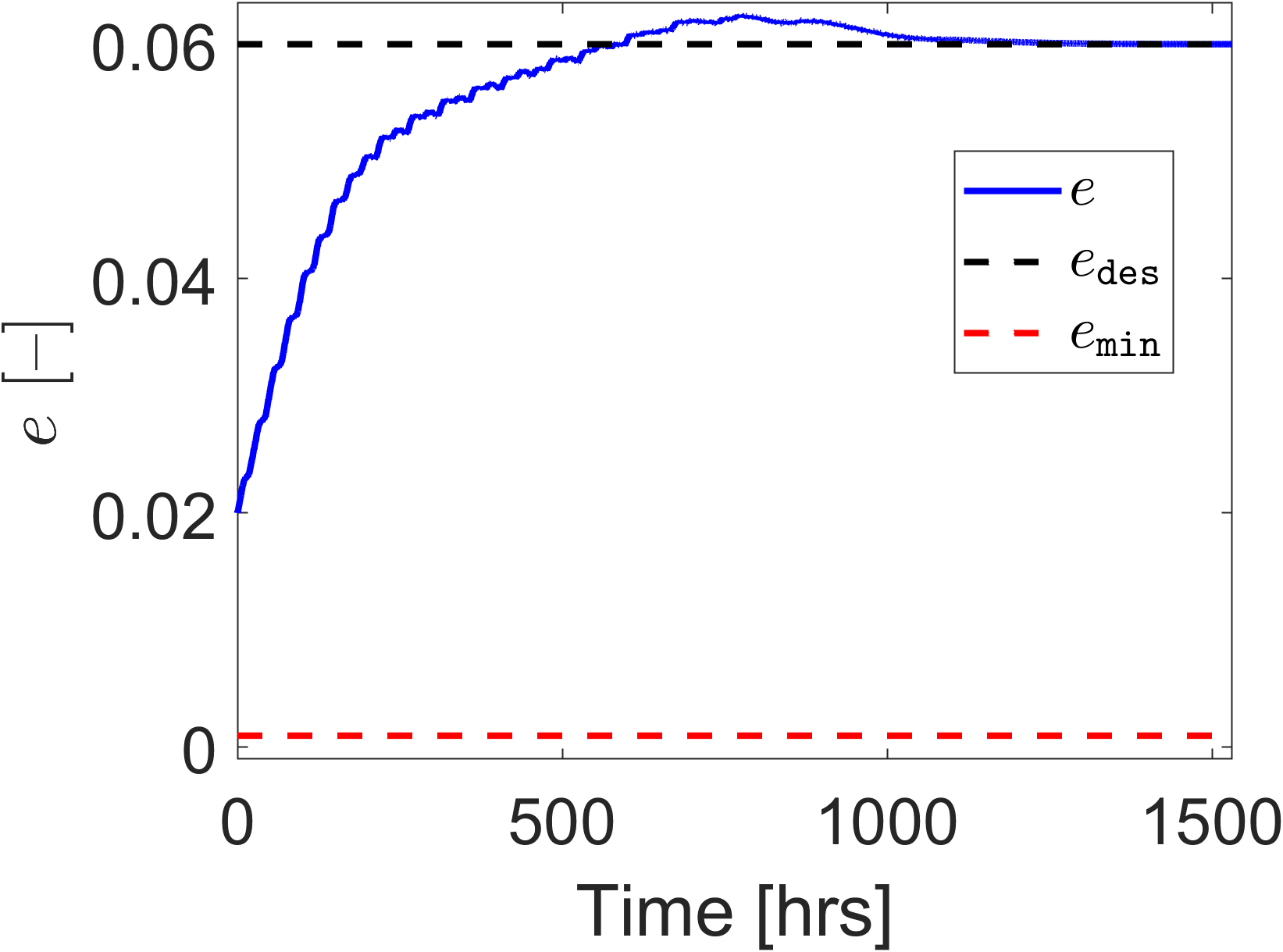}
	\caption{Semi-major axis (left) and eccentricity (right) during transition from a lower semi-major axis orbit to a higher semi-major axis orbit with Lorentz force actuation.}
	\label{fig:alorentzup}
\end{figure}

\begin{figure}[h!]
 	\centering
 	\includegraphics[width=0.4\linewidth]{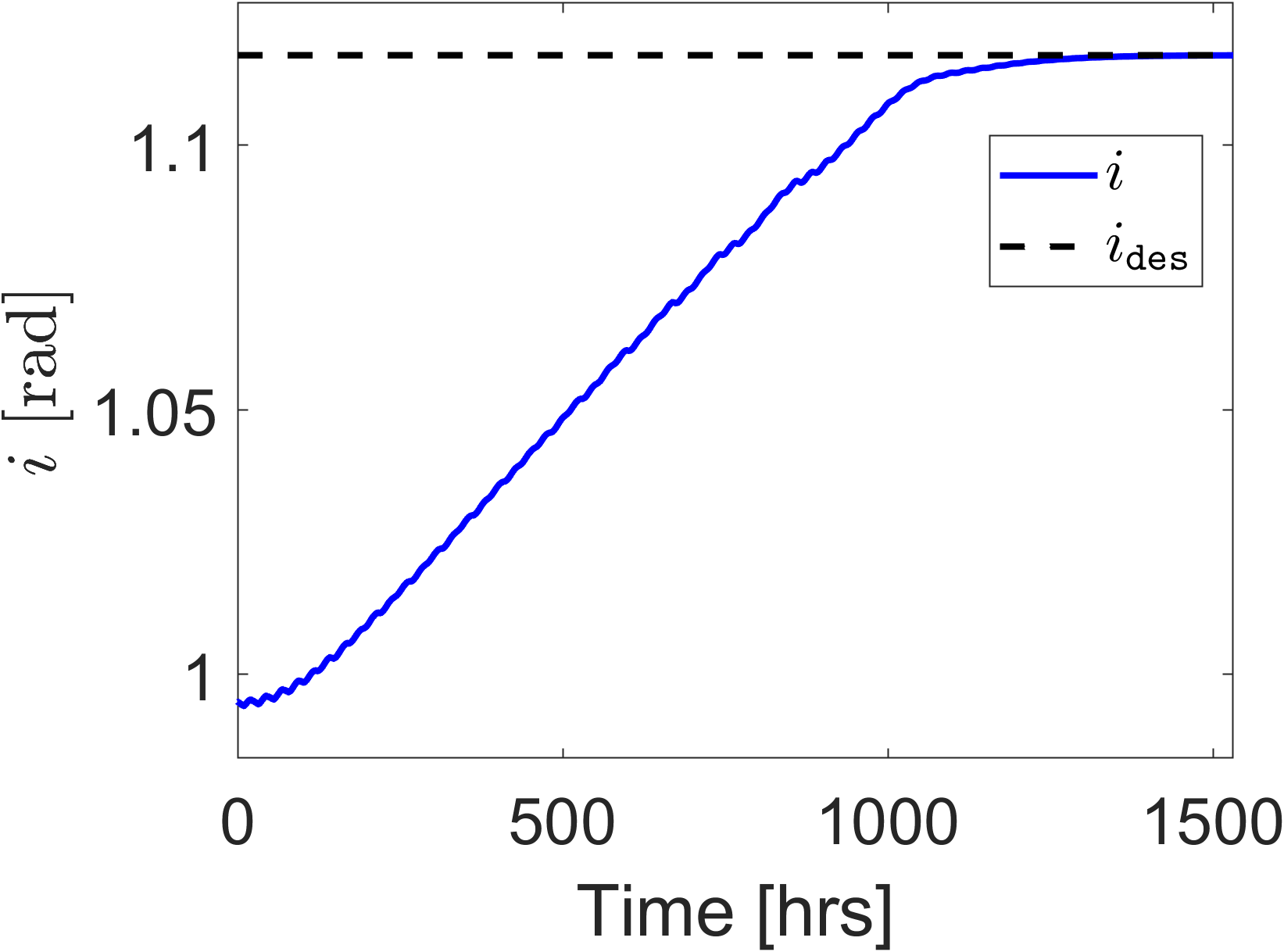}
 	\includegraphics[width=0.4\linewidth]{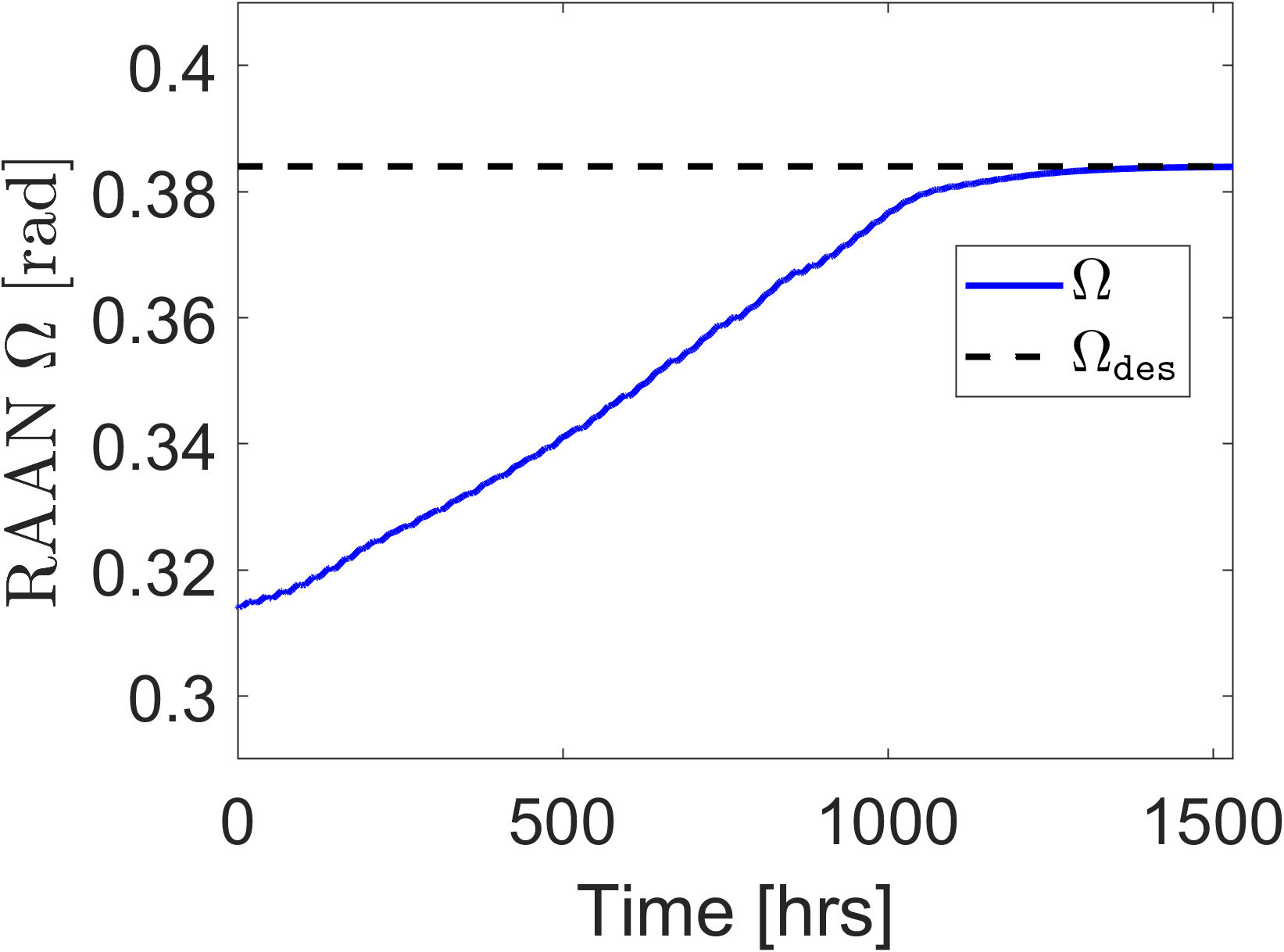}
 	\caption{Inclination (left) and RAAN (right) during transition from a lower semi-major axis orbit to a higher semi-major axis orbit with Lorentz force actuation.}
 	\label{fig:ilorentzup}
 \end{figure}

\begin{figure}[h!]
	\centering
	\includegraphics[width=0.4\linewidth]{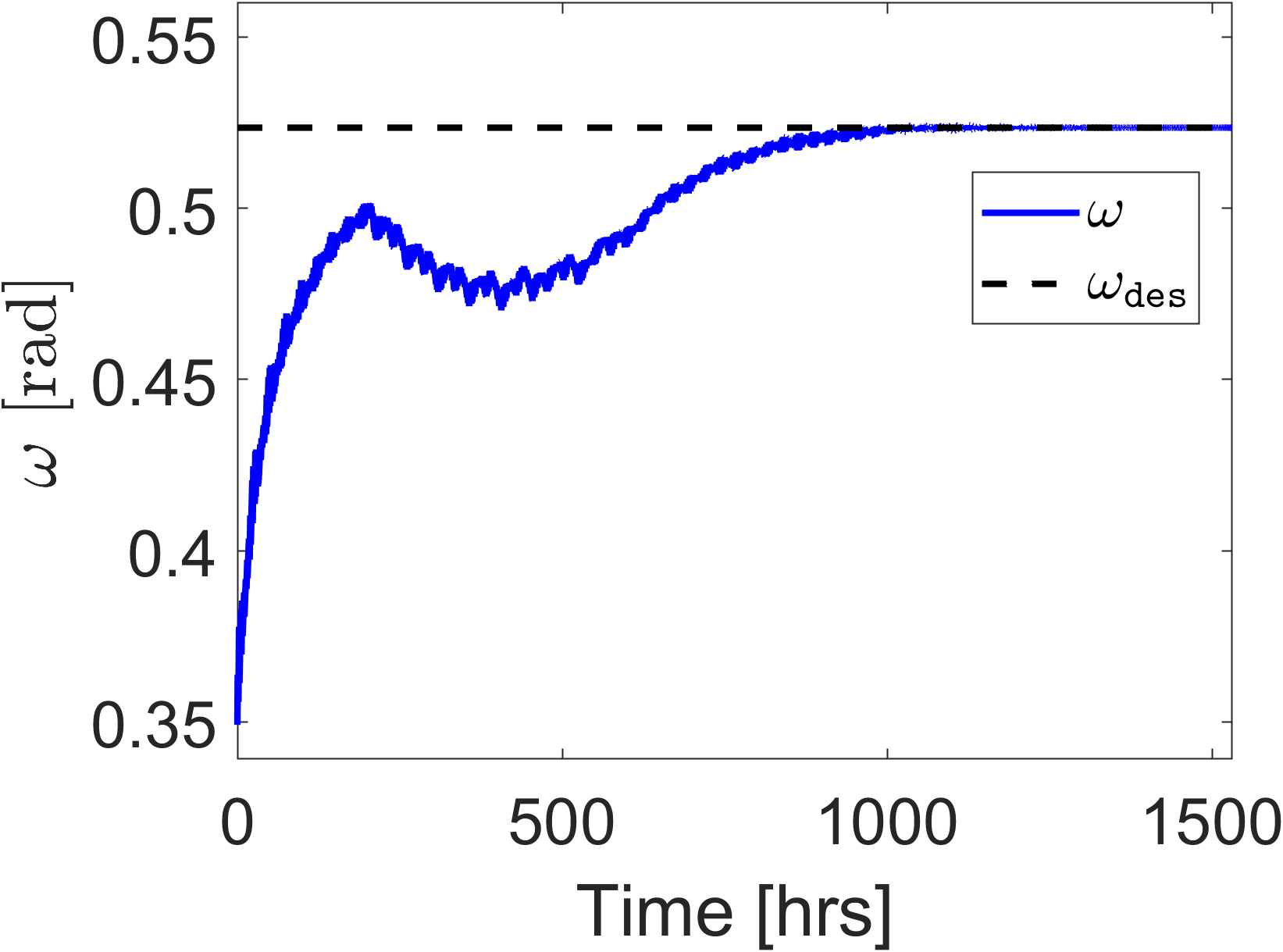}
	\includegraphics[width=0.4\linewidth]{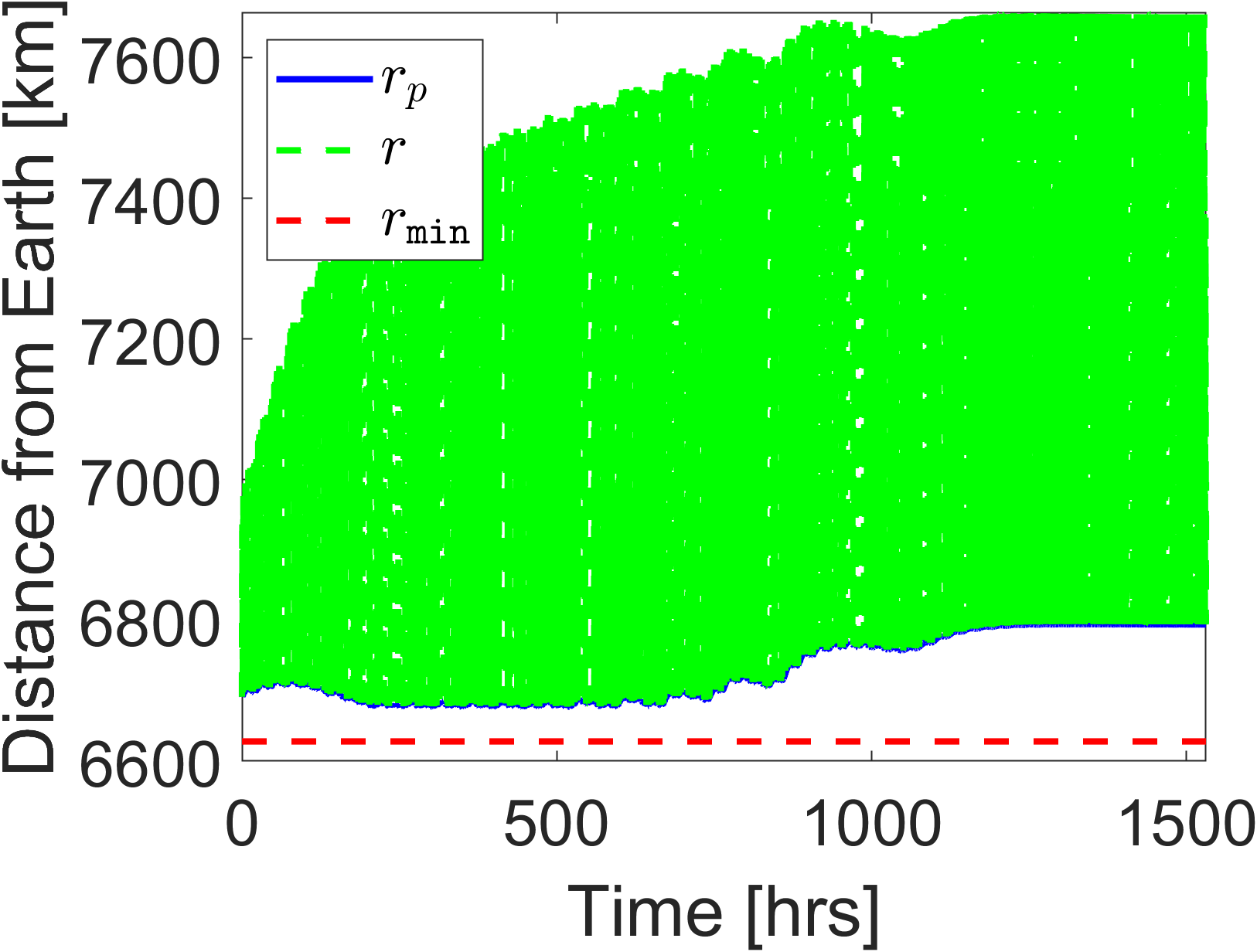}
	\caption{Argument of periapsis (left) and radius of periapsis (right) during transition from a lower semi-major axis orbit to a higher semi-major axis orbit with Lorentz force actuation.}
	\label{fig:argumentofperiapsislorentzup}
\end{figure}

\begin{figure}[h!]
	\centering
	\includegraphics[width=0.4\linewidth]{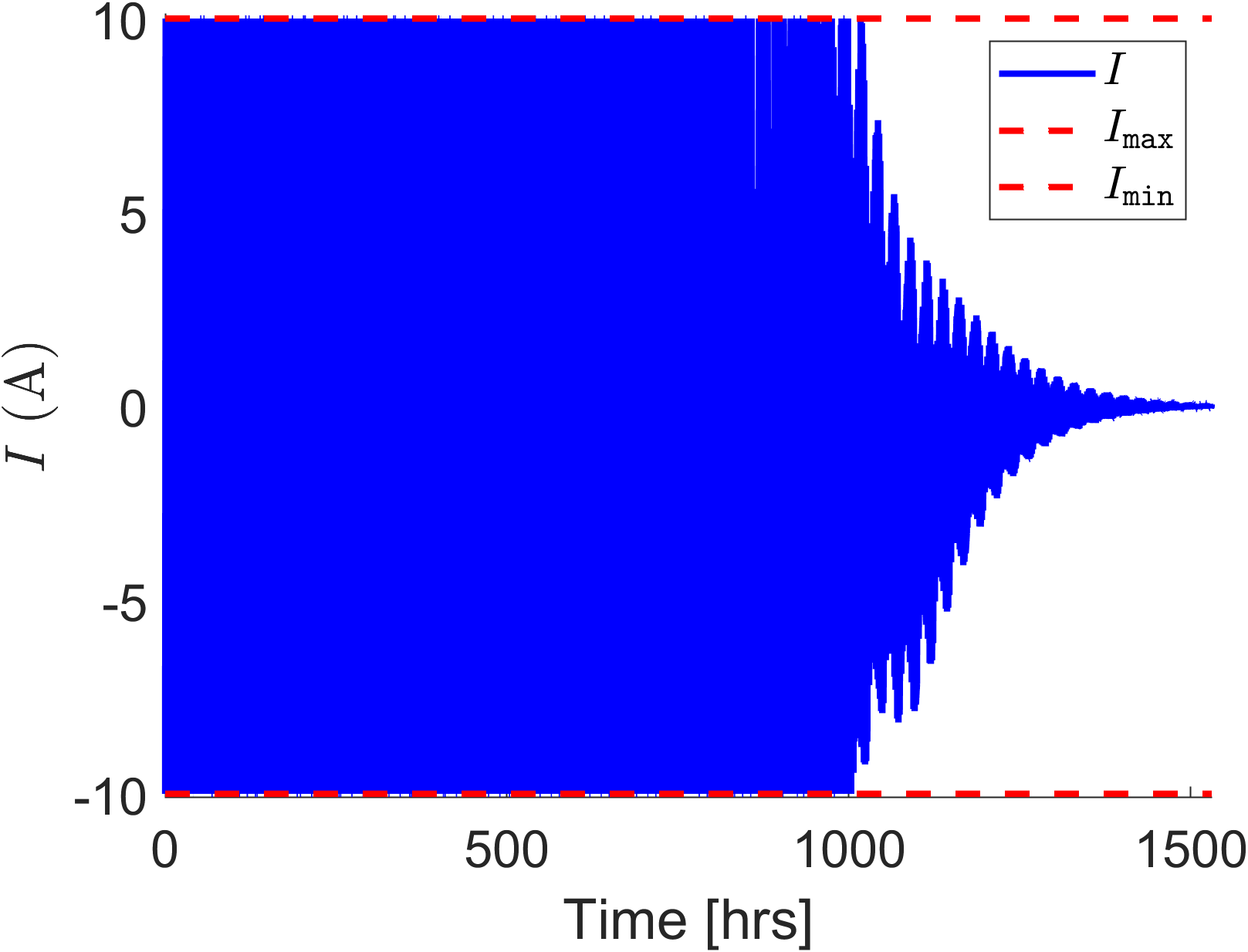}
	\includegraphics[width=0.4\linewidth]{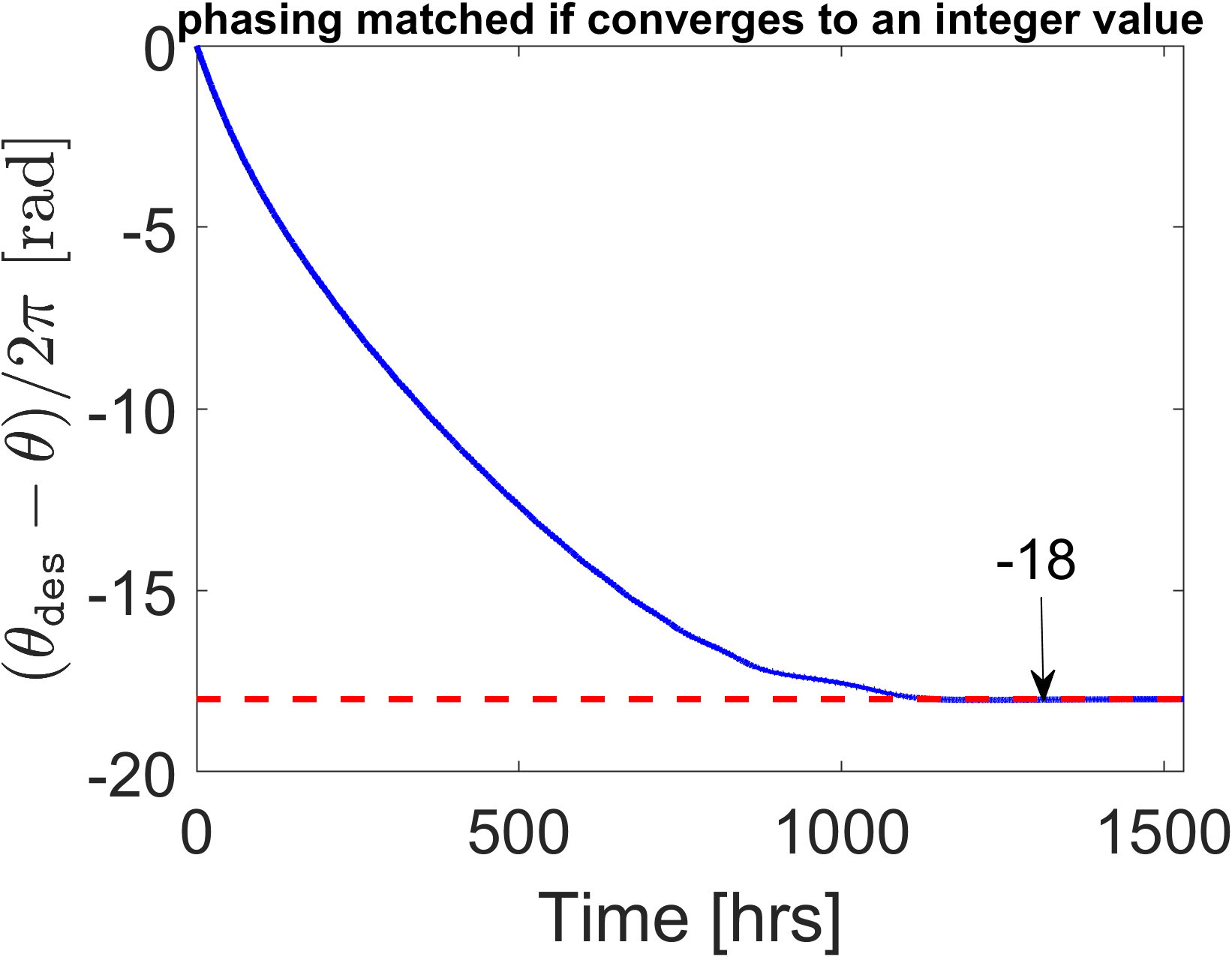}
	\caption{Current (left) and true anomaly difference (right) during transition from a lower semi-major axis orbit to a higher semi-major axis orbit with Lorentz force actuation.}
	\label{fig:currentlorentzup}
\end{figure}

Figures~\ref{fig:alorentzup}-\ref{fig:currentlorentzup} illustrate the maneuver to transition from the lower semi-major axis orbit to the higher semi-major axis orbit using Lorentz force-based actuation.  The maneuver is successfully completed with all state and control constraints satisfied.  Figure~\ref{fig:currentlorentzup}-right shows the time history of $\frac{\theta_{\tt des}(t)-\theta(t)}{2 \pi}$.  As it converges to an integer value $-18$, this indicates the true anomaly matching.  This matching is achieved by altering the semi-major axis reference command indicated in Figure~\ref{fig:alorentzup}-left by the red dashed line.

\subsection{Orbital Maneuvers to Lower Orbit with Lorentz Force Propulsion}

\begin{figure}[h!]
	\centering
	\includegraphics[width=0.4\linewidth]{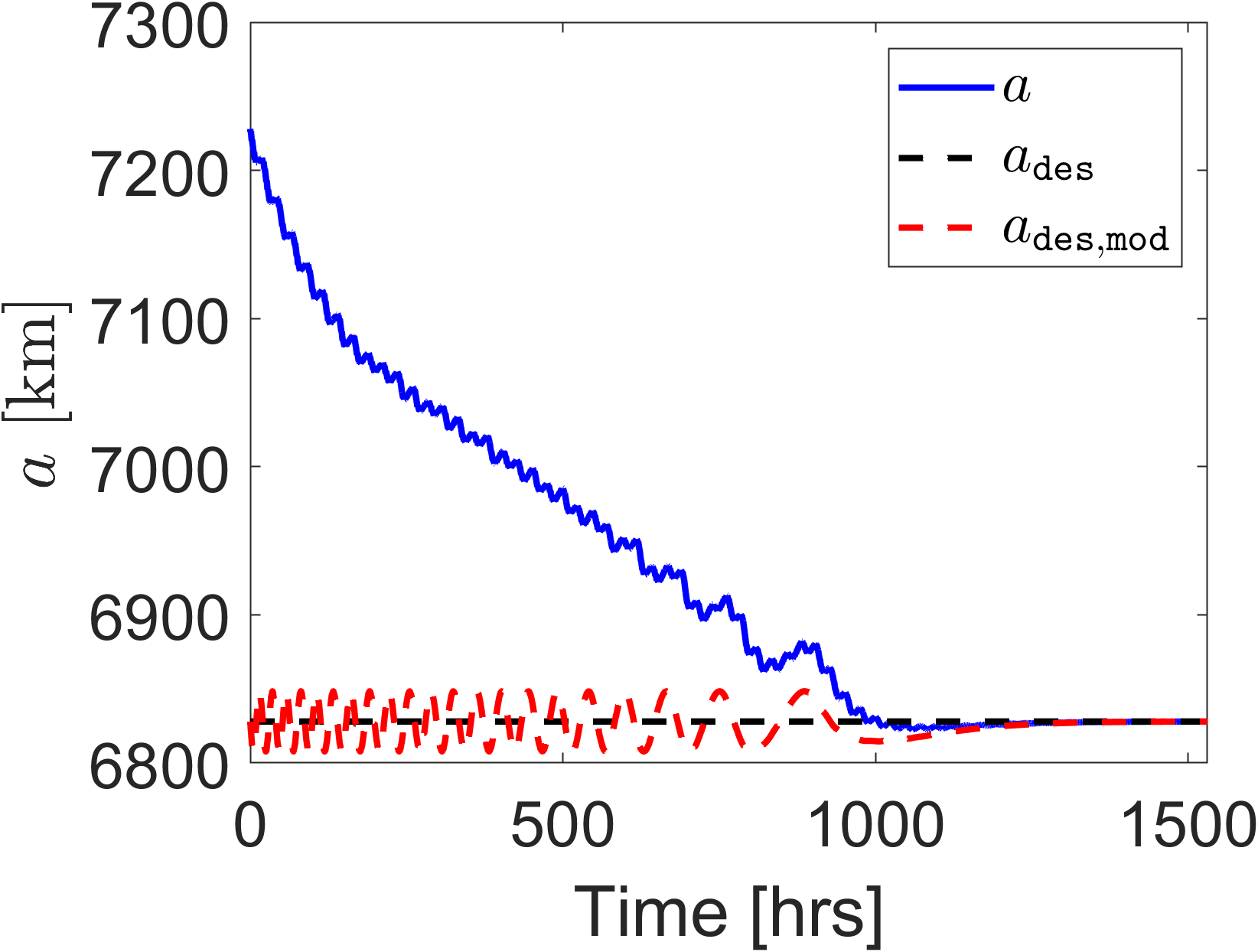}
	\includegraphics[width=0.4\linewidth]{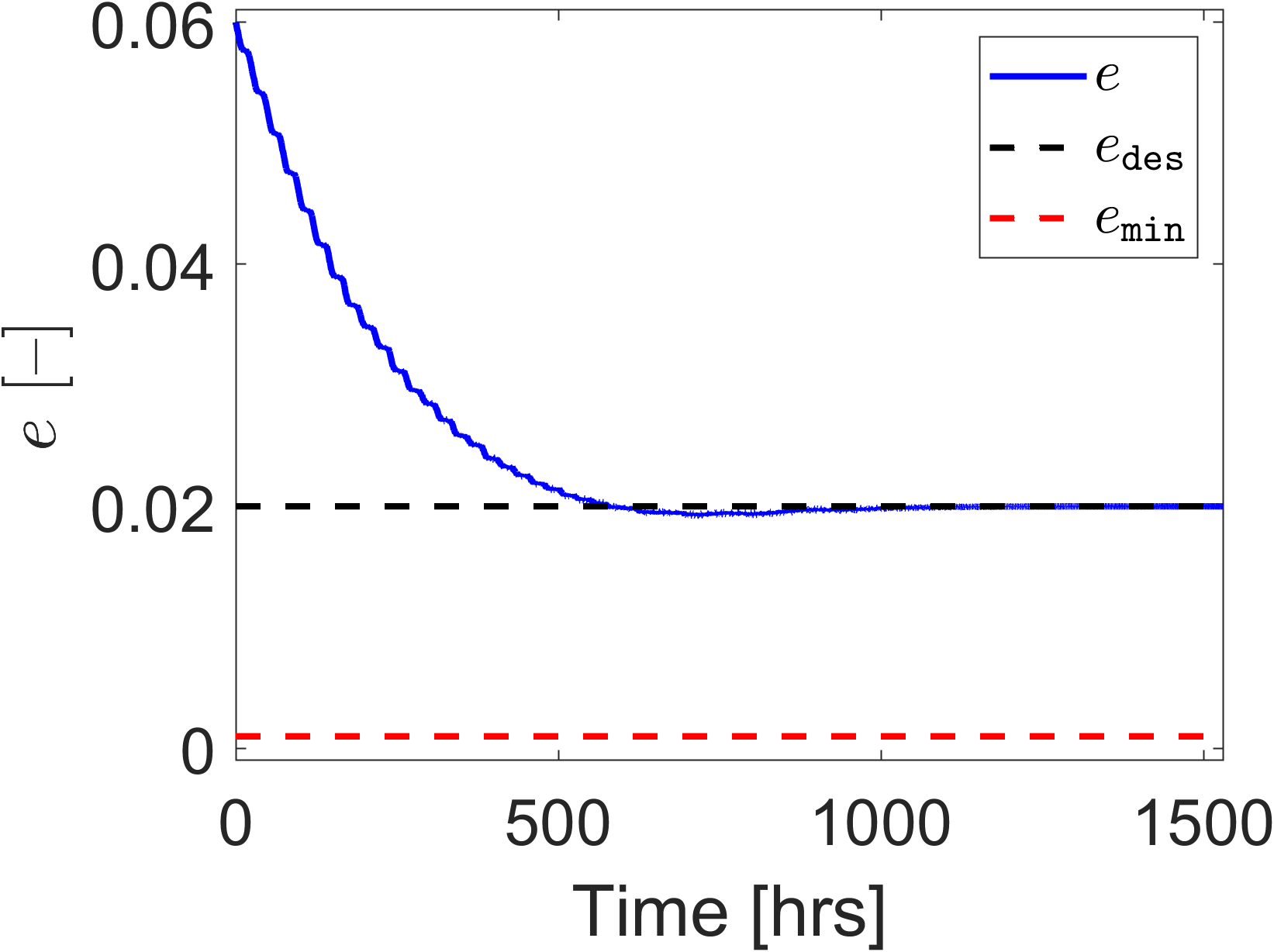}
	\caption{Semi-major axis (left) and eccentricity (right) during transition from a higher semi-major axis orbit to a lower semi-major axis orbit with Lorentz force actuation.}
	\label{fig:alorentzdown}
\end{figure}

\begin{figure}[h!]
	\centering
	\includegraphics[width=0.4\linewidth]{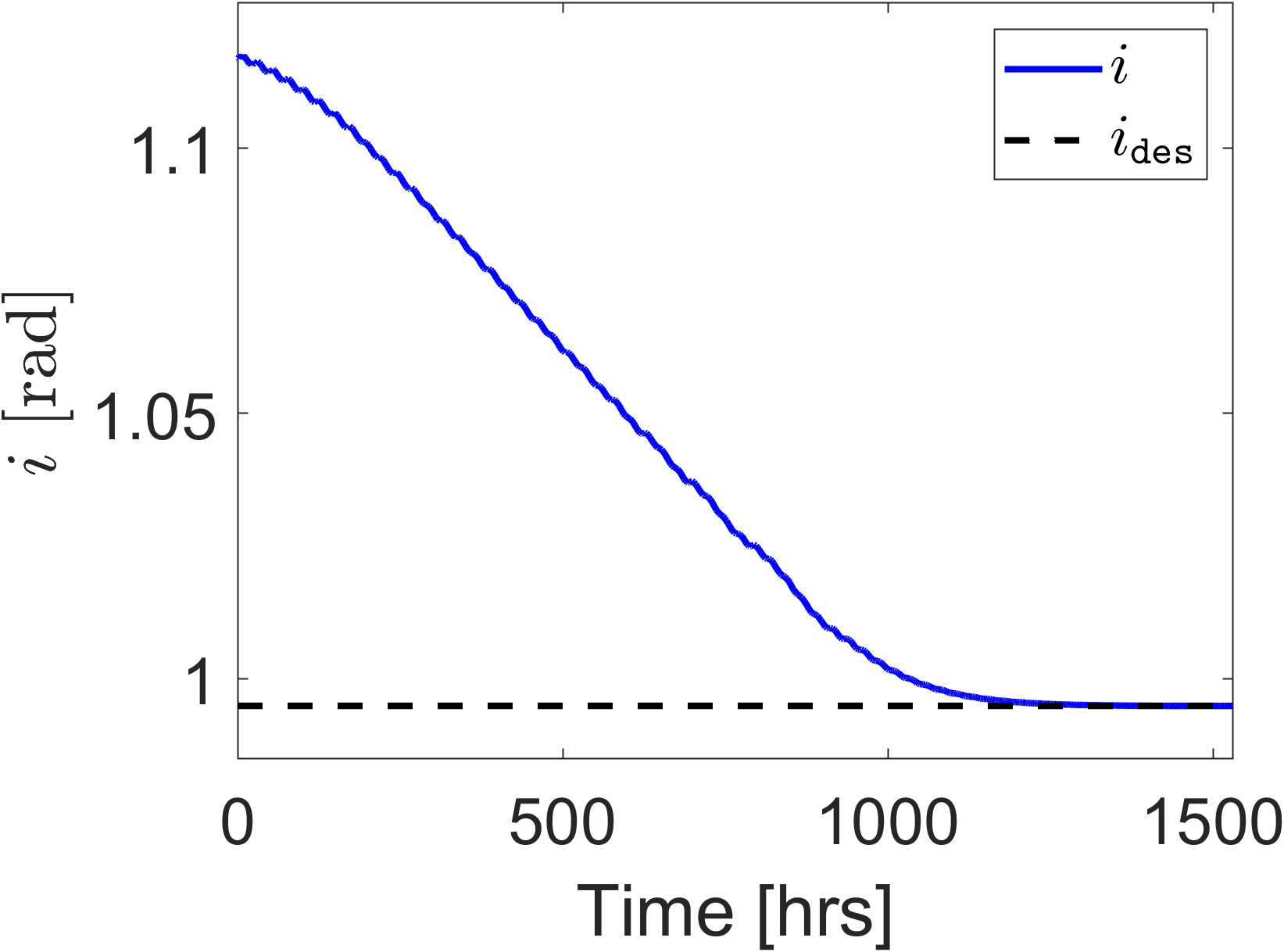}
    \includegraphics[width=0.4\linewidth]{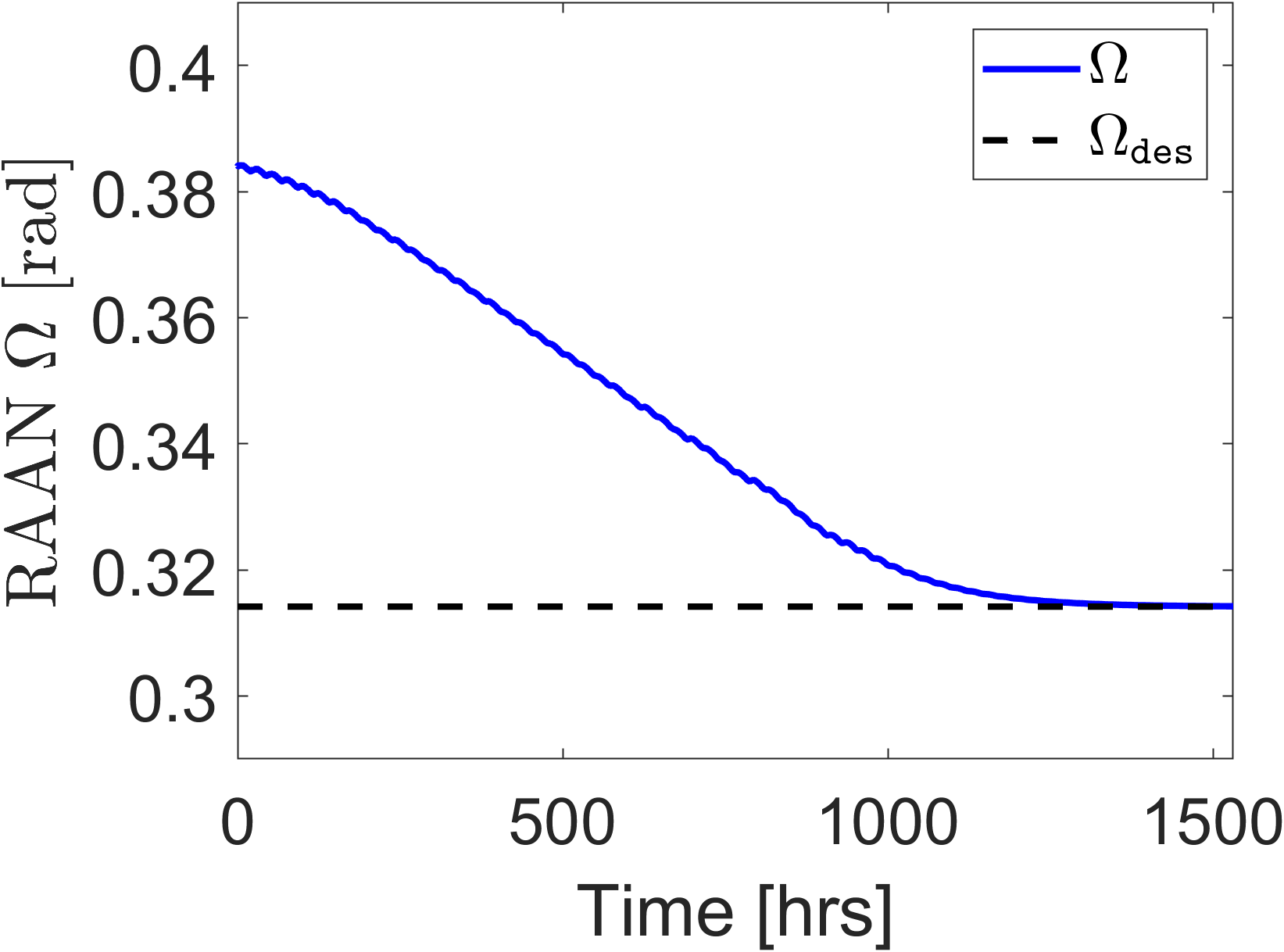}
	\caption{Inclination (left) and RAAN (right) during transition from a higher semi-major axis orbit to a lower semi-major axis orbit with Lorentz force actuation.}
	\label{fig:ilorentzdown}
\end{figure}

\begin{figure}[h!]
	\centering
	\includegraphics[width=0.4\linewidth]{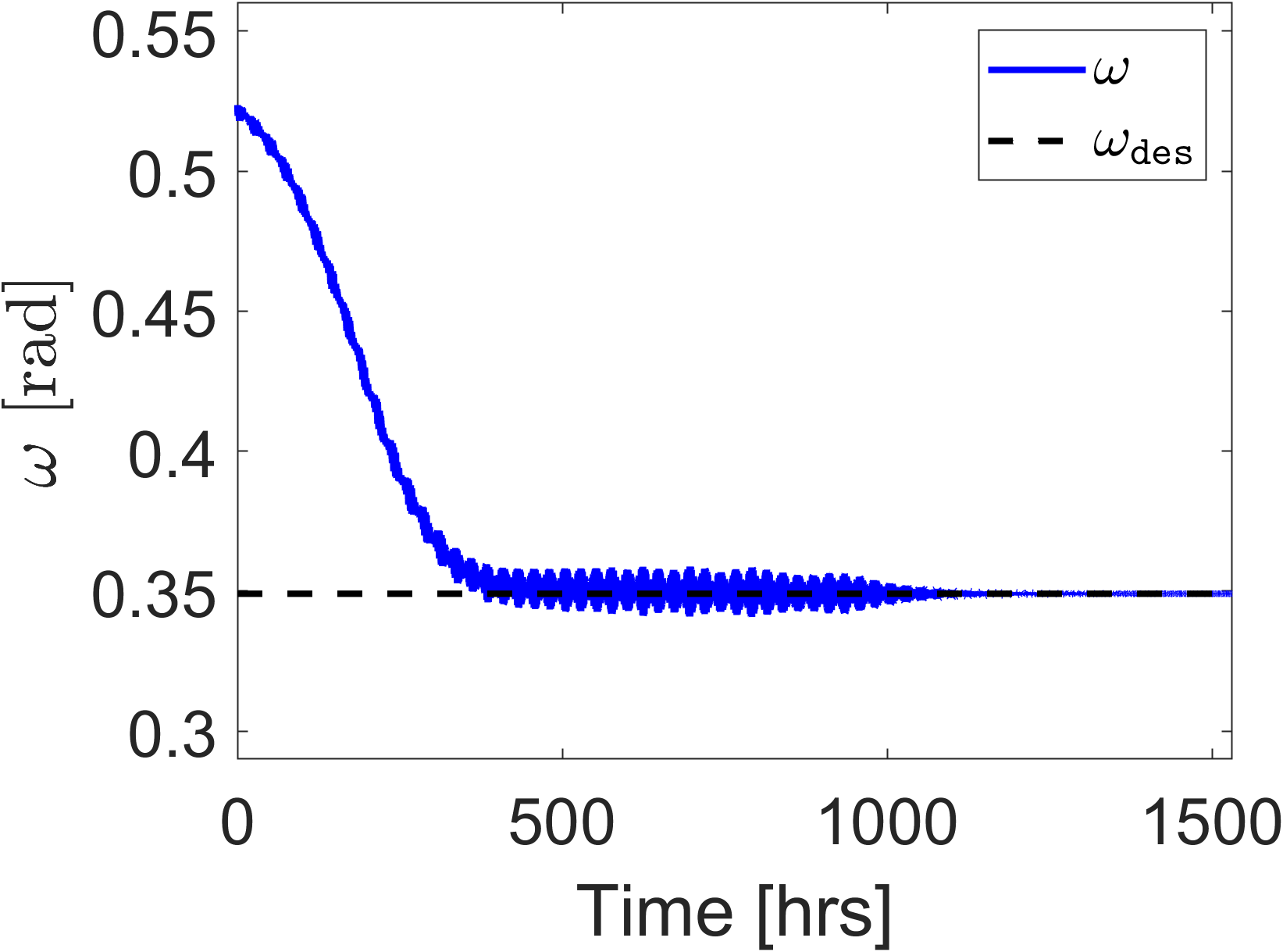}
	\includegraphics[width=0.4\linewidth]{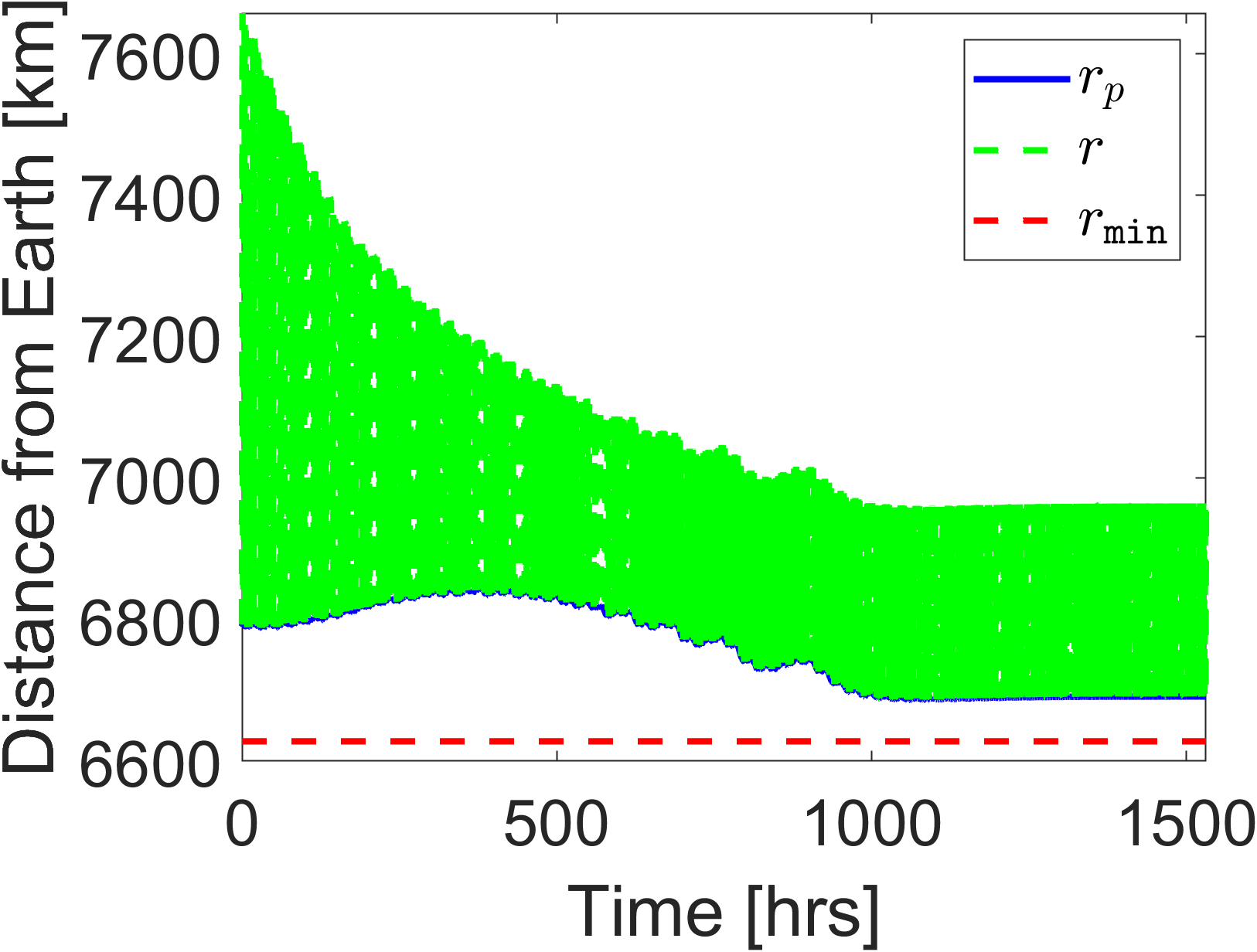}
	\caption{Argument of periapsis (left) and radius of periapsis (right) during transition from a higher semi-major axis orbit to a lower semi-major axis orbit with Lorentz force actuation.}
\label{fig:argumentofperiapsislorentzdown}
\end{figure}

\begin{figure}[h!]
	\centering
	\includegraphics[width=0.4\linewidth]{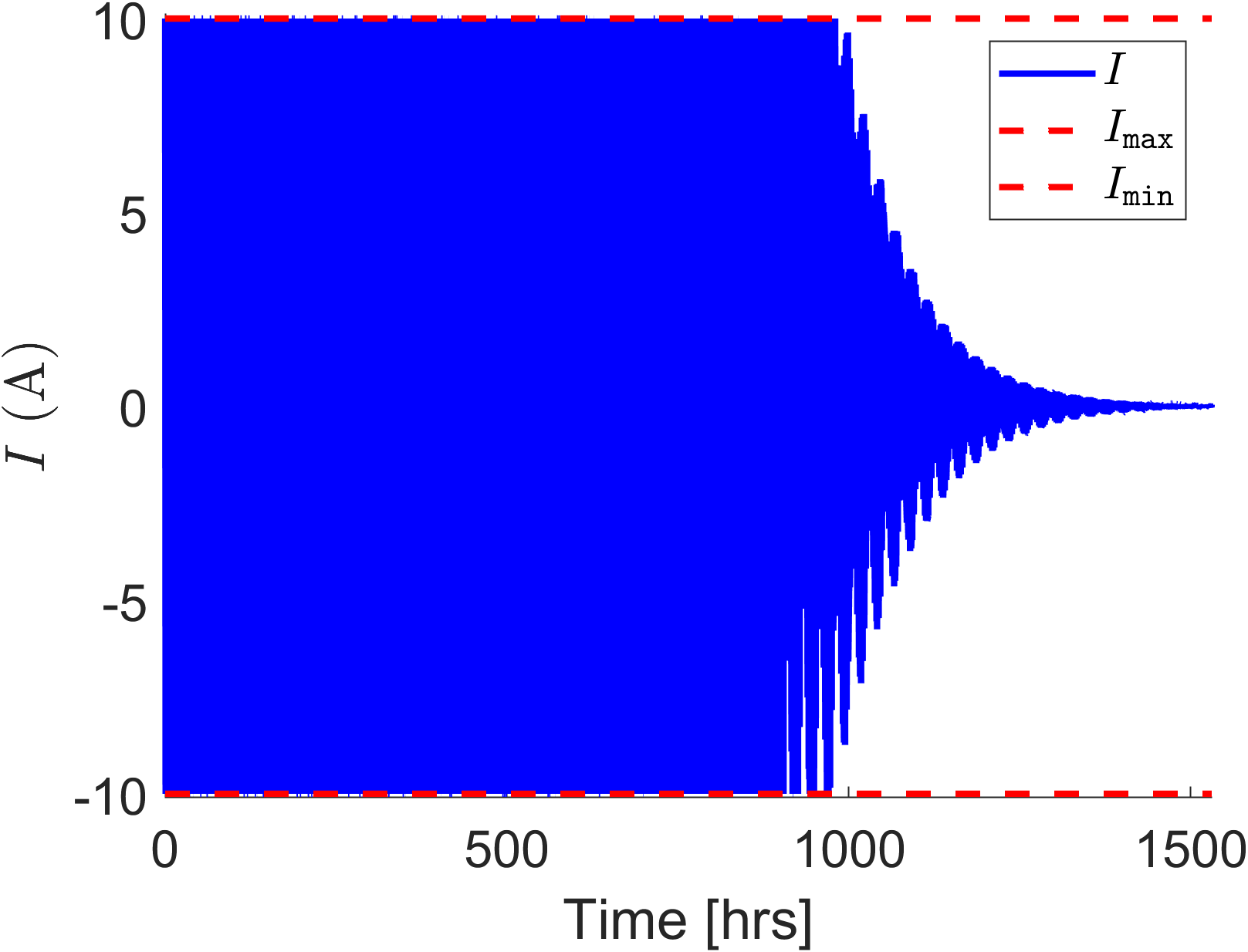}
	\includegraphics[width=0.4\linewidth]{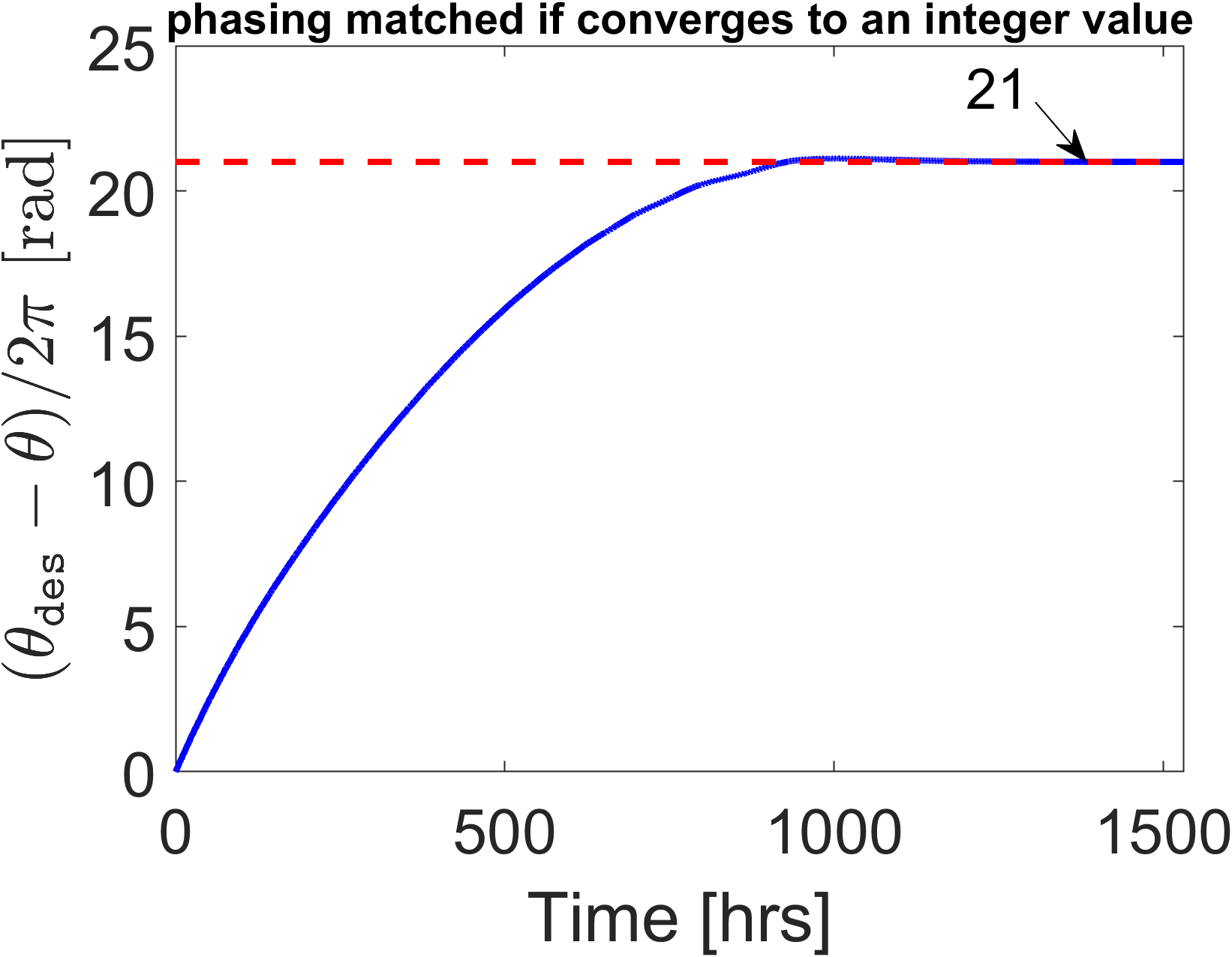}
	\caption{Current (left) and true anomaly difference (right) during transition from a higher semi-major axis orbit to a lower semi-major axis orbit with Lorentz force actuation.}
	\label{fig:currentlorentzdown}
\end{figure}

 Figures~\ref{fig:alorentzdown}-\ref{fig:currentlorentzdown} illustrate the maneuver to transition from the higher semi-major axis orbit to the lower semi-major axis orbit using Lorentz force-based actuation.  The maneuver is successfully completed with all state and control constraints satisfied.  Figure~\ref{fig:currentlorentzdown}-right shows the time history of $\frac{\theta_{\tt des}(t)-\theta(t)}{2 \pi}$.  As it converges to an integer value $21$, this indicates the true anomaly matching.  This matching is achieved by altering of the semi-major axis reference command indicated in Figure~\ref{fig:alorentzdown}-left by the red dashed line.

\section{Conclusion}
With the augmentation of an outer loop feedback controller for alterning the semi-major axis to a Lyapunov feedback law, rendezvous maneuvers can be performed during which the target semi-major axis, eccentricity, inclination, RAAN and argument of the periapsis and (time-dependent) true anomaly are achieved.  The control law can also handle state and control constraints through the augmentation of the barrier function terms to the Lyapunov function, saturation and the use of the reference governor to ensure convergence.  While convergence has been observed in a variety of simulations, rigorous conditions for  closed-loop stability for our cascade controller remain the subject for future analysis.

\section*{Acknowledgments}
The first author acknowledges the support of the Air Force Office of Scientific Research Grant number 
FA9550-23-1-0678 and National Science Foundation Awards 2515359 and 2502857.

\bibliography{sample}

@article{bombardelli2013deorbiting,
  title={Deorbiting performance of bare electrodynamic tethers in inclined orbits},
  author={Bombardelli, Claudio and Zanutto, Denis and Lorenzini, Enrico},
  journal={Journal of Guidance, Control, and Dynamics},
  volume={36},
  number={5},
  pages={1550--1556},
  year={2013},
  publisher={American Institute of Aeronautics and Astronautics}
}

@techreport{cosmo1997tethers,
  title={Tethers in Space Handbook},
  author={Cosmo, Mario L and Lorenzini, Enrico C},
  year={1997},
  url={https://ntrs.nasa.gov/citations/19980018321}
}

@book{levin2007dynamic,
  title={Dynamic Analysis of Space Tether Missions},
  author={Levin, E.},
  volume={126},
  year={2007},
  publisher={Univelt Incorporated}
}

@article{kolmanovsky2002speed,
	title={Speed-gradient approach to torque and air-to-fuel ratio control in DISC engines},
	author={Kolmanovsky, Ilya  and Druzhinina, Maria and Sun, Jing},
	journal={IEEE {T}ransactions on {C}ontrol {S}ystems {T}echnology},
	volume={10},
	number={5},
	pages={671--678},
	year={2002},
	publisher={IEEE}
}

@article{kolmanovsky1996switched,
  title={Switched mode feedback control laws for nonholonomic systems in extended power form},
  author={Kolmanovsky, Ilya and Reyhanoglu, Mahmut and McClamroch, N Harris},
  journal={Systems \& Control Letters},
  volume={27},
  number={1},
  pages={29--36},
  year={1996},
  publisher={Elsevier}
}

@book{gurfil2016celestial,
	title={Celestial Mechanics and Astrodynamics: {T}heory and Practice},
	author={Gurfil, Pini and Seidelmann, P Kenneth},
	volume={436},
	year={2016},
	publisher={Springer}
}

@inproceedings{marley2020kinematic,
  title={A kinematic hybrid feedback controller on the unit circle suitable for orientation control of ships},
  author={Marley, Mathias and Skjetne, Roger and Teel, Andrew R},
  booktitle={Proceedings of 59th IEEE Conference on Decision and Control (CDC)},
  pages={1523--1529},
  year={2020}
}

@inproceedings{dongare2022integrated,
  title={Integrated guidance and control of driftless control-affine systems with control constraints and state exclusion zones},
  author={Dongare, Abhijit and Sanyal, Amit K and Kolmanovsky, Ilya and Viswanathan, Sasi Prabhakaran},
  booktitle={Proceedings of American Control Conference (ACC)},
  pages={3893--3898},
  year={2022}
}

@article{garone2024constrained,
	title={On constrained Lyapunov stabilization of spacecraft orbits using Gauss Variational Equations},
	author={Garone, Emanuele and Kolmanovsky, Ilya},
	journal={Journal of Spacecraft and Rockets},
	volume={61},
	number={4},
	pages={1128--1132},
	year={2024},
	publisher={American Institute of Aeronautics and Astronautics}
}

@inproceedings{semeraro2023constrained,
  author    = {Semeraro, Simone and Kolmanovsky, Ilya and Garone, Emanuele},
  title     = {On constrained feedback control of spacecraft orbital transfer maneuvers},
  booktitle = {Proceedings of the 33rd AAS/AIAA Space Flight Mechanics Meeting},
  year      = {2023},
  address   = {Austin, Texas},
  month     = {January 15--19},
  note      = {Paper AAS 23-292},
  url       = {https://arxiv.org/abs/2407.16389}
}

@article{semeraro2024reference,
	title={Reference governor for constrained spacecraft orbital transfers},
	author={Semeraro, Simone and Kolmanovsky, Ilya and Garone, Emanuele},
	journal={Advanced Control for Applications: Engineering and Industrial Systems},
	volume={6},
	number={1},
	pages={e179},
	year={2024},
	publisher={Wiley Online Library}
}

@inproceedings{Kolmanovsky2024CommandGovernors,
  author    = {Ilya Kolmanovsky and Torbj{\o}rn Cunis and Emanuele Garone},
  title     = {Command governors based on bilevel optimization for constrained spacecraft orbital transfer},
  booktitle = {AIAA SciTech 2024 Forum},
  year      = {2024},
  publisher = {American Institute of Aeronautics and Astronautics},
  doi       = {10.2514/6.2024-0097}
}

@inproceedings{Lantukh2017EnhancedQLaw,
  author    = {Demyan V. Lantukh and Christopher L. Ranieri and Marc D. DiPrinzio and Peter J. Edelman},
  title     = {Enhanced Q-Law Lyapunov control for low-thrust transfer and rendezvous design},
  booktitle = {Proceedings of the AAS/AIAA Astrodynamics Specialist Conference, Paper AAS 17-589},
  year      = {2017},
  publisher = {Univelt, Inc.},
  volume    = {162}
}

@phdthesis{naasz2002classical,
  title={Classical element feedback control for spacecraft orbital maneuvers},
  author={Naasz, Bo James},
  year={2002},
  school={Virginia Tech}
}

@inproceedings{Petropoulos2004QLaw,
  author    = {A.E. Petropoulos},
  title     = {Low-thrust orbit transfers using candidate Lyapunov functions with a mechanism for coasting},
  booktitle = {AIAA/AAS Astrodynamics Specialist Conference},
  year      = {2004},
  pages={5089},
  number    = {AIAA 2004-5089},
  url       = {https://arc.aiaa.org/doi/10.2514/6.2004-5089}
}

@techreport{Petropoulos2005Refined,
  author      = {A.E. Petropoulos},
  title       = {Refinements to the Q-law for low-thrust orbit transfers},
  institution = {NASA Jet Propulsion Laboratory},
  year        = {2005},
  volume      = {120}, 
  number      = {JPL IOM 343R-05-02},
  url         = {https://ntrs.nasa.gov/citations/20210001725},
publisher={Pasadena, CA : Jet Propulsion Laboratory, National Aeronautics and Space Administration}
}

@article{garone2017reference,
	title={Reference and command governors for systems with constraints: {A} survey on theory and applications},
	author={Garone, Emanuele and Di Cairano, Stefano and Kolmanovsky, Ilya},
	journal={Automatica},
	volume={75},
	pages={306--328},
	year={2017},
	publisher={Elsevier},
}

\end{document}